% --------------------- Reglages generaux -------------------------
% --------------------- que vous pouvez modifier ------------------
%---------------------- sans faire de degats ----------------------

\hsize=13.50cm
\vsize=18cm
\parindent=12pt   
\parskip=0pt
%\parindent=0
\pageno=1

%%%%%%%%%%%%%%%%%%%%%%%%%%%%%%%%%%%%%%%%%%%%%%%%%%%%%%%%%%%%%
%
%  Font Bardot
%

% -------------------- Debut des macros privees -------------------
\catcode`\@=11
% --------------------- Les fontes --------------------------------
% --------------------- Les fontes --------------------------------

\font\eightrm=cmr8         \font\eighti=cmmi8
\font\eightsy=cmsy8        \font\eightbf=cmbx8
\font\eighttt=cmtt8        \font\eightit=cmti8
\font\eightsl=cmsl8        \font\sixrm=cmr6
\font\sixi=cmmi6           \font\sixsy=cmsy6
\font\sixbf=cmbx6

% Fontes AMS

\font\eightgoth=eufm10 at 8pt      \font\eightbboard=msbm10 at 8 pt
      \font\sevenbboard=msbm7
\font\sixgoth=eufm7 at 6 pt        \font\fivegoth=eufm5

\font\eightcyr=wncyr10 at 8 pt      
      
\font\sixcyr=wncyr10 at 6 pt

% Nouvelles familles pour les maths

\newfam\gothfam           \newfam\bboardfam
\newfam\cyrfam
%%%%%%%%%%%%%%%%%%%%%%%%%%%%%%%%%%%%%%%%%%%%%%%%%%%%%%%%%%%%
\font\titre=cmbx10 at 16pt
\font\Auteur=cmcsc10 at 10pt
\font\soustitre=cmbx10 at 12pt

%%%%%%%%%%%%%%%%%%%%%%%%%%%%%%%%%%%%%%%%%%%%%%
\def\eightpoint{%
   \textfont0=\eightrm \scriptfont0=\sixrm \scriptscriptfont0=\fiverm
   \def\rm{\fam\z@\eightrm}%
   \textfont1=\eighti  \scriptfont1=\sixi  \scriptscriptfont1=\fivei
   \def\oldstyle{\fam\@ne\eighti}\let\old=\oldstyle
   \textfont2=\eightsy \scriptfont2=\sixsy \scriptscriptfont2=\fivesy
   \textfont\gothfam=\eightgoth \scriptfont\gothfam=\sixgoth
   \scriptscriptfont\gothfam=\fivegoth
   \def\goth{\fam\gothfam\eightgoth}%
   \textfont\cyrfam=\eightcyr \scriptfont\cyrfam=\sixcyr
   \scriptscriptfont\cyrfam=\sixcyr
   \def\cyr{\fam\cyrfam\eightcyr}%
   \textfont\bboardfam=\eightbboard \scriptfont\bboardfam=\sevenbboard
   \scriptscriptfont\bboardfam=\sevenbboard
   \def\bb{\fam\bboardfam}%
   \textfont\itfam=\eightit
   \def\it{\fam\itfam\eightit}%
   \textfont\slfam=\eightsl
   \def\sl{\fam\slfam\eightsl}%
   \textfont\bffam=\eightbf \scriptfont\bffam=\sixbf
   \scriptscriptfont\bffam=\fivebf
   \def\bf{\fam\bffam\eightbf}%
   \textfont\ttfam=\eighttt
   \def\tt{\fam\ttfam\eighttt}%
   \abovedisplayskip=9pt plus 3pt minus 9pt
   \belowdisplayskip=\abovedisplayskip
   \abovedisplayshortskip=0pt plus 3pt
   \belowdisplayshortskip=3pt plus 3pt
   \smallskipamount=2pt plus 1pt minus 1pt
   \medskipamount=4pt plus 2pt minus 1pt
   \bigskipamount=9pt plus 3pt minus 3pt
   \normalbaselineskip=9pt
   \setbox\strutbox=\hbox{\vrule height7pt depth2pt width0pt}%
   \let\bigf@nt=\eightrm     \let\smallf@nt=\sixrm
   \normalbaselines\rm}
%
%
%%%%%%%%%%%%%%%%%%%%%%%%%%%%%%%%%%%%%%%%%%%%%%%%%%%%%%%%%%%
\input amssym.tex
\input amssym.def
\let\dsp=\displaystyle
\def\BR{{\Bbb R}}
\def\BN{{\Bbb N}}
\let\Br=\BR\let\pa=\partial\def\ga{\gamma}
\let\le=\leq\let\ge=\geq\let\ne=\neq\let\pv=\enskip\let\to=\rightarrow
\def\eps{\varepsilon}\def\La{\Lambda}\def\la{\lambda}
\def\dfrac#1#2{{\displaystyle\strut #1 \over\displaystyle\mathstrut #2}}
\def\Sum{\displaystyle\sum}
\def\conj#1{\overline{#1}}\def\sconj#1{\underline{#1}}
\def\abs#1{\left\vert #1\right\vert}
\def\scal#1{\left\langle\,#1\,\right\rangle}
\def\Scal#1#2{\left\langle\,#1\,,\,#2\,\right\rangle}
\def\dxi{{\rm d}\xi}\def\dy{{{\rm d}y}}\def\ds{{{\rm
   d}s}}\def\drm{{{\rm d}}}\def\dx{{\rm d}x} 
\def\dtau{{{\rm d}\tau}}\def\th{{\rm th}}\def\loc{{\rm loc}}\def\det{{\rm det}}
\def\hc{\hfill\cr}\def\tvi{\vrule height 12pt depth 5pt width 0pt}
\def\moins{\,\setminus}
\def\Lim{\displaystyle\lim}\def\vide{\,\emptyset}
\def\Id{{\rm Id}}
\def\SGN{\mathop{\rm sgn}\nolimits}
\def\sgnt{\SGN t}
\def\Int{\displaystyle\int}
\def\Norme#1{\left\Vert #1\right\Vert}\def\Normes#1{\Big\Vert #1\Big\Vert}
\def\norme#1{\Vert\,#1\,\Vert}
\def\Cup{\dsp\mathop\bigcup}
\def\og{\leavevmode\raise.3ex\hbox{$\scriptscriptstyle\langle\!\langle$}}
\def\fg{\leavevmode\raise.3ex\hbox{$\scriptscriptstyle\,\rangle\!\rangle$}}
\def\sfrac#1#2{{\scriptstyle\mathstrut#1 \over\scriptstyle\mathstrut#2}}
\def\dt{{{\rm d}t}}
\def\Supp{\mathop{\rm supp}\nolimits}
\def\supp{\Supp}
\def\erom{{\rm e}}
\def\Im{\mathop{\rm Im}}\def\Re{\mathop{\rm Re}}
\def\hbk{\hfill\break\null}

\def\vpt{\medskip}
\font\rond=pzdr at 12pt
\def\rondun{\ \hbox{\rond\char"AC}}\def\rondtwo{\ \hbox{\rond\char"AD}}
\def\rondtrois{\ \hbox{\rond\char"AE}}\def\rondq{\ \hbox{\rond\char"AF}}

\def\cqfd{\hskip 6bp\raise -8pt \hbox{\vrule\vbox to6pt% modif : -5
{\hrule width 6pt\vfill\hrule}\vrule}\par
\vskip 2pt plus 2pt minus 1pt}
\def\Intd{\Int\!\!\!\!\!\Int}
\def\cc#1{\hfill\enskip#1\enskip\hfill}
%%%%%%%%%%%%%%%%%%%%%%%%%%%%%%%%%%%%%%%%%%%%%%%%%%%%%%%%%%
%
%  Def ajoutees
%
\def\CV{{\cal V}}
\def\CF{{\cal F}}
\def\CD{{\cal D}}
\def\vpd{\medskip}

\def\BC{{\Bbb C}}
\def\gcap{\Auteur}
\def\bun{{1\mkern -5mu \rm I}}
\def\bunt{\bun_{[0,T]}}
%
%
% Fin  Def ajoutees
%
%%%%%%%%%%%%%%%%%%%%%%%%%%%%%%%%%%%%%%%%%%%%%%%%%%%%%%%%%%%%%%%%%%%%
\magnification 1200
%%%%%%%%%%%%%%%%%%%%%%%%%%%%%%%%%%%%%%%%%%%%%%%%%%%%%%%%%%%%%%%%%%%%%%%%%%%
%
%   Debut texte
%
%%%%%%%%%%%%%%%%%%%%%%%%%%%%%%%%%%%%%%%%%%%%%%%%%%%%%%%%%%%%%%%%%%%%%%%%%%%
{\centerline {\titre The Kato smoothing effect for Schr{\"o}dinger}
\medskip
{\centerline {\titre
equations with 
unbounded potentials}
\medskip
{\centerline {\titre
in exterior domains}
\bigskip

{\centerline {by}

\vskip 0,5cm
{\centerline{\Auteur  Luc ROBBIANO%
{\eightpoint\footnote{$^1$}{Laboratoire 
de Math\'ematiques de Versailles,  Universit\'e de Versailles St
Quentin, CNRS, 45,
Avenue des Etats-Unis, 78035 Versailles, FRANCE,
  e-mail : Luc.Robbiano@math.uvsq.fr}}}

{\centerline{\Auteur  Claude ZUILY%
{\eightpoint\footnote{$^2$}{Univ. Paris-Sud, UMR 8628,
Orsay, F-91405 , e-mail :
Claude.Zuily@math.u-psud.fr}}{\eightpoint\footnote{$^3$}{CNRS, Orsay,
F-91405}}}}

\vskip 0,3cm

{\centerline {\hbox to 3cm{\hrulefill}}}

\vskip 20pt

\baselineskip=14pt

\noindent {\soustitre 1. Introduction}
\vpt
The Kato $\dfrac 12-$ smoothing effect for Schr{\"o}dinger equations has
received much attention during the last years. See
Constantin-Saut [C-S], Sj{\"o}lin [Sj], {Vega} [V],
 {Yajima} [Y] for the case of the flat Laplacian in $\BR^d$.
It has been successively extended to variable coefficients  operators
by {Doi} (see [D1], [D2]) and to perturbations of such operators by 
potentials  growing at most quadratically at infinity (see {Doi}
[D3]). The aim of this paper is to consider exterior boundary value
problems for variable coefficients operators with unbounded potentials.
The case of potentials decaying at infinity has been considered by
 {Burq} [B1] using resolvent estimates.

Our main smoothing estimate is proved by contradiction. The idea of
proving estimates by contradiction (with the appropriate technology)
goes back to   {Lebeau} [L] and it has been subsequently used with
success by several authors, (see e.g.   {Burq} [B2]).

In this paper, some ideas of   {G{\'e}rard-Leichtnam} [G-L],   {Burq}
[B3] and   {Miller} [Mi] will be also used.

Let us briefly outline how this method applies here. Assuming that our
estimate is false gives rise, after renormalization, to a sequence
$\big(u_k\big)$ which is bounded in
$L^2_{loc}\big([0,T]\times\BR^d\big)$. To a subsequence we associate a
microlocal semi-classical defect measure $\mu$ in the sense of   
G{\'e}rard [G]. Then, roughly speaking, there are three main steps in the
proof. First $\mu$ does not vanish identically. Moreover $\mu$ vanishes
somewhere (in the incoming region). Finally the support of $\mu$ is
invariant by the generalized bicharacteristic flow (in the sense of
  {Melrose-Sj{\"o}strand} [M-S]). Since one of our assumptions (the non
trapping condition) ensures that the backward generalized flow always
meet the incoming region (where $\mu$ vanishes) we obtain a
contradiction thus proving the desired estimate.

Let us now describe more precisely the content of each section.

In the next one we describe the assumptions and state the main result of
this paper. In the third section we begin our contradiction argument and
we show in the next one how to obtain a bounded sequence in
$L^2\big([0,T], L^2_{loc}(\Br^d)\big)$. Then in the fifth section we
introduce the semi classical defect measure $\mu$ and we state without
proof the invariance of its support by the generalized Melrose-Sj{\"o}strand
bicharacteristic flow.  In the next section we show that $\mu$ does not
vanish identically while in section seven we show that $\mu$ vanishes in
the incoming region. In the section eight we end the proof of our main
result by achieving a contradiction. Finally in the appendix (section
nine) we recall the geometrical framework we prove the invariance of the
support of $\mu$ and we end by proving some technical Lemmas used in the
preceeding sections.

Aknowledgments~: The authors would like to thank Nicolas Burq for useful
discussions at an early stage of this work.

\vpt
\noindent {\soustitre 2. Statement of the result}
\vpt
Let $K$ be a compact obstacle in $\BR^d$ whose complement $\Omega$ is a
connected open set with ${\cal C}^\infty$ boundary $\pa\Omega$.

Let $P$ be a second order differential operator of the form
$$P=\Sum_{j,k=1}^d D_j\big(a^{jk}(x)D_k\big)+V(x)\quad,\quad D_j=\dfrac
1i\dfrac\pa{\pa x_j}\leqno{(2.1)}$$
whose coefficients $a^{jk}$ and $V$ are assumed (for simplicity) to be in
${\cal C}^\infty(\conj\Omega)$, real valued and $a^{jk}=a^{kj}$,
$1\le j,k\le d$

We shall set
$$p(x,\xi)=\Sum_{j,k=1}^d a^{jk}(x)\xi_j\xi_k\leqno{(2.2)}$$
and we shall assume that
$$\exists\, c>0~:~p(x,\xi)\ge c\abs\xi^2,\hbox{ for $x$ in $\conj\Omega$
and $\xi$ in $\BR^d$}.\leqno{(2.3)}$$
To express the remaining assumptions on the coefficients we introduce
the metric
$$g=\dfrac{\dx^2}{\scal x^2}+\dfrac{\dxi^2}{\scal\xi^2}\leqno{(2.4)}$$
where $\scal\cdot=\big(1+\abs\cdot^2\big)^{1/2}$ and we shall denote by
$S_\Omega(M,g)$ the H{\"o}rmander's class of symbols if $M$ is a weight.
Then $a\in S_\Omega(M,g)$ iff $a\in{\cal
C}^\infty\big(\conj\Omega\times\Br^d\big)$ and for all $\alpha$, $\beta$
in $\BN^d$ one can find $C_{\alpha,\beta}>0$ such that
$$\abs{D_x^\beta D_\xi^\alpha a(x,\xi)}\le C_{\alpha\beta}M(x,\xi)\scal
x^{-\abs\beta}\scal\xi^{-\abs\alpha}$$
for all $x$ in $\conj\Omega$ and $\xi$ in $\BR^d$.

Next we assume
$$\left\{\eqalign{
(i)\pv &a^{jk}\in S_\Omega(1,g),\pv \nabla_x a^{jk}(x)=o\Big(\dfrac
1{\abs x}\Big),\pv \abs x\to +\infty,\pv 1\le j,k\le d.\hc
(ii)\pv&V\in S_\Omega\big(\scal x^2,g\big),\pv V\ge-C_0\hbox{ for some
positive constant } C_0.\hc}\right.\leqno{(2.5)}$$
Under the assumptions (2.3), (2.5) the operator $P$ is essentially self
adjoint on $\big\{u\in{\cal C}_0^\infty(\conj\Omega)~:
u_{|\pa\Omega}=0\big\}$. We shall denote by $P_D$ ($D$ means Dirichlet)
its self adjoint extension.

Let us describe now our geometrical assumptions. We shall assume
$$\left\{\eqalign{
&\hbox{the generalized bicharacteristic flow (in the sense of}\hc
&\hbox{Melrose-Sj{\"o}strand) is not backward
trapped.}\hc}\right.\leqno{(2.6)}$$
This assumption needs some explanations. Let $M=\Omega\times\BR_t$. Let
us set  $T^*_bM=T^*M\moins\{0\}\cup T^*\pa M\moins\{0\}$. We have a
natural restriction map $\pi~:~T^*\BR^{d+1}_{|\conj M}\moins\{0\}\to
T^*_b M$ which is the identity on $T^*\Br^{d+1}_{|M}\moins\{0\}$.

Let $\Sigma=\big\{(x,t,\xi,\tau)\in T^*\BR^{d+1}\moins\{0\}$~:
$x\in\conj\Omega$, $t\in\big[0,T \big]$, $\tau+p(x,\xi)=0\big\}$ and
$\Sigma_b=\pi(\Sigma)$. For $a\in\Sigma_b$ the generalized
bicharacteristic $\Gamma(t,a)$ lives in $\Sigma_b$ (see section 9.1 for
details). Then (2.6) means the following.

For any $a$ in $\Sigma_b$ there exists $s_0$ such that for all $s\le s_0$
we have $\Gamma(s,a)\subset T^*M\moins\{0\}$, then 
$\Gamma(s,a)=\big(x(s),t,\xi(s),\tau\big)$ where $\big(x(s),\xi(s)\big)$
is the usual flow of $p$ and $\Lim_{s\to -\infty} \abs{x(s)}=+\infty$.

We shall need another assumption on the flow whose precise meaning will
be given in the appendix, section 9.1, Definition 9.3.
$$\left\{\eqalign{
&\hbox{The  bicharacteristics  have no contact}\hc
&\hbox{of infinite order with the boundary
$\pa\Omega$}\hc}\right.\leqno{(2.7)}$$

Now we set
$$\Lambda_D=\Big(\big(1+C_0\big)\Id+P_D\Big)^{1/2}\leqno{(2.8)}$$
which is well defined by the functionnal calculus of self adjoint
positive operators.

We shall consider the problem
$$\left\{\eqalign{
&i\,\dfrac{\pa u}{\pa t}+P_Du=0,\hc
&u_{|t=0}=u_0,\hc
&u_{|\pa\Omega\times\BR_t}=0,\hc}\right.\leqno{(2.9)}$$
where $u_0\in L^2(\Omega)$.

Then we can state our main result.

\proclaim Theorem 2.1. 
Let $T>0$, $\chi   \in{\cal C}_0^\infty(\conj\Omega)$,
$s\in\big[-1,1\big]$. Let $P$ be defined by (2.1) satisfying the
assumptions (2.3), (2.5), (2.6) and (2.7). Then one can  find a positive
constant $C(T,\chi   , s)=C$ such that
$$\Int_0^T\Big\|\chi   \Lambda^{s+\sfrac 12}_Du(t)\Big\|^2_{L^2(\Omega)}\dt\le
C\Norme{\Lambda^s_Du_0}^2_{L^2(\Omega)}\leqno{(2.10)}$$
for all $u_0$ in ${\cal C}_0^{\infty}(\Omega)$, where $u$ denotes the
solution of (2.9).

 Here are some remarks

\proclaim Remarks 2.2. \hfill\break
{$(i)$}{Theorem 2.1 can be extended to operators of the form
$$P=\Sum_{j,k=1}^d
\left(D_j-b_j(x)\right)a^{jk}(x)\left(D_k-b_k(x)\right)+V(x)$$
where $b_j\in S_\Omega\big(\scal x,g\big)$.}\hfill\break
{$(ii)$}{In the case $\Omega=\BR^d$ the above result has been
proved by   {Doi} [D3].}\hfill\break
{$(iii)$}{Without lack of generality one may assume $s=0$ in the
theorem.\hbk
Moreover working with $\widetilde u(t)=\erom^{-i(1+C_0)t}u(t)$
one may assume $V\ge 1$ in (2.5)$\,(ii)$ and 
$\Lambda_D=P_{D}^{1/2}$
which we will assume in that follows.}

\vpt
\noindent {\soustitre 3. The contradiction argument}
\vpt
Our goal is to begin the proof by contradiction of Theorem 2.1. We shall
first consider a version of the estimate which is localized in
frequency.

Let $T>0$ and $I=\big]0,T\big[$. Let $\theta\in{\cal C}_0^\infty(\Br)$
be such that $\Supp\theta\subset \Big\{t~:~\dfrac 12\le\abs t\le
2\Big\}$.

\proclaim Theorem 3.1.
Let $\chi   _0\in{\cal C}_0^\infty\big(\BR^d\big)$ be fixed. There exists
$C>0$, $h_0>0$ such that for all $h$ in $\big]0,h_0\big[$ we have
$$\Int_0^T\Norme{\chi
_0\theta\big(h^2P_D\big)P_D^{1/4}u(t)}^2_{L^2(\Omega)}\le 
C\Norme{u_0}^2_{L^2(\Omega)}\leqno{(3.1)}$$
for all $u_0\in L^2(\Omega)$.

Here $\theta\big(h^2P_D\big)$ is defined by the functionnal calculus of
selfadjoint operators.

Recall that $K$ is our compact obstacle. We take $R_0\ge 1$ so large
that
$$K\subset\big\{x\in\BR^n~:~\abs x<R_0\big\}.\leqno{(3.2)}$$
Let $R_1>R_0$ be such that
$\Supp\chi   _0\subset\big\{x\in\Br^d~:~\abs x<R_1\big\}$. Let
$\chi   _1\in{\cal C}_0^\infty\big(\BR^d\big)$ be such that $0\le \chi   _1\le
1$ and
$$\left\{\eqalign{
&\chi   _1(x)=1\pv\hbox{if}\pv \abs x\le R_1+2\hc
&\Supp\chi   _1\subset\big\{x~:~\abs x\le
R_1+3\big\}.\hc}\right.\leqno{(3.3)}$$
Then $\chi   _0\chi   _1=\chi   _0$. Moreover let us set
$$\theta_1(t)=t^{\sfrac 14}\theta(t),\pv \theta_2(t)=t^{-\sfrac
14}\theta(t).\leqno{(3.4)}$$
It is easy to see that (3.1) will be implied by the following estimate.
$$\left\{\eqalign{
&\exists C>0, \exists h_0>0~:~ \forall h\in\big]0,h_0\big[, \forall
u_0\in{\cal C}^\infty_0(\Omega),\hc
&\Int_0^T\Big\|\chi   _1h^{-\sfrac
12}\theta_1\big(h^2P_D\big)u(t)\Big\|^2_{L^2(\Omega)}\dt\le
C\Norme{u_0}^2_{L^2}.\hc}\right.\leqno{(3.5)}$$
We shall prove (3.5) by contradiction. Assuming it is false, taking
$h_0=\dfrac 1k$, $C=k$, $k\in\BN^*$, we deduce sequences
$\big(h_k\big)\to 0$, $u^0_k\in{\cal C}^\infty_0(\Omega)$, such that
$$\Int_0^T\Big\|\chi   _1h_k^{-\sfrac
12}\theta_1\big(h^2_kP_D\big)u_k(t)\Big\|^2_{L^2(\Omega)}\dt>k
\Norme{u^0_k}^2_{L^2}.$$  
It follows that the left hand side does not vanish. Therefore if we set
$$\left\{\eqalign{
\alpha^2_k&=\Int_0^T\Normes{\chi    h_k^{-\sfrac
12}\theta_1\big(h^2_kP_D\big)u_k(t)}^2_{L^2(\Omega)}\dt>0,\hc
\widetilde u\,^0_k&=\dfrac 1{\alpha_k} u^0_k,\pv\widetilde u_k=\dfrac
1{\alpha_k} u_k,\hc
w_k&=h^{-\sfrac 12}_k\theta_1\big(h^2_kP_D\big)\widetilde
u_k,\hc}\right.\leqno{(3.6)}$$
we see that
$$\left\{\eqalign{
(i)\quad &\Int_0^T\Norme{\chi   _1w_k(t)}^2_{L^2(\Omega)}\dt=1,\hc
(ii)\quad &\Norme{\widetilde u\,^0_k}_{L^2(\Omega)}<\dfrac
1k\cdotp\hc}\right.\leqno{(3.7)}$$

\vpt
\noindent {\soustitre 4. The sequence $\big(w_k\big)$ is bounded in
$
L^2\big(\BR ,L^2_{loc}(\BR^d)\big)$}
\vpt
We shall prove in this section the following result.

\proclaim Proposition 4.1.
    For any $\chi   \in {\cal C}^\infty_0\big(\Br^d\big)$ one can find a
positive constant $C$ such that
 $$\Int_0^T\Norme{\chi    w_k (t)}^2_{L^2(\Omega)}\dt\,\le\,C\leqno{(4.1)}$$
 for all $k\ge 1$.
 
\vpt
\noindent {\bf Proof}
\vpt
We begin by extending to the whole $\BR^d$ the operator $P$ given in
(2.1).

Let $\chi   _2\in{\cal C}^\infty_0\big(\Br^d\big)$ be such that
$0\le\chi   _2\le 1$ and
$$\chi   _2(x)=1\hbox{ if } \abs x\le R_0,\pv \chi   _2(x)=0\hbox{ if } \abs
x\ge R_0+1.\leqno{(4.2)}$$
Then we set for $x\in\Br^d$,
$$\widetilde P=\Sum_{j,k=1}^d
D_j\big(\chi   _2\delta_{jk}D_k\big)+\Sum_{j,k=1}^d
D_j\left(\big(1-\chi   _2\big)a^{jk}(x)
D_k\right)+\chi   _2+\big(1-\chi   _2\big)V\leqno{(4.3)}$$
where $\delta_{jk}$ denotes the Kronecker symbol.

The principal symbol of $\widetilde P$ is
$$\left\{\eqalign{\widetilde p(x,\xi)&=\Sum_{j,k=1}^d \widetilde
a^{jk}(x)\xi_j\xi_k\hc
\widetilde
a^{jk}(x)&=\chi_2(x)\delta_{jk}+\big(1-\chi_2(x)\big)a^{jk}(x)
\hc}\right.\leqno{(4.4)}$$ 
According to conditions (2.2), (2.5), (2.6) we have the following,
$$\left\{\eqalign{
(i)\quad &\widetilde P=P\hbox{ if } \abs x\ge R_0+1,\hc
(ii)\quad &\widetilde p(x,\xi)\ge\widetilde c\abs\xi^2, x\in\BR^d,
\xi\in\BR^d, \widetilde c>0,\hc
(iii)\quad &\widetilde a^{jk}\in S_{\BR^d}(1,g)=S(1,g),
\nabla_x\widetilde a^{jk}(x)=o\big(\abs x^{-1}\big), \abs x\to
+\infty,\hc
(iv)\quad &\widetilde V=\chi   _2+\big(1-\chi   _2\big)V\in S\big(\scal
x^2,g\big)\hbox{ and } \widetilde V\ge 1,\hc
(v)\quad &\hbox{The flow of $\widetilde p$ in non
trapping.}\hc}\right.\leqno{(4.5)}$$
Let $\chi   _3\in{\cal C}^\infty_0\big(\Br^d\big)$ be such that
$0\le\chi   _3\le 1$ and with $R_1$ defined in (3.3),
$$\chi   _3(x)=1\hbox{ if } \abs x\le R_1+1,\pv \chi   _3(x)=0\hbox{ if }\abs
x\ge R_1+2.\leqno{(4.6)}$$
Then, according to (3.3), we have for $\alpha\ne 0$,
$$\Supp\pa^\alpha\chi   _3\subset\big\{x~:~R_1+1\le\abs x\le
R_1+2\big\}\subset\big\{x~:~\chi   _1(x)=1\big\}.\leqno{(4.7)}$$
Moreover since $R_1>R_0$ (see (3.2)) we have
$$\Supp\big(1-\chi   _3\big)\subset\big\{x~:~\abs
x>R_0+1\big\}.\leqno{(4.8)}$$
It follows from (4.5)$\,(i)$ that
$$P=\widetilde P\ \hbox{ on }\Supp\big(1-\chi   _3\big).\leqno{(4.9)}$$
Now with $w_k$ defined in (3.6) we set,
$$U_k=\big(1-\chi   _3\big)w_k.\leqno{(4.10)}$$
Then we have
$$\left\{\eqalign{
&\big(D_t-\widetilde P\big)U_k=G_k\hc
&G_k=\big[\widetilde P, \chi   _3\big]w_k\hc
&U_k(0)=\big(1-\chi   _3\big)h_k^{-\sfrac
12}\theta_1\big(h^2_kP_D\big)\widetilde u_k^0.\hc}\right.\leqno{(4.11)}$$
According to conditions $(ii)$ to $(v)$ in (4.5) we may apply Theorem
2.8 in [D3] with $s=-\dfrac 12$. It follows that for any $\chi   \in{\cal
C}_0^\infty\big(\BR^d\big)$ and any $\nu>0$ we have,
$$\Int_0^T\Norme{\chi
U_k(t)}^2_{L^2}\dt\,\le\,C_\nu\Big(\Big\|E_{-\sfrac 
12}U_k(0)\Big\|^2_{L^2}+\Int_0^T\Big\|\scal
x^{\sfrac{1+\nu}2}E_{-1}G_k(t)\Big\|^2_{L^2}\dt\Big)\leqno{(4.12)}$$
where $E_s$ is the pseudo-differential operator with symbol
$e_s(x,\xi)=\big(1+\widetilde p(x,\xi)+\abs x^2\big)^{s/2}$ which
belongs to $S\left(\big(\abs\xi+\scal x\big)^s,g\right)$.

To handle the first term in the right hand side of (4.12) we shall need
the following Lemma.

\proclaim Lemme 4.2.
    Let $Q=P_D^{\sfrac 12}\big(1-\chi   _3\big)A_{-1}$ where $A_{-1}\in{\cal
O}p S\left(\big(\abs\xi+\scal x\big)^{-1},g\right)$. Then $Q$ is bounded
from $L^2\big(\BR^d\big)$ to $L^2(\Omega)$.

\vpt
\noindent {\bf Proof}
\vpt
Let $\CV(\Omega)=\big\{u\in H^1_0(\Omega)~: V^{\sfrac 12}u\in
L^2(\Omega)\big\}$ endowed with the norm
$\norme u^2_{\CV(\Omega)}=\norme u^2_{H^1}+\norme{V^{\sfrac
12}u}^2_{L^2}$. It is well known that $\CV(\Omega)$ is the domain
of $P_D^{\sfrac 12}$ and that $\norme u_{\CV(\Omega)}$ is
equivalent to $\norme{P^{\sfrac 12}_Du}_{L^2(\Omega)}$. Moreover since
$\abs V^{\sfrac 12}\le C\scal x$ we have,
$$\Norme{\big(1-\chi   _3\big)f}_{\CV(\Omega)}\le C\left(\norme
f_{H_1(\BR^d)}+\Norme{\scal x f}_{L^2(\Br^d)}\right)$$
whenever the right hand side is finite.

It follows that we can write
$$\eqalign{\norme{Qu}_{L^2(\Omega)}
\le C_1\Norme{\big(1-\chi   _3\big)A_{-1}u}_{\CV(\Omega)}
&\le C_2\left(\Norme{A_{-1}u}_{H^1(\BR^d)}+\Norme{\scal x
A_{-1}u}_{L^2(\BR^d)}\right)\cr
&\le C_3\norme u_{L^2(\BR^d)}}$$
\hfill\cqfd

Now let us set $\rondun=\Big\|E_{-\sfrac 12}U_k(0)\Big\|^2_{L^2(\BR^d)}$.
According to (4.11) we have,
$$\rondun\le C\Big\|E_{-\sfrac 12}\big(1-\chi   _3\big)h^{-\sfrac
12}_k\theta_1\big(h^2_kP_D\big)\widetilde
u^0_k\Big\|^2_{L^2}=C\Big\|E_{-\sfrac 12}\big(1-\chi   _3\big)P_D^{\sfrac
14}\theta\big(h^2_kP_D\big)\widetilde u\,^0_k\Big\|^2_{L^2}$$
Introducing $S=E_{-\sfrac 12}\big(1-\chi   _3\big)P_D^{\sfrac 14}$ we can
write
$$\eqalign{
\rondun&\le C\left(S\theta\big(h^2_kP_D\big)\widetilde u\,^0_k,
S\theta\big(h^2_kP_D\big)\widetilde u\,^0_k\right)
=\left(\theta\big(h^2_kP_D\big)S^*S\theta\big(h^2_kP_D\big)\widetilde
u^0_k, \widetilde u\,^0_k\right)\hc
\rondun&\le
C\Norme{\theta\big(h^2_kP_D\big)S^*S\theta\big(h^2_kP_D\big)\widetilde
u^0_k}_{L^2}\Norme{\widetilde u\,^0_k}_{L^2}\hc}$$
Now $$S^*S=P_D^{-\sfrac 14}P_D^{\sfrac
12}\big(1-\chi   _3\big)A_{-1}\big(1-\chi   _3\big)P_D^{\sfrac
14}=P_D^{-\sfrac 14}Q\big(1-\chi   _3\big)P_D^{\sfrac 14}$$ where $Q$ has
been defined in Lemme 4.2 and $A_{-1}=E^*_{-\sfrac 12}E_{-\sfrac 12}$.

Using (3.4) we obtain
$$\rondun\le C\Normes{h_k^{\sfrac
12}\theta_2\big(h^2_kP_D\big)Qh_k^{-\sfrac
12}\big(1-\chi   _3\big)\theta_1\big(h^2_kP_D\big)\widetilde
u^0_k}_{L^2}\norme{\widetilde u\,^0_k}_{L^2}$$
Since the operators $\theta_j\big(h^2_kP_D\big)$, $j=1,2$, are uniformly
bounded in $L^2(\Omega)$, using Lemma 4.2 and (3.7) we obtain
$$\Normes{E_{-\sfrac 12} U_k(0)}^2_{L^2(\BR^d)}\le C\norme{\widetilde
u\,^0_k}_{L^2(\Omega)}^2\le C\leqno{(4.13)}$$
with a uniforn constant $C>0$.

We claim now that we have (see (4.11)), uniformly in $k\ge 1$,
$$\Int_0^T\Normes{\scal
x^{\sfrac{1+\nu}2}E_{-1}G_k(t)}^2_{L^2(\Br^d)}\dt={\cal
O}(1).\leqno{(4.14)}$$
By (4.7) we can write $G_k=\big[\widetilde P,\chi   _3\big]\chi   _1w_k$.
Moreover the symbolic calculus shows that the symbol of $\big[\widetilde
P,\chi   _3\big]$ belongs to $S\left(\dfrac{\scal\xi^2}{\scal
x\scal\xi},g\right)$.

It follows that the symbol of $\scal
x^{\sfrac{1+\nu}2}E_{-1}\big[\widetilde P,\chi   _3\big]$ belongs to
$S(M,g)$ where
$$M(x,\xi)=\dfrac{\scal x^{\sfrac{1+\nu}2}}{\scal
x+\scal\xi}\cdot\dfrac{\scal\xi^2}{\scal x\scal\xi}\le C.$$
This operator is therefore $L^2$ bounded, so using (3.7) we obtain
$$\Int_0^T\Normes{\scal
x^{\sfrac{1+\nu}2}E_{-1}G_k(t)}^2_{L^2(\BR^d)}\dt\,\le\, 
C'\Int_0^T\Norme{\chi   _1w_k(t)}^2_{L^2(\Omega)}\dt\,\le\,C'$$
which proves (4.14).

Using (4.10), (4.12), (4.13) and (4.14) we conclude that
$$\Int_0^T\Norme{\chi   \big(1-\chi   _3\big)w_k(t)}^2_{L^2(\Omega)}\dt={\cal
O}(1).\leqno{(4.15)}$$
Since by (3.7) we have
$\Int_0^T\Norme{\chi   _1w_k(t)}^2_{L^2(\Omega)}\dt=1$ and since by (3.3)
and (4.6) we have $\chi   _1+\big(1-\chi   _3\big)\ge 1$ we obtain (4.1). The
proof of Proposition 4.1 is complete.\hfill\cqfd

\vpt
\noindent {\soustitre 5. The measure $\mu$ and its properties}
\vpt
We shall set
$$\left\{\eqalign{
&\sconj w_k(t)=\bun_\Omega w_k(t),\hc
&W_k=\bun_{[0,T]}\,\sconj w_k.\hc}\right.\leqno{(5.1)}$$
It follows from Proposition 4.1 that the sequence $\big(W_k\big)$ is
bounded in $L^2\left(\BR_t,L^2_{loc}\big(\BR^d\big)\right)$.

Now to a symbol $a=a(x,t,\xi,\tau)\in{\cal
C}^\infty_0\big(T^*\BR^{d+1}\big)$ we associate the semi-classical
pseudo-differential operator (pdo) by the formula
$$\eqalign{{\cal O}p(a)\big(x,t,hD_x,h^2D_t\big)v(x,t)&=\cr
(2\pi
h)^{-(d+1)}\Intd&\erom^{i\left(\sfrac{x-y}h\xi+\sfrac{t-s}{h^2}\tau\right)}
\varphi(y)
a(x,t,\xi,\tau)v(y,s)\dy\ds\drm\xi\dtau\cr} \leqno{(5.2)}$$
where $\varphi\in{\cal C}^\infty_0\big(\BR^d\big)$ is equal to one on a
neighborhood of the $x$-projection of the support of $a$.

We note that by the symbolic calculus the operator ${\cal O}p(a)$ is,
modulo operators bounded in $L^2$ by ${\cal O}\big(h^\infty\big)$,
independant of the function $\varphi$. The following result is
classical and introduces the notion of semi-classical defect measure.

\proclaim Proposition 5.1.
There exists a subsequence $\big(W_{\sigma(k)}\big)$ and a Radon measure
$\mu$ on $T^*\BR^{d+1}$ such that for every $a\in{\cal
C}^\infty_0\big(T^*\BR^{d+1}\big)$ one has
$$\Lim_{k\to
+\infty}\left({\cal
O}p(a)\big(x,t,h_{\sigma(k)}D_x,h^2_{\sigma(k)}
D_t\big)W_{\sigma(k)},W_{\sigma(k)}\right)_{L^2(\BR^{d+1})}=\scal{\mu,a}.$$

Here are the two main properties of the measure $\mu$ which will be used
later on.

\proclaim Theorem 5.2.
    The support of $\mu$ is contained in the set
$$\Sigma=\big\{(x,t,\xi,\tau)\in
T^*\BR^{d+1}\moins\{0\}~:~x\in\conj\Omega, t\in\big[0,T\big]\hbox{ and }
\tau+p(x,\xi)=0\big\}.$$

\vpt 
\noindent {\bf Proof}
\vpt
See section 9.2 in the appendix.\hfill\cqfd

To state the propagation result let us recall some notations. Let
$M=\Omega\times\BR_t$. We set
$$T^*_bM=T^*M\moins\{0\}\cup T^*\pa M\moins\{0\}.$$
We have a natural application of restriction
$$\pi~:~T^*\BR^{d+1}_{|\conj M}\moins\{0\}\to T^*_bM$$
which is the identity on $T^*\BR^{d+1}_{|M}\moins\{0\}$ (see section 9.1
for details).

With $\Sigma$ defined in Theorem 5.2 we set $\Sigma_b=\pi(\Sigma)$.
%
%\noindent 
The measure $\mu$ has its support in $\Sigma\subset
T^*\BR^{d-1}_{|\conj 
M}\moins\{0\}$ while for $\zeta\in\Sigma_b$ the generalized
bicharacteristic $\Gamma(t,\zeta)$ lives in $\Sigma_b$.

Then we can state an important result of this paper.

\proclaim Theorem 5.3.
 Let $\zeta\in\Sigma_b$ and $s_1,s_2\in\BR$. Then we have
 $$\pi^{-1}\big(\Gamma(s_1,\zeta)\big)\cap\Supp\mu=\vide\iff
\pi^{-1}\big(\Gamma(s_2,\zeta)\big)\cap\Supp\mu=\vide$$

For the proof, see the appendix, section 9.2.

\vpt
\noindent {\soustitre 6. The measure $\mu$ does not vanish
identically}
\vpt
The purpose of this section is to prove the following results.

Let $A\ge 1$, $R\ge 1$, $\psi_A\in{\cal C}^\infty_0(\BR)$,
$\Phi_R\in{\cal C}^\infty_0(\Br)$ be such that $0\le\psi_A,\Phi_R\le 1$
and
$$\psi_A(\tau)=1\hbox{ if } \abs\tau\le A,\pv \Phi_R(t)=1\hbox{ if } \abs
t\le R.\leqno{(6.1)}$$

\proclaim Proposition 6.1.
    There exist positive constants $A_0$, $R_0$, $k_0$ such that
 $$\Int_\BR\Norme{\psi_A\big(h^2_kD_t\big)\Phi_R
 \big(h^2_k\Delta\big)\bun_{[0,T]}\chi_1\sconj
 w_k(t)}^2_{L^2(\Br^d)}\dt\ge \dfrac  
12, \leqno{(6.2)}$$
 when $A\ge A_0$, $R\ge R_0$, $k\ge k_0$. Here $\chi   _1$, $\sconj w_k$
have been defined in (3.3), (5.1) and $\Delta$ is the usual Laplacien.

  %\newpage
\proclaim Corollary 6.2.
    The measure $\mu$ defined in Proposition 5.1 does not vanish
identically.

\vpt    
\noindent {\bf Proof}
\vpt
Let $\widetilde\chi   _1\in{\cal C}^\infty_0\big(\BR^d\big)$ be such that
$\widetilde \chi   _1=1$ on $\Supp\chi   _1$.
Let $\varphi=\varphi(t)\in{\cal
C}_0^\infty(\Br)$ and
$a(x,t,\xi,\tau)=\varphi(t)\chi_1(x)
\psi^2_A(\tau)\Phi^2_R\big(\abs\xi^2\big)\chi   _1$.  
It follows from (6.2) that
$$\left(a\big(x,t,hD_x,h^2D_t\big)\widetilde\chi_1\bun_{[0,T]}\sconj
w_k(t),\bun_{[0,T]}\sconj w_k\right)_{L^2(\BR^{d+1})}\ge\dfrac 13.$$
Since the left hand side with the subsequence $\sigma(k)$ tends to
$\scal{\mu,a}$ when $k\to +\infty$ the Corollary follows.
\hfill\cqfd

\vpt
\noindent {\bf Proof of Proposition 6.1}
\vpt
We shall need the following Lemma.

\proclaim Lemma 6.3.
   Let $\theta\in{\cal C}^\infty_0(\Br)$, $\chi   \in{\cal
C}^\infty_0\big(\BR^d\big)$. Then there exists $C>0$ such that
$$\eqalign{
(i)\quad
&\Norme{\big[\theta\big(h^2P_D\big),\chi   \big]u}^2_{L^2(\Omega)}\le
C\norme{hu}^2_{L^2(\Omega)},\hbox{ \hskip 4cm} \hc
(ii)\quad &\Norme{\pa_j\theta\big(h^2P_D\big)u}^2_{L^2(\Omega)}\le
C\Norme{h^{-1}u}_{L^2(\Omega)},\hc
(iii)\quad&
\Norme{\pa_j\big[\theta\big(h^2P_D\big),\chi   \big]u}_{L^2(\Omega)}\le
C\norme u_{L^2(\Omega)},\hc}$$
for all $j=1,\cdots,d$, $h>0$ and $u\in L^2(\Omega)$.

\vpt
 \noindent {\bf Proof}
\vpt
See the Appendix, section 9.3.\hfill\cqfd

Let us set
$$\left\{\eqalign{
I&=\left(\Id-\psi_A\big(h^2D_t\big)\right)\bun_{[0,T]}\chi   _1w_k\hc
\widetilde\psi(\tau)&=\dfrac{1-\psi(\tau)}\tau.\hc}\right.\leqno{(6.3)}$$
Then $\widetilde\psi\in L^\infty(\Br)$ and
$\abs{\widetilde\psi(\tau)}\le\dfrac 1A$ for all $\tau\in\BR$.

Now we can write
$I=\widetilde\psi\big(h^2_kD_t\big)h^2_kD_t\big(\bun_{[0,T]}\chi_1w_k\big)$.
Using (3.6) and the fact that \break $D_t\widetilde u_k=P_D\widetilde u_k $ we
deduce that
$$\left\{\eqalign{
I&=\rondun+\rondtwo\hc
\rondun&={1\over i}\widetilde\psi\big(h^2_kD_t\big)\chi_1h^2_k
\big(w_k(0)\delta_{t=0}-w_k(T)\delta_{t=T}\big)\hc 
\rondtwo&=-\widetilde\psi\big(h^2_kD_t\big)\bun_{[0,T]}\chi_1h^2_kP_Dh^{-\sfrac
12}_k\theta_1\big(h^2_kP_D\big)\widetilde u_k.\hc}\right.\leqno{(6.4)}$$

\vpt
{\bf  Estimate of \rondun}
\vpt
If $a\in\BR$, we have
$\widetilde\psi\big(h^2_kD_t\big)\delta_{t=a}=\conj\CF
\left(\widetilde\psi\big(h^2_k\tau\big)\erom^{-ia\tau}\right)$, 
so by Parseval formula we have
$$\Norme{\widetilde\psi\big(h^2_kD_t\big)\delta_{t=a}}^2_{L^2(\BR)}
=c_n\Int\abs{\widetilde\psi\big(h^2_k\tau\big)}^2\dtau=c_nh^{-2}_k
\Int\dfrac{\abs{\psi_A(\tau)-1}^2}{\abs\tau^2}\dtau$$
It follows from (3.6) that
$$\Int_\BR\norme{\rondun}^2_{L^2(\Omega)}\dt\,\le\,
C\,h^4_k\,h^{-2}_k\,h^{-1}_k\left(\Norme{\widetilde
u_k(0)}^2_{L^2}+\Norme{\widetilde u_k(T)}^2_{L^2}\right)$$
Applying the energy estimate and (3.7) we obtain
$$\Int_\BR\norme\rondun^2_{L^2}\dt\,\le\,C\,h_k\Norme{\widetilde
u^0_k}^2_{L^2}=o(1).\leqno{(6.5)}$$

\vpt
{\bf  Estimate of \rondtwo}
\vpt
Let $\widetilde\theta\in{\cal C}^\infty_0(\Br)$  be such that
$\widetilde\theta=1$ on the support of $\theta_1$. Then we can write
with $\widetilde\theta_1(t)=t\,\widetilde\theta(t)$
$$\rondtwo=-\widetilde\psi\big(h^2_kD_t\big)\left[\chi_1,\widetilde
\theta_1\big(h^2_kP_D\big)\right]\bun_{]0,T]}h^{-\sfrac 
12}_k\theta_1\big(h^2_kP_D\big)\widetilde
u_k-\widetilde\psi\big(h^2_kD_t\big)\widetilde\theta_1
\big(h^2_kP_D\big)\bun_{[0,T]}\chi   _1w_k(t)$$ 
Using Lemma 6.3$\,(i)$ and the fact that,
$$\Norme{\widetilde\psi\big(h^2_kD_t\big)}_{L^2(\BR)\to L^2(\Br)}={\cal
O}\left(\dfrac 1A\right),\pv
\Norme{\widetilde\theta_1\big(h^2_kP_D\big)}_{L^2(\Omega)\to L^2(\Omega)}={\cal
O}(1)$$
uniformly in $k$ we obtain,
$$\Int_\BR\norme\rondtwo^2_{L^2(\Omega)}\dt\le\dfrac
CA\left(\Int_0^T\Normes{h^{\sfrac 12}_k\widetilde
u_k(t)}^2_{L^2(\Omega)}\dt+\Int_0^T\Norme{\chi_1w_k(t)}^2_{L^2(\Omega)}
\dt\right)$$
Using the energy estimate and (3.7) we deduce that
$$\Int_\BR\norme\rondtwo_{L^2(\Omega)}^2\dt=o(1)+{\cal O}\left(\dfrac
1A\right).\leqno{(6.6)}$$
Taking $k$ and $A$ sufficiently large and using (3.7), (6.3), (6.4),
(6.5), (6.6) we obtain
$$\Int_\BR
\Norme{\psi_A\big(h^2_kD_t\big)\bun_{[0,T]}\chi_1w_k(t)}^2_{L^2(\Omega)}\ge  
\dfrac 12.\leqno{(6.7)}$$
Now with $\Phi_R$ defined in (6.1) we set
$$\hbox{II}=h^{-\sfrac
12}_k\left(\Id-\Phi_R\big(h^2_k\Delta\big)\right)\psi_A
\big(h^2_kD_t\big)\bun_{[0,T]}\chi_1\sconj 
w_k(t).
\leqno{(6.8)}$$
Since $\Supp\big(1-\Phi_R(t)\big)\subset\big\{t\in\BR~:~\abs t\ge
R\big\}$ we have by Fourier transform
$$\Norme{\big(\Id-\Phi_R\big(h^2\Delta\big)\big)
h^{-1}_kv}^2_{L^2(\BR^d)}\le\dfrac  
CR\Sum_{j=1}^d\Norme{\pa _jv}^2_{L^2(\BR^d)}\pv,\pv v\in
H^1\big(\Br^d\big).
\leqno{(6.9)}$$
Now by (3.6) we have $h^{-\sfrac 12}_k\sconj w_k=h^{-1}_k\sconj v_k$,
$v_k=\theta_1\big(h^2_kP_D\big)\widetilde u_k$.

Thus applying (6.9) we obtain
$$\Int_\BR\norme{\hbox{II}}^2_{L^2(\BR^d)}\dt\le\dfrac
CR\Sum_{j=1}^d\Int_\BR\Norme{\pa_j\psi_A\big(h^2D_t
\big)\bun_{[0,T]}\chi_1\sconj 
v_k}^2_{L^2(\BR^d)}.$$
Since $v_k\in H^1_0(\Omega)$ we have
$$\pa_j\big(\chi   _1\sconj
v_k\big)=\pa_j\big(\bun_\Omega\chi_1
v_k\big)=\bun_\Omega\pa_j\big(\chi   _1v_k\big)$$  
It follows that
$$\Int_\BR\norme{\hbox{II}}^2_{L^2(\Br^d)}\dt\,\le\,\dfrac
CR\Sum_{j=1}^d\Int_\BR\Norme{\pa_j\psi_A\big(h^2_kD_t\big)\bun_{[0,T]}\chi
_1v_k}^2_{L^2(\Omega)}\dt$$ 
Let $\widetilde\theta\in{\cal C}^\infty_0(\Br)$ be such that
$\widetilde\theta=1$ near the support of $\theta_1$. Then
$\big(1-\widetilde\theta(t)\big)\theta_1(t)=0$. We first consider
$$\rondun=\Int_\BR\Norme{\pa_j
\widetilde\theta\big(h^2_kP_D\big)\psi_A\big(h^2_k  
D_t\big)\bun_{[0,T]}\chi   _1v_k}^2_{L^2(\Omega)}\dt.\leqno{(6.10)}$$
Using Lemma 6.3$\,(ii)$ we obtain
$$\rondun\le
C\Int_\BR\Norme{\psi_A\big(h^2_kD_t\big)h^{-1}_k\bun_{[0,T]}\chi
_1v_k(t)}^2_{L^2(\Omega)}\dt\le 
C'\Int_0^T\Normes{h^{-\sfrac 12}_k\chi   _1w_k(t)}^2_{L^2(\Omega)}\dt$$
so by (3.6) we obtain
$$\rondun\le C'h^{-1}_k.\leqno{(6.11)}$$
It remains to consider
$$\rondtwo=\Int_\BR\Norme{\pa_j
\left(\Id-\widetilde\theta\big(h^2_kP_D\big)\right)  
\psi_A\big(h^2_kD_t\big)\bun_{[0,T]}\chi 
_1v_k(t)}^2_{L^2(\Omega)}\dt.\leqno{(6.12)}$$ 
Since $v_k(t)=\theta_1\big(h^2_kP_D\big)\widetilde u_k$ and
$\left(\Id-\widetilde\theta
\big(h^2_kP_D\big)\right)\theta_1\big(h^2_kP_D\big)=0$ 
we obtain
$$\rondtwo\le\Int_\BR\Norme{\pa_j\left[\widetilde\theta\big(h^2_kP_D\big),
\chi_1\right]\psi_A\big(h^2_kD_t\big)\bun_{[0,T]}\widetilde\chi
_1v_k(t)}^2_{L^2(\Omega)}\dt$$ 
where $\widetilde\chi   _1\in{\cal C}^\infty_0\big(\conj\Omega\big)$,
$\widetilde\chi   _1=1$ on $\Supp\chi   _1$.

By Lemma 6.3$\,(iii)$ we obtain
$$\rondtwo\le C\Int_\BR\Norme{\psi_A\big(h^2_k
D_t\big)\bun_{[0,T]}\widetilde\chi_1\theta_1\big(h^2_kP_D\big)\widetilde
u_k(t)}^2_{L^2(\Omega)}\dt.$$
Since the operator $\psi_A\big(h^2_kD_t\big)$ is uniformly $L^2$ bounded
we obtain by the energy estimate
$$\rondtwo\le C'\Int_0^T\Norme{\widetilde u_k(t)}^2_{L^2}\dt={\cal
O}(1).$$
It follows from (6.10), (6.11), (6.12) that
$$\Int_\BR\Norme{\hbox{II}}^2_{L^2}\dt\,\le\,\dfrac
CR\big(h^{-1}_k+{\cal O}(1)\big).\leqno{(6.13)}$$
Using (6.8) we deduce that
$$\Int_\BR\Norme{\left(\Id-\Phi_\BR\big(h^2\Delta\big)\right)
\psi_A\big(h^2D_t\big)\bun_{[0,T]}\chi   _1\sconj
w_k(t)}^2_{L^2(\Br^d)}\dt\,\le\,\dfrac CR\left(1+{\cal
O}\big(h_k\big)\right).\leqno{(6.14)}$$
Taking $R$ sufficiently large and using (6.7) we obtain
$$\Int_\BR\Norme{\Phi_\BR\big( h^2\Delta\big) \psi_A\big( h^2D_t\big)
\bun_{[0,T]}\chi_1\sconj  
w_k(t)}^2_{L^2}\dt\,\ge\,\dfrac 13,$$
which is (6.2). The proof of Proposition 6.1 is complete.\hfill\cqfd

\vpt
\noindent {\soustitre 7.  The measure $\mu$ vanishes in
the incoming set}
\vpt
We pursue here our reasoning by contradiction in proving that the
measure $\mu$ vanishes in the incoming set. Let us state the main result
of this section.

Let $\widetilde P$ be the operator defined by (4.3) satisfying the
conditions (4.5).

\proclaim Theorem 7.1.
    Let $m_0=\big(x_0,t_0,\xi_0,\tau_0\big)\in T^*\BR^{d+1}$ be such
$\xi_0\ne 0$,
$\tau_0+\widetilde p\big(x_0,\xi_0\big)=0$, $\abs{x_0}\ge 3R_0,\pv
\Sum_{j,k=1}^d \widetilde a^{jk}\big(x_0\big)x_{0j}\xi_{0k}\le
-3\delta\abs{x_0}\abs{\xi_0}$
 (for some $\delta>0$ small enough). Then $m_0\notin\Supp\mu$

The rest of this section will be devoted to the proof of this result. Il
will be a consequence of an estimate which will be proved in
constructing an appropriate escape function and will require several
Lemmas.

\proclaim Lemma 7.2.
    Let us set $e_0(x,\xi)=\Sum_{j,k=1}^d \widetilde
a^{jk}(x)x_j\dfrac{\xi_k}{\scal\xi}$. Then there exist positive
constants $R$, $C_0$, $C_1$ such that
 $$H_{\tilde p} e_0(x,\xi)\ge C_0\abs\xi-C_1,\quad
\forall(x,\xi)\in T^*\BR^d,\quad \abs x\ge R$$
 where $H_{\tilde p}$ denotes the Hamiltonian field of $\widetilde
P$.

\vpt
\noindent {\bf Proof}
\vpt
It is an easy computation which uses the conditions $(ii)$ and $(iii)$
in (4.5).\hfill\cqfd

\proclaim Lemma 7.3.
    Under condition $(v)$ in (4.5) there exist $e\in S\big(\scal x,g\big)$
and positive constants $C$, $C'$, $R'$ such that
 $$\eqalign{
(a)\quad &H_{\tilde p}e(x,\xi)\ge C\abs\xi-C',\quad
\forall(x,\xi)\in T^*\BR^d,\hc
(b)\quad &e(x,\xi)=e_0(x,\xi)\pv\hbox{if}\pv\abs x\ge R'.\hc}$$

\vpt
 \noindent {\bf Proof}\pv See   {Doi} [D3].\hfill\cqfd
\vpt
The symbol $e$ is an escape function. However it is not adapted to our
situation because its Poisson bracket with our potentiel $\widetilde V$
(see (4.5)$\,(iv)$) belongs to $S\Big(\dfrac{\scal
x^2}{\scal\xi},g\Big)$ so does not correspond to an operator bounded in
$L^2$ which will be required later on. We shall describe below a
construction by   {Doi} [D3] which will take care of this problem.

Let $\psi\in{\cal C}^\infty(\Br)$ be such that $0\le\psi\le 1$ and
$$\psi(t)=1\hbox{ if } t\ge 2\eps,\quad
\Supp\psi\subset\big[\eps,+\infty\big[,\quad \psi'(t)\ge 0\ \
\forall t\in\BR,\leqno{(7.1)}$$
where $\eps>0$ is a small constant choosen later on.

We set
$$\left\{\eqalign{
\psi_0(t)&=1-\psi(t)-\psi(-t)=1-\psi(|t|)\hc
\psi_1(t)&=\psi(-t)-\psi(t)=-(\sgnt)\psi\big(\abs t\big).\hc}\right.
\leqno{(7.2)}$$
Then $\psi_j\in{\cal C}^\infty(\Br)$, $j=1,2$, and we have
$$\left\{\eqalign{
\psi'_0(t)&=-(\sgnt)\psi'\big(\abs t\big)\hc
\psi'_1(t)&=-\psi'\big(\abs t\big).\hc}\right.
\leqno{(7.3)}$$
Let $\chi   \in{\cal C}^\infty(\Br)$ be such that $0\le\chi   \le 1$ and
$$\chi   (t)=1\hbox{ if } t\le\dfrac\rho 2,\quad \chi   (t)=0\hbox{ if }
t\ge\rho,\quad \rho>0\hbox{ small}.
\leqno{(7.4)}$$
With $e$ defined in Lemma 7.3 we set
$$-\la=\Big(\dfrac e{\scal x}\psi_0\Big(\dfrac e{\scal
x}\Big)-\left(M_0-\scal e^{-\nu}\right)\psi_1\Big(\dfrac e{\scal
x}\Big)\Big)\chi   \Big(\dfrac{\scal x}{\sqrt{\widetilde
p(x,\xi)}}\Big)
\leqno{(7.5)}$$
where $\nu>0$ is an arbitrary small constant and $M_0$ is a sufficiently
large constant.

\proclaim Lemma 7.4
{\rm (  {Doi} [D3], Lemma 8.3)}.
$$\eqalign{
& (i)\,\la\in S(1,g)\cr
&(ii)\,\big[\widetilde P, {\cal O}p^w(\la)\big]-\dfrac
1i\big(H_{\tilde p}\la\big)^w\in{\cal O}p^w S(1,g)\cr
&(iii)\,\hbox{There exists}\ M_0>0\  \hbox{such that for any}\
\nu>0\ \hbox{there 
exist}\ C>0,\ C'>0\cr
&\quad \quad \, \hbox{such that}\
 -H_{\tilde p}\la(x,\xi)\ge C\scal x^{-1-\nu}\big(\abs
x+\abs\xi\big)-C',\quad \forall(x,\xi)\in T^*\BR^d.\cr}\leqno{(7.6)}$$

We must now localize this escape function near the incoming set.

We shall need the following Lemma. Let us set
$$a(x,\xi)=\Sum_{j,k=1}^d\widetilde a^{jk}(x)x_j\xi_k\leqno{(7.7)}$$
Let $\xi_0\ne 0$ be defined in Theorem 8.1.

\proclaim Lemma 7.5. 
    There exists a symbol $\Phi\in S(1,g)$ such that $0\le\Phi\le 1$
    and
$$\eqalign{
&(i)\,\Supp\Phi\subset\Big\{(x,\xi)\in T^*\BR^d~:~\abs x\ge
2R_0,\pv a(x,\xi)\le-\dfrac\delta 2\abs x\,\abs\xi,\pv
\abs\xi\ge\dfrac{\abs{\xi_0}}4\Big\},\cr
&(ii)\,\Big\{(x,\xi)~:~\abs x\ge\dfrac 52R_0,\quad
a(x,\xi)\le-\delta\abs x\,\abs\xi,\pv
\abs\xi\ge\dfrac{\abs{\xi_0}}
2\Big\}\subset\big\{(x,\xi)~:~\Phi(x,\xi)=1\big\},\cr
&(iii)\,\Phi(x,h\xi)=\Phi(x,\xi)\ \hbox{\sl when\ }
\abs{h\xi}\ge\dfrac{\abs{\xi_0}}2\ \hbox{ and\ } 0<h\le 1,\cr
&(iv)\,H_{\tilde p}\Phi(x,\xi)\le 0\hbox{\  on the support of\ }
\la,\cr
&(v)\,\la(x,\xi)\ge 0\hbox{\ on the support of\ } \Phi.\cr}$$

\vpt    
\noindent {\bf Proof}
\vpt
Let $\varphi_j$, $j=1,2,3$, be such $\varphi_j\in{\cal C}^\infty(\Br)$,
$0\le \varphi_j\le 1$ and
$$\left\{\eqalign{
\varphi_1(s)&=0\hbox{ if } s\le R_0,\pv\varphi_1(s)=1\hbox{ if }
s\ge\dfrac 52R_0,\pv \varphi_1\hbox{ increasing},\hc
\varphi_2(s)&=0\hbox{ if } s\ge -\dfrac 12\delta,\pv
\varphi_2(s)=1\hbox{ if } s\le -\delta,\pv \varphi_2\hbox{
decreasing},\hc
\varphi_3(s)&=0\hbox{ if } s\le\dfrac 14\abs{\xi_0},\pv
\varphi_3(s)=1\hbox{ if } s\ge \dfrac 12\abs{\xi_0}.\hc}\right.
\leqno{(7.8)}$$
Let us set
$$\Phi(x,\xi)=\varphi_1\big(\abs
x\big)\varphi_2\left(\dfrac{a(x,\xi)}{\abs
x\,\abs\xi}\right)\varphi_3\big(\abs\xi\big).
\leqno{(7.9)}$$
Then $(i)$ and $(ii)$ follow immediatly. Now if
$\abs{h\xi}\ge\dfrac{\abs{\xi_0}}2$ then
$\abs{\xi}\ge\dfrac{\abs{\xi_0}}{2h}\ge\dfrac{\abs{\xi_0}}2$ so
$\varphi_3\big(h\abs\xi\big)=\varphi_3\big(\xi\big)=1$ and $(iii)$
follows.

Let us prove $(iv)$. We have
$$\left\{\eqalign{
&H_{\tilde p}\Phi(x,\xi)=\rondun+\rondtwo+\rondtrois,\hc
&\rondun=\varphi'_1\big(\abs x\big)H_{\tilde p}\abs
x\varphi_2\left(\dfrac a{\abs
x\,\abs\xi}\right)\varphi_3\big(\abs\xi\big),\hc
&\rondtwo=\varphi_1\big(\abs x\big)\varphi'_2\left(\dfrac a{\abs
x\,\abs\xi}\right)H_{\tilde p}\left(\dfrac a{\abs
x\,\abs\xi}\right)\varphi_3\big(\abs\xi\big),\hc
&\rondtrois=\varphi_1\big(\abs x\big)\varphi_2\left(\dfrac a{\abs
x\,\abs\xi}\right)\varphi'_3\big(\abs\xi\big)H_{\tilde p} \abs
\xi.\hc}\right.
\leqno{(7.10)}$$
According to (7.4) and (7.5) we have $\displaystyle\widetilde
p(x,\xi)\ge\dfrac1{\rho^2}\scal x^2\ge\dfrac 1{\rho^2}$ on the support
of $\la$. Therefore we can choose~$\rho$ so small that $\abs\xi>\dfrac
12\abs{\xi_0}$ on the support of $\la$. It
follows that $\rondtrois=0$ on this set. Now an easy computation shows
that $H_{\tilde p}\abs x=\dfrac{2a(x,\xi)}{\abs x}$ when $\abs x\ge
R_0$ which implies that
$$\rondun=2\varphi'_1\big(\abs x\big)\dfrac a{\abs
x}\varphi_2\left(\dfrac a{\abs x\,\abs\xi}\right)\varphi_3\big(\abs
\xi\big).$$
On the support of $\varphi_2$ we have $a\le-\dfrac 12\delta\abs
x\,\abs\xi$. Since $\varphi'_1\ge 0$, $\varphi_2\ge 0$, $\varphi_3\ge
0$, we conclude that
$$\rondun\le 0.\leqno{(7.11)}$$
Let us look to \rondtwo. First of all we have on the support of $\Phi$
$$H_{\tilde p}\left(\dfrac a{\abs x\,\abs\xi}\right)=\dfrac 1{\abs
x\,\abs\xi}H_{\tilde p}a+a\ H_{\tilde p}\left(\dfrac1{\abs
x\,\abs\xi}\right).
\leqno{(7.12)}$$
Since we have (see (4.5)) $\big(\widetilde a^{jk}(x)\big)\ge C\ \Id$,
$\abs{\nabla_x\,\widetilde a^{jk}(x)}=o\big(\abs x^{-1}\big)$ as $\abs
x\to +\infty$ and $\abs x\ge R_0$ on the support of $\Phi$, taking $R_0$
large enough we obtain by an easy computation
$$H_{\tilde p} a(x,\xi)\ge C_0\abs\xi^2\hbox{ on } \Supp\Phi.
\leqno{(7.13)}$$
We also obtain
$$H_{\tilde p}\left(\dfrac 1{\abs
x\,\abs\xi}\right)=-2\dfrac{a(x,\xi)}{\abs\xi\,\abs x^3}+o\left(\dfrac
1{\abs x^2}\right)\hbox{ as } \abs x\to +\infty.
\leqno{(7.14)}$$
It follows from (7.12), (7.13), (7.14) and $\abs a\le C\abs x\,\abs\xi$
that
$$H_{\tilde p}\left(\dfrac a{\abs x\,\abs\xi}\right)\ge
C_0\dfrac{\abs\xi}{\abs x}-2\dfrac{a^2(x,\xi)}{\abs\xi\,\abs
x^3}+o(1)\dfrac{\abs\xi}{\abs x}.$$
On the support of $\varphi'_2\left(\dfrac a{\abs x\,\abs\xi}\right)$ we
have, by (8.8), $-\delta\le\dfrac a{\abs x\,\abs \xi}\le-\dfrac
12\delta$. It follows that $\abs a\le\delta\abs x\,\abs\xi$ so
$$-2\dfrac{a^2}{\abs \xi\,\abs x^3}\ge -\delta^2\dfrac{\abs\xi}{\abs
x}.$$
Moreover on the support of $\varphi_1\big(\abs x\big)$ we have $\abs
x\ge R_0$. So taking $R_0$ large enough and $\delta$ small, we obtain
$$H_{\tilde p}\left(\dfrac a{\abs
x\,\abs\xi}\right)\ge\dfrac{C_0}2\dfrac{\abs\xi}{\abs x}.$$
Since, by (7.5), we have $\varphi'_2\left(\dfrac a{\abs
x\,\abs\xi}\right)\le 0$, we conclude that
$$\rondtwo\le 0.\leqno{(7.15)}$$
The claim $(iv)$ in Lemma 7.5 follows then from (7.11), (7.15) and
(7.10) since $\rondtrois=0$ on $\Supp\la$.

Finally let us look to the claim $(v)$.

On the support of $\Phi$ we have $\abs x\ge R_0$, $\abs\xi\ge\dfrac
14\abs{\xi_0}$ and $a(x,\xi)\le-\dfrac 12\delta\abs x\,\abs\xi$. It
follows that $\scal x\le\sqrt2\abs x$, $\scal\xi\le C\abs\xi$ and
$a(x,\xi)\le -C'\delta\scal x\scal\xi$. Moreover since $\abs x\ge R_0$,
taking $R_0$ large enough, we deduce from Lemma 7.3 that
$e(x,\xi)=e_0(x,\xi)=\dfrac{a(x,\xi)}{\scal\xi}$ by (7.7). It follows
that $\dfrac e{\scal x}\le -C'\delta$ which implies that $\dfrac{\abs
e}{\scal x}\ge C'\delta$. Using (7.1), (7.2) and taking $\eps\ll\delta$
we see that $\psi_0\left(\dfrac e{\scal x}\right)=0$ and
$\psi_1\left(\dfrac e{\scal x}\right)\ge 0$. It follows from (7.5) that
$-\la=-\big(M_0-\scal e^{-\nu}\big)\psi_1\left(\dfrac e{\scal
x}\right)\chi   \le 0$.

The proof of Lemma 7.5 is complete.\hfill\cqfd

\proclaim Corollary 7.6.
    Let $\la_1=\Phi^2\la$ where $\la$ has been defined in Lemma 7.4. Then
$$\eqalign{
&(i)\,\la_1\in S(1,g),\cr
&(ii)\,\big[\widetilde P,\la^w_1\big]-\dfrac 1i{\cal
O}p^w\big(H_{\tilde p}\la_1\big)\in {\cal O}p^w S(1,g),\cr
&(iii)\,\hbox{ There exist two positive constants } C,\ C' \hbox{ such
that }\cr
&\quad\quad\,-H_{\tilde p}\la_1\ge C\scal x^{-1-\nu}\Phi^2(x,\xi)\big(\abs
x+\abs \xi\big)-C'\Phi^2(x,\xi).\cr}\leqno{(7.16)}$$

\vpt
\noindent {\bf Proof}
\vpt
$(i)$, $(ii)$ follow from Lemma 7.4 and the fact that $\Phi^2\in
S(1,g)$. Let us look to $(iii)$. We have
$$-H_{\tilde p}\la_1=\big(-H_{\tilde p}\la\big)\Phi^2-2\la\Phi
H_{\tilde p}\Phi.$$
By Lemma 7.5 we have $H_{\tilde p}\Phi\le 0$ on $\Supp\la$ and
$\la\ge 0$ on $\Supp\Phi$.
It follows that $-2\la\Phi H_{\tilde p}\Phi\ge 0$. Thus $(iii)$
follows from (7.6).\hfill\cqfd

Let now $\big(x_0\,,\xi_0\big)$ be given as in Theorem 7.1. We set
$$V_{(x_0,\xi_0)}=\big\{(x,\xi)\in
T^*\BR^d~:~\abs{x-x_0}+\abs{\xi-\xi_0}\le\eps_0\big\}.$$
Since $\abs{x_0}\ge 3R_0$, $a\big(x_0,\xi_0\big)\le
-3\delta\abs{x_0}\,\abs{\xi_0}$ we can take $\eps_0$ so small that we
will have
$$V_{(x_0,\xi_0)}\subset\left\{(x,\xi)~: \abs x\ge\dfrac 52 R_0,\quad
a(x,\xi)\le-\delta\abs
x\,\abs\xi,\quad\abs\xi\ge\dfrac{\abs{\xi_0}}2\right\}.$$
It follows from Lemma 7.5$\,(ii)$ that
$$V_{(x_0,\xi_0)}\subset\big\{(x,\xi)\in
T^*\BR^d~:~\Phi(x,\xi)=1\big\}.\leqno{(7.17)}$$
Let $b\in{\cal C}^\infty_0\big(V_{(x_0,\xi_0)}\big)$ be such that
$b\big(x_0,\xi_0\big)=1$. It follows from (7.17) that one can find $C>0$
such that
$$\abs{b(x,\xi)}\le C\Phi(x,\xi),\quad \forall(x,\xi)\in
T^*\BR^d. \leqno{(7.18)}$$
Therefore we will have $\abs{b(x,h\xi)}\le C\Phi(x,h\xi)$ for all
$(x,\xi)$ in $T^*\BR^d$ and all $h\in\big]0,1\big]$. Now on the support
of $b(x,h\xi)$ we have $h\abs\xi\ge\dfrac{\abs{\xi_0}}2$ so it follows
from Lemma 7.5$(iii)$ that $\Phi(x,h\xi)=\Phi(x,\xi)$. Therefore
$$\left\{\eqalign{
&\hbox{There exists $C>0$ such that}\hc
&\abs{b(x,h\xi)}\le C\Phi(x,\xi),\quad \forall(x,\xi)\in T^*\BR^d,\quad
\forall h\in\big]0,1\big].\hc}\right.
\leqno{(7.19)}$$
We deduce from Corollary 7.6$\,(iii)$ that
$$-H_{\tilde p}\la_1\ge\, C\scal
x^{1-\nu}\abs{b(x,h\xi)}^2\abs\xi-C',\quad \forall (x,\xi)\in
T^*\BR^d,\quad \forall h\in\big]0,1\big]
\leqno{(7.20)}$$
Let now $m_0=\big(x_0,t_0,\xi_0,\tau_0\big)$ be as in Theorem 7.1.

Let $\varphi_0\in{\cal C}^\infty_0(\Br)$, $\psi\in{\cal
C}^\infty_0(\Br)$, $\varphi_1\in{\cal C}^\infty_0\big(\Br^d\big)$ be
such that
$$\left\{\eqalign{&\varphi_0\big(t_0\big)\ne 0,\quad
\psi\big(\tau_0\big)\ne 0,\hc
&\varphi_1(x)=1\hbox{ if } \abs x\le\dfrac 43 R_0,\quad
\Supp\varphi_1\subset\left\{x~:~\abs x\le\dfrac 32 R_0\right\}.\hc}\right.
\leqno{(7.21)}$$
Let
$$b_1(x,\xi)=b(x,\xi)\abs\xi^{\sfrac 12}.\leqno{(7.22)}$$
Then $b_1\in{\cal
C}^\infty_0\left(B\big(\big(x_0,\xi_0\big),\eps_0\big)\right)$ and
$b_1\big(x_0,\xi_0\big)\ne 0$.

Finally let us recall for convenience that in (3.6) and (5.1) we have
set
$$\left\{\eqalign{
W_k(t)&=\bun_{[0,T]}\,\sconj w_k(t),\quad \sconj w_k(t)=\bun_\Omega
w_k(t),\hc
w_k(t)&=h^{-\sfrac 12}_k v_k(t),\quad
v_k(t)=\theta_1\big(h^2_kP_D\big)\widetilde u_k(t).\hc}\right.
\leqno{(7.23)}$$

\proclaim Lemma 7.7.
    We have 
$$\Int_\BR\Norme{\varphi_0(t)\psi\big(h^2_kD_t\big)b_1\big(x,h_kD_x\big)
\big(1-\varphi_1(x)\big)W_k(t)}^2_{L^2\big(\BR^d\big)}\dt=o(1)
\hbox{ as } k\to +\infty$$

\vpt
\noindent {\bf Proof}
\vpt
With $\la_1$ defined in Corollary 7.6 we set
$$N(t)=\left(\big(M-\big(1-\varphi_1\big)\la_1^w\big(1-
\varphi_1\big)\big)v_k(t),  
v_k(t)\right)_{L^2(\Omega)},$$
where $M$ is a large constant and $\la^w_1$ the Weyl quantization of the
symbol $\la_1\in S(1,g)$. Then there exists $C>0$ such that for $k\ge
1$,
$$\Norme{v_k(t)}^2_{L^2(\Omega)}\le C\,N(t).\leqno{(7.24)}$$
Setting $\Lambda=M-\big(1-\varphi_1\big)\la^w_1\big(1-\varphi_1\big)$
and $(.,.)=(.,.)_{L^2(\Omega)}$ we
can write
$$\dfrac\drm\dt N(t)=\left(\Lambda\dfrac\drm\dt v_k(t),
v_k(t)\right)+\left(\Lambda v_k(t),\dfrac\drm\dt v_k(t)\right)$$
Since $\dfrac{\drm v_k}{\dt}=i P_D v_k$ and $P_D$ is self adjoint in
$L^2(\Omega)$ we have
$$\dfrac\drm\dt N(t)=-i\left(\big[P_D,\Lambda\big]v_k(t),v_k(t)\right)$$
Now on the support of $1-\varphi_1$ we have $\abs x\ge\dfrac 43 R_0$. It
follows from (4.5)$(i)$ that $\big[P_D,\Lambda\big]=\big[\widetilde
P,\Lambda\big]$. Therefore we have,
$$\left\{\eqalign{
\dfrac\drm\dt N(t)&=\rondun+\rondtwo+\rondtrois\quad\hbox{where,}\hc
\rondun&=-i\left(\big[\widetilde
P,\varphi_1\big]\la^w_1\big(1-\varphi_1\big)v_k(t),v_k(t)\right),\hc
\rondtwo&=i\left(\big[\widetilde
P,\la^w_1\big]\big(1-\varphi_1\big)v_k(t),\big(1-\varphi_1\big)v_k(t)
\right),\hc 
\rondtrois&=-\left(\big(1-\varphi_1\big)\la_1^w\big[\widetilde
P,\varphi_1\big]v_k(t),v_k(t)\right).\hc}\right.
\leqno{(7.25)}$$
Since $\la_1=\la\Phi^2$,  Lemma 7.5 shows that the support of $\la_1$ is
contained in $\big\{x~:~\abs x\ge 2R_0\big\}$. By (7.21) the Poisson
bracket $\big(\widetilde p,\varphi_1\big\}$ has its support in
$\left\{x~:~\dfrac 43 R_0\le\abs x\le\dfrac 32 R_0\right\}$. It follows
that $\big\{\widetilde p,\varphi_1\big\}\la_1\equiv\ 0$ from which we
deduce that $\big[\widetilde P,\varphi_1\big]\la_1^w$, is a zero$^{\th}$
order operator. It follows that
$$\abs\rondun+\abs\rondtrois\le C\Norme{v_k(t)}^2_{L^2(\Omega)}.
\leqno{(7.26)}$$
Using the sharp Garding inequality, $(ii)$ in Corollary 7.6 and (7.20)
we see that

$$\rondtwo=-\left(-\big(H_{\tilde
p}\la_1\big)^w\big(1-\varphi_1\big)v_k(t),\big(1-\varphi_1\big)v_k(t)\right)+
{\cal
O}\left(\Norme{v_k(t)}^2_{L^2(\Omega)}\right)$$
$$\rondtwo\le -\Normes{\scal
x^{-\sfrac{1+\nu}2}b\big(x,h_kD_x\big)(-\Delta)^{\sfrac
14}\big(1-\varphi_1\big)v_k(t)}^2_{L^2(\BR^d)}+C\Norme{v_k(t)}^2_{L^2(\Omega)}.
\leqno{(7.27)}$$
It follows from (7.24), (7.25), (7.26), (7.27) and (7.22) that
$$N(t)+\Int_0^t\Normes{\scal
x^{-\sfrac{1+\nu}2}b_1\big(x,h_kD_x\big)\big(1-\varphi_1\big)\sconj
w_k(t)}^2_{L^2(\Br^d)}\dt\,\le\,C\Int_0^t N(s)\ds+N(0).
\leqno{(7.28)}$$
Using the Gronwall inequality we see that $N(t)\le N(0) e^{CT}$. 

Now
$N(0)\le C\Norme{v_k(0)}^2_{L^2(\Omega)}\le C'\Norme{\widetilde
u^0_k}^2_{L^2(\Omega)}$. Thus using again (7.28) we obtain
$$\Int_0^T\Norme{b_1\big(x,h_kD_x\big)\big(1-\varphi_1\big)\sconj
w_k(t)}^2_{L^2(\Omega)}\dt=0(1)
\leqno{(7.29)}$$
by (3.7), since $\scal x^{-\sfrac{1+\nu}2}\approx 1$ on the support of
$b_1\big(x,h_kD_x\big)$.

Now since $\varphi_0(t)\psi\big(h^2_kD_t\big)$ is bounded in $L^2(\Br)$
Lemma 7.7 follows from (7.29) and (7.23).

\hfill\cqfd

\vpt
\noindent {\bf End of the proof of Theorem 7.1}
\vpt
Applying Lemma 7.7 to the subsequence $\big(W_{\sigma(k)}\big)$ and
using Proposition 5.1 we see that $\scal{\mu,a}=0$ with
$$a(x,t,\xi,\tau)=\left[\big(1-\varphi_1(x)\big)
\varphi_0(t)\psi(\tau)b_1(x,\xi)\right]^2$$   
Since by (7.21), (7.22) we have $a\big(x_0,t_0,\xi_0,\tau_0\big)\ne 0$
we conclude that $m_0\notin\Supp\mu$.

The proof of Theorem~7.1 is thus
complete.\hfill\cqfd

\vpt
 \noindent {\soustitre 8.   End of the proofs of Theorem 3.1 and 2.1}
\vpt
{\bf 8.1 End of the proof of Theorem 3.1}
\vpt
According to Corollary 6.2 we will reach to a contradiction if we show
that the measure $\mu$ vanishes identically. Recall that
$$\Supp\mu \subset\Sigma=\big\{(x,t,\xi,\tau)\in T^*\BR^{d+1},\pv
x\in\conj\Omega,\pv t\in\big[0,T\big],\pv \tau+p(x,\xi)=0\big\}.$$
Let $m=(x,t,\xi,\tau)\in\Sigma$ and $a=\pi(m)\in\Sigma_b$. The
assumption (2.6) tell us that we can find $s_0\in\Br$ such that for all
$s\le s_0$ we have $\Gamma(s,a)\subset T^*M\moins\{0\}$,
$\Gamma(s,a)=\big(x(s),t,\xi(s),\tau\big)$ where $\big(x(s),\xi(s)\big)$
is the usual flow of~$p$ and $\Lim_{s\to -\infty} \abs{x(s)}=+\infty$.

Then we have the following Lemma.

\proclaim Lemma 8.1.
  One can find $s_1\le s_0$ such that with the notations of Theorem
  7.1
$$\eqalign{
&(i)\,\abs{x\big(s_1\big)}\ge 3R_0,\cr
&(ii)\,\Sum_{j,k=1}^d \widetilde a^{jk}\big(x\big(s_1\big)\big)
x_j\big(s_1\big)\xi_k\big(s_1\big)\le
-3\delta\abs{x\big(s_1\big)}\abs{\xi\big(s_1\big)}.\cr}$$

Let us assume this Lemma for a moment.

Since $\tau+\widetilde
p\big(x\big(s_1\big),\xi\big(s_1\big)\big)=\tau+p\big(x\big(s_1\big),
\xi\big(s_1\big)\big)=0$ (because $\Gamma(s,a)\subset\Sigma_b$)
we deduce from Theorem 7.1 that
$\big(x\big(s_1\big),t,\xi\big(s_1\big),\tau\big)=\Gamma\big(s_1,a\big)=
\pi^{-1}\big(\Gamma\big(s_1,a\big)\big)\notin\Supp\mu$
($\pi$ is the identity on $T^*\Br^{d+1}_M$)$\tvi$.
By Theorem 5.3 we have
$\pi^{-1}\big(\Gamma(0,a)\big)=\pi^{-1}(a)\cap\Supp\mu=\vide$.
Since $m\in\pi^{-1}(a)$ it follows that $m\notin\Supp\mu$. Therefore
$\Supp\mu=\vide$ which contradicts Corollary~6.2 and proves Theorem 3.1.

\vpt
\noindent {\bf Proof of Lemma 8.1}
\vpt
Since $\Lim_{s\to -\infty}\abs{x(s)}=+\infty$ we can find
$\widetilde s_0$ such that
$$\abs{x(s)}\ge 3R_0\pv\hbox{for}\pv s\le\widetilde s_0.
\leqno{(8.1)}$$
Let us set for $s\in\big]-\infty,\widetilde s_0\big]$
$$\left\{\eqalign{
F(s)&=F_1(s)+F_2(s),\hc
F_1(s)&=\Sum_{j,k=1}^d a^{jk}\big(x(s)\big)x_j(s)\xi_k(s),\hc
F_2(s)&=3\delta\abs{x(s)}\,\abs{\xi(s)}.\hc}\right.
\leqno{(8.2)}$$
Let us remark that since $\abs{x(s)}\ge 3R_0$ we have
$a^{jk}\big(x(s)\big)=\widetilde a^{jk}\big(x(s)\big)$.

We have
$$\left\{\eqalign{
\dot x_j(s)&=2\Sum_{k=1}^d a^{jk}\big(x(s)\big)\xi_k(s)\hc
\dot \xi_j(s)&=-\Sum_{p,q=1}^d \dfrac{\pa a^{pq}}{\pa
x_j}\big(x(s)\big)\xi_p(s)\xi_q(s)\hc}\right.
\leqno{(8.3)}$$
and by assumption (2.5), $\abs{\nabla_x a^{j}(x)}=o\left(\dfrac1{\abs
x}\right)$ as $\abs x\to +\infty$.

Using (8.1), (8.3), (2.5), the ellipticity condition (2.3) and taking
$R_0$ large enough we find by an easy computation that
$$\dfrac\drm\ds F_1(s)\ge C\abs{\xi(s)}^2,\quad
s\in\big]-\infty,\widetilde s_0\big]
\leqno{(8.4)}$$
for some fixed constant $C>0$.

Using again (8.3) and the same arguments we see easily that
$$\dfrac\drm\ds F_2(s)\le C'\delta\abs{\xi(s)}^2.
\leqno{(8.5)}$$
It follows from (8.4), (8.5) and (8.2), taking $\delta$ small enough,
that for $s\in\big]-\infty,\widetilde s_0\big]$ we have
$$\dfrac\drm\ds F(s)\ge C_0\abs{\xi(s)}^2\ge
C_1p\big(x(s),\xi(s)\big)=C_1p\left(x\big(\widetilde
s_0\big),\xi\big(\widetilde s_0\big)\right)\ge
C_2\abs{\xi\big(\widetilde s_0\big)}^2.$$
Integrating this inequality between $s$ and $\widetilde s_0$ we obtain
$$F(s)\le F\big(\widetilde s_0\big)+C_2\abs{\xi\big(\widetilde
s_0\big)}^2\big(s-\widetilde s_0\big).$$
Since the right hand side tends to $-\infty$ when $s$ goes to $-\infty$
we can find $s_1\le\widetilde s_0$ such that $F(s)\le 0$ when $s\le
s_1$.\hfill\cqfd

\vpt
{\bf 8.2 End of the proof of Theorem 2.1}
\vpt
We shall need the following Lemma.

\proclaim Lemma 8.2.
    Let $\theta\in{\cal C}^\infty_0(\Br)$, $\chi   _0\in{\cal
C}^\infty_0\big(\conj\Omega\big)$. There exists $C>0$ such that
$$\Normes{\Big[\theta\big(h^2P_D\big), \chi   _0P_D^{\sfrac
14}\Big]v}_{L^2(\Omega)}\le C\,h^{\sfrac 12}\norme v_{L^2(\Omega)}$$
for every $h\in\big]0,1\big]$ and $v\in L^2(\Omega)$.

\vpt
\noindent {\bf Proof}
\vpt
See section 9.3\hfill\cqfd

Now it is classical that one can find $\psi$, $\theta$ in ${\cal
C}^\infty_0(\Br)$ such that
$$\left\{\eqalign{
&\Supp\psi\subset\big\{t~:~\abs t\le 1\big\},\quad
\Supp\theta\subset\left\{t~:~\dfrac 12\le \abs t\le
2\right\}\quad\hbox{and}\hc
&\psi(t)+\Sum_{p=0}^{+\infty} \theta\big(2^{-p}t\big)=1\hbox{ for all }
t\in\BR.\hc}\right.$$
By the functionnel calculus we see easily that
$$\left\{\eqalign{
&\psi\big(P_D\big)+\Sum_{p=0}^{+\infty}
\theta\big(2^{-p}P_D\big)=\Id\quad\hbox{and}\hc
&\norme v^2_{L^2(\Omega)}\le
C\left(\Norme{\psi\big(P_D\big)v}^2_{L^2(\Omega)}+\Sum_{p=0}^{+\infty}
\Norme{\theta\big(2^{-p}P_D\big)v}^2_{L^2(\Omega)}\right)\hc}\right.
\leqno{(8.6)}$$
Let $u(t)=e^{itP_D}u_0$. Using (8.6) we see that
$$\left\{\eqalign{
&\Normes{\chi   _0P_D^{\sfrac 14}u(t)}^2_{L^2(\Omega)}\le
C\big(\rondun+\rondtwo\big)\quad\hbox{where}\hc
&\rondun=\Normes{\psi\big(P_D\big)\chi   _0P_D^{\sfrac
14}u(t)}^2_{L^2(\Omega)}\hc
&\rondtwo=\Sum_{p=0}^{+\infty}
\Normes{\theta\big(2^{-p}P_D\big)\chi   _0P_D^{\sfrac 14}
u(t)}^2_{L^2(\Omega)}\hc}\right.
\leqno{(8.7)}$$
We have
$$\rondun \le 2\Big(\Normes{\Big[\psi\big(P_D\big),\chi   _0P_D^{\sfrac
14}\Big]u(t)}^2_{L^2}+\Normes{\chi   _0\psi\big(P_D\big)P_D^{\sfrac 14}
u(t)}^2_{L^2}\Big)$$
Using Lemma 8.2 with $h=1$, the fact that the operator
$\psi\big(P_D\big)P_D^{\sfrac 14}$ is $L^2(\Omega)$ bounded and the
energy estimate we deduce that
$$\rondun\le C\norme{u_0}^2_{L^2(\Omega)}.\leqno{(8.8)}$$
On the other hand we have
$$\left\{\eqalign{
\rondtwo&\le C\big(\rondtrois+\rondq\big)\hc
\rondtrois&=\Sum_{p=0}^{+\infty}
\Normes{\Big[\theta\big(2^{-p}P_D\big),\chi   _0P_D^{\sfrac
14}\Big]u(t)}^2_{L^2(\Omega)}\hc
\rondq&=\Sum_{p=0}^{+\infty}\Normes{\chi   _0P_D^{\sfrac
14}\theta\big(2^{-p}P_D\big)u(t)}^2_{L^2(\Omega)}\hc}\right.
\leqno{(8.9)}$$
Using Lemma 8.2 we can write
$$\rondtrois\le C\left(\Sum_{p=0}^{+\infty}
2^{-p}\right)\norme{u(t)}^2_{L^2(\Omega)},$$
so by the energy estimate
$$\rondtrois \le C\Norme{u_0}^2_{L^2(\Omega)}.
\leqno{(8.10)}$$
To handle the term \rondq\ we use the Theorem 3.1. Let
$\widetilde\theta\in{\cal C}^\infty_0(\Br)$ be such that\hfill\break
$\Supp\widetilde\theta\subset\left\{t~:~\dfrac 13\le\abs t\le
3\right\}$, $\widetilde \theta(t)=1$ on the support of $\theta$ and
$0\le\widetilde\theta\le 1$. Then
$\theta\big(2^{-p}P_D\big)u(t)=\theta\big(2^{-p}P_D\big)\erom^{itP_D}
\widetilde\theta\big(2^{-p}P_D\big)u_0$.

It follows from Theorem 3.1 that
$$\Int_0^T\rondq\, \dt\,\le\,
C\Sum_{p=0}^{+\infty}\Norme{\widetilde\theta\big(2^{-p}P_D\big)
u_0}^2_{L^2(\Omega)}$$ 
Now we have
$$\Sum_{p=0}^{+\infty}
\left[\widetilde\theta\big(2^{-p}t\big)\right]^2\le\left(\Sum_{p=0}^{+\infty}
\widetilde\theta\big(2^{-p}t\big)\right)^2\le
M_0,\quad\hbox{for all } t\in\BR$$
It follows that the operator $\Sum_{p=0}^{+\infty}
\left[\widetilde\theta\big(2^{-p}P_D\big)\right]^2$ is $L^2(\Omega)$
bounded, therefore
$$\Int_0^T\rondq\, \dt\,\le\,C\Norme{u_0}^2_{L^2(\Omega)}.\leqno{(8.11)}$$
It follows from (8.7), (8.8), (8.9), (8.10) and (8.11) that
$$\Int_0^T\Normes{\chi   _0P_D^{\sfrac 14}\erom^{itP_D}
u_0}^2_{L^2(\Omega)}\dt\,\le\,C\Norme{u_0}^2_{L^2(\Omega)}$$
which is the claim in Theorem 2.1. The proof is complete.\hfill\cqfd

\vpt
\noindent {\soustitre 9.   Appendix}
\vpt
{\bf 9.1 The geometrical framework}
\vpt
We recall here the definition of the generalized bicharacteristic flow
in the sense of Melrose and Sj{\"o}strand. For this purpose we follow
  {H{\"o}rmander} [H{\"o}].

Let $M=\Omega\times\BR_t$. We set $T^*_bM=T^*M\moins \{0\}\cup T^*\pa
M\moins\{0\}$. We have a natural restriction %\break
${\rm map\,}\pi~:~T^*\BR^{d+1}_{|\conj M}\to T^*_bM$ (which will be
describe more 
precisely in local coordinates below) which is the identity on
$T^*\BR^{d+1}_{|M}\moins\{0\}$.

With $p$ defined in (2.2) we introduce the characteristic set
$$\Sigma=\left\{(x,t,\xi,\tau)\in T^*\BR^{d+1},\pv x\in\conj\Omega,\pv
t\in\big[0,T\big],\pv \tau+p(x,\xi)=0\right\},$$
and we set $\Sigma_b=\pi(\Sigma)$.

\proclaim Definition 9.1.
    Let $\zeta\in T^*\pa M\moins\{0\}$. We shall say that
$$\eqalign{
&(i)\,\zeta\hbox{ is elliptic (or } \zeta\in{\cal E})\hbox{ iff }
\zeta\notin\Sigma_b\cr
&(ii)\,\zeta\hbox{ is hyperbolic (or }\zeta\in {\cal H})\hbox{ iff }
\#\left\{\pi^{-1}(\zeta)\cap\Sigma\right\}=2\cr
&(iii)\,\zeta \hbox{ is glancing (or } \zeta\in{\cal G}\hbox{ iff }
\#\left\{\pi^{-1}(\zeta)\cap\Sigma\right\}=1.\cr}$$

Let us describe $\pi$ and these sets in local coordinates. As we said
before, $\pi$ is the identity map on $T^*\BR^{d+1}_{|M}\moins\{0\}$.

Near any point of $\pa M$ we can use the geodesical coordinates where
$M$ is given by
$$\left\{\big(x_1,x',t\big)\in\BR\times\BR^{d-1}\times\BR~:~x_1>0\right\},$$
$\pa M$ is given by $\left\{\big(x_1,x',t\big)~:~x_1=0\right\}$ and
$\tau+p(x,\xi)$ is transformed to
$\xi^2_1+r\big(x_1,x',\xi'\big)+\tau$.

In these coordinates if $\rho\in T^*\BR^{d+1}_{|\pa M}\moins\{0\}$ then
$\rho=\big(0,x',t,\xi_1,\xi',\tau\big)$ and
$\pi(\rho)=\big(x',t,\xi',\tau\big)\in T^*\pa M\moins\{0\}$.

Now let $\zeta=\big(x',t,\xi',\tau\big)\in T^*\pa M\moins\{0\}$. Then
$$\left\{\eqalign{
\zeta\in{\cal E}&\iff r\big(0,x',\xi'\big)+\tau>0,\hc
\zeta\in{\cal H}&\iff r\big(0,x',\xi'\big)+\tau<0,\hc
\zeta\in{\cal G}&\iff r\big(0,x',\xi'\big)+\tau=0.\hc}\right.
\leqno{(9.1)}$$
When $\zeta\in{\cal H}$ then
$\pi^{-1}(\zeta)\cap\Sigma=\left\{\big(0,x',t,\xi_1^\pm,\xi',\tau\big)\right\}$
where
$$\xi_1^\pm=\pm\bigg(-\left(r\big(0,x',\xi'\big)+\tau\right)\bigg)^{\sfrac
12}.
\leqno{(9.2)}$$
When $\zeta\in{\cal G}$ then
$\pi^{-1}(\zeta)\cap\Sigma=\big\{\big(0,x',t,0,\xi',\tau\big)\big\}$.

For the purpose of the proofs it is important to decompose the set
${\cal G}$ of glancing points into several subsets. The following
definition is given in local coordinates but could be written in an
intrinsic way (see [H]). We shall set
$$r_0\big(x',\xi'\big)=r\big(0,x',\xi'\big)
\leqno{(9.3)}$$
and $H_{r_0}$ will denote the Hamilton field of $r_0$ namely
$\displaystyle H_{r_0}=\dfrac{\pa r_0}{\pa \xi'}\dfrac{\pa}{\pa
x'}-\dfrac{\pa r_0}{\pa x'}\dfrac\pa{\pa \xi'}\cdotp$

\proclaim Definition 9.2.
    Let $\zeta=\big(x',t,\xi',\tau\big)\in{\cal G}$. We shall say that
$$\eqalign{
&(i)\,\zeta\hbox{ is diffractive (or } \zeta\in{\cal G}_d)\hbox{ iff }
\dfrac{\pa r}{\pa x_1}\big(0,x',\xi'\big)<0,\cr
&(ii)\,\zeta\hbox{ is gliding (or }\zeta\in{\cal G}_g)\hbox{ iff }
\dfrac{\pa r}{\pa x_1}\big(0,x',\xi'\big)>0,\ 
\hbox{ and we set } {\cal G}^2={\cal G}_d\cup{\cal G}_g,\cr
&(iii)\,\zeta\hbox{ belongs to } {\cal G}^k,\  k\ge 3\hbox{ , iff }\cr
&\quad\quad\,H^j_{r_0}\left(\dfrac{\pa r}{\pa
x_1}_{|x_1=0}\right)(\zeta)=0,\ 0\le j<k-2,\
H^{k-2}_{r_0}\left(\dfrac{\pa r}{\pa x_1}_{|x_1=0}\right)(\zeta)\ne
0.\cr}$$

We can now give the meaning of the assumption made in (2.7).

\proclaim Definition 9.3.
   We shall say that the bicharacteristics have no contact of
infinite order with the boundary if
$${\cal G}=\Cup_{k=2}^{+\infty} {\cal G}^k.$$

We are going now to make a brief description of the generalized
bicharacteristic flow and we refer to [M-S] or [H{\"o}] for more details.

First of all we introduce some notations.

We shall denote by $\ga(s)=\big(x(s),\xi(s)\big)$ the usual
bicharacteristic of $p$ in $T^*\Omega$ defined by
$$\big(\dot
x(s),\dot \xi(s)\big)=\Big(\dfrac{\pa
p}{\pa\xi}\big(\ga(s)\big),-\dfrac{\pa p}{\pa
x}\big(\ga(s)\big)\Big).$$
We shall denote by $\ga_g(s)=\left(x'_g(x),\xi'_g(s)\right)$ the gliding
ray in $T^*\pa\Omega$ defined in the geodesic coordinates by the
equations
$$\left({\dot
x}\,'_g(s),{\dot\xi}\,'_g(s)\right)=\Big(\dfrac{\pa r_0}{\pa
\xi'}\big(\ga_g(s)\big),-\dfrac{\pa  r_0}{\pa
x'}\big(\ga_g(s)\big)\Big)$$
where $r_0$ has been introduced in (9.3).

\noindent The generalized flow lives in $\Sigma_b\subset T^*_bM$ and for
$a\in\Sigma_b$ is denoted by $\Gamma(s,a)$. Since $\Sigma_b$ is the
disjoint union of $\Sigma_b\cap T^*M$, $\Sigma_b\cap{\cal H}$,
$\Sigma_b\cap{\cal G}_d$, $\Sigma_b\cap{\cal G}_g$ and
$\Sigma_b\cap\Bigg(\Cup_{k\ge 3}{\cal G}^k\Bigg)$ we shall consider
separatly the case where a belongs to each set. Moreover each
description of $\Gamma(s,a)$ holds for small $\abs s$.

\vpt
{\bf Case 1}~: $a\in\Sigma_b\cap T^*M$
\vpt
Here $a=(x,t,\xi,\tau)$ where $x\in\Omega$, $t\in\big[0,T\big]$,
$\tau+p(x,\xi)=0$. Then for small $\abs s$ we have
$$\Gamma(s,a)=\big(x(s),t,\xi(s),\tau\big)\subset T^*M$$
where $\big(x(s),\xi(s)\big)$ is the bicharacteristic of $p$ starting
from the point $(x,\xi)$.

\vpd
{\bf Case 2}~: $a\in\Sigma_b\cap{\cal H}$
\vpt
\noindent In the geodesic coordinates we have
$a=\big(x',t,\xi',\tau\big)$ and the 
equation $\xi^2_1+r\big(0,x',\xi'\big)+\tau=0$ has two distinct roots
$\xi^+_1$, $\xi^-_1$ described in (9.2). For $s>0$ (resp. $s<0$) let
$\ga^+(s)=\big(x^+(s),\xi^+(s)\big)$ (resp.
$\ga^-(s)=\big(x^-(s),\xi^-(s)\big)$ be the bicharacteristic of $p$
starting for $s=0$ at the point $\big(0,x',\xi^+_1,\xi'\big)$ (resp.
$\big(0,x',\xi^-_1,\xi'\big)$). They are contained in $T^*\Omega$ for
small $\abs s\ne 0$. Then $\Gamma(0,a)=a$ and
$$\Gamma(s,a)=\left\{\eqalign{
&\big(x^+(s),t,\xi^+(s),\tau\big),\quad 0<s<\eps,\hc
&\big(x^-(s),t,\xi^-(s),\tau\big),\quad -\eps<s<0.\hc}\right.$$

Here $\Gamma(s,a)\subset T^*M$ for $s\ne 0$.

\vpd
{\bf Case 3}~: $a\in\Sigma_b\cap {\cal G}_d$
\vpt
Here $a=\big(x',t,\xi',\tau\big)$ and the equation
$\xi^2_1+r\big(0,x',\xi'\big)+\tau=0$ has a double root $\xi_1=0$.
Let $\ga(s)=\big(x(s),\xi(s)\big)$ be the flow of $p$ starting when
$s=0$ at the point $\big(0,x',\xi_1=0,\xi'\big)$. Then we have
$$\Gamma(s,a)=\big(x(s),t,\xi(s),\tau\big)\subset T^*M,\quad 0<\abs
s<\eps.$$

\vpd
{\bf Case 4}~: $a\in\Sigma_b\cap{\cal G}_g$
\vpt
As above $a=\big(x',t,\xi',\tau\big)$ and $\xi_1=0$ is a double root.
Let $\ga_g(s)=\big(x'_g(s),\xi'_g(s)\big)$ be the gliding ray starting
when $s=0$ at the point $\big(x',\xi'\big)$. Then  we have
$$\Gamma(s,a)=\left(x'_g(s),t,\xi'_g(s),\tau\right)\subset T^*\pa
M,\quad \abs s<\eps.$$

\vpd
{\bf Case 5}~: $a\in \Sigma _b \cap \left(\Cup_{k=3}^{+\infty}
{\cal G}^k\right)$
\vpt
Let $a=\big(x',t,\xi',\tau\big)$. Let
$\ga_g(s)=\big(x'_g(s),\xi'_g(s)\big)$ be the gliding ray starting when
$s=0$ at the point $\big(x',\xi'\big)$. Then (see Theorem 24.3.9 in
[H{\"o}]) one can find $\eps>0$ such that with $I=\big]0,\eps\big[$ we have
either $\ga_g(s)\in{\cal G}_g$, $\forall s\in I$ and then
$\Gamma(s,a)=\left(x'_g(s), t, \xi'_g(s),\tau\right)\subset T^*\pa M$,
$\forall s\in I$, or $\ga_g(s)\in {\cal G}_d$, $\forall s\in I$ and
then
$\Gamma(s,a)=\big(x(s),t,\xi(s),\tau\big)\subset T^*M$, $\forall
s\in I$, where $\big(x(s),\xi(s)\big)$ is the bicharacteristic of $p$
starting when $s=0$ at the point $\big(0,x',\xi_1=0,\xi'\big)$.

The same discussion is independently valid for $-\eps<s<0$.

\proclaim Remark 9.4.
 Let $a\in\Sigma_b$ and $\Gamma(t,a)$ be the generalized
bicharacteristing starting for $t=0$ at the point $a$. Then the above
discussion shows that one can find $\eps>0$ such that for $0<\abs
t\le\eps$ we have $\Gamma(t,a)\subset T^*M\cup{\cal G}_g$.
  Let us note (see [M-S]) that the maps $s\mapsto\Gamma(s,a)$ and
$a\mapsto\Gamma(s,a)$ are continuous, the later when $T^*_bM$ is endowed
with the topologie induced by the projection $\pi$. Moreover we have the
usual relation $\Gamma(t+s,a)=\Gamma\big(t,\Gamma(s,a)\big)$ for $s,t$
in
$\BR$.

\vpt
{\bf 9.2 Proofs of Theorem 5.2 and Theorem 5.3}
\vpt
{\bf a) Proof of Theorem 5.2}
\medskip
According to (5.1), (5.2) it is obvious that
$$\Supp\mu\subset\left\{(x,t,\xi,\tau)\in
T^*\Br^{d+1}~:~x\in\conj\Omega\hbox{ and } t\in\big[0,T\big]\right\}$$
Therefore it remains to show that if
$m_0=\big(x_0,t_0,\xi_0,\tau_0\big)$ with $x_0\in\conj\Omega$,
$t_0\in\big[0,T\big]$, but $\tau_0+p\big(x_0,\xi_0\big)\ne 0$ then
$m_0\notin\Supp\mu$.

\vpd
{\bf Case 1}~: assume $x_0\in\Omega$
\medskip
Let $\varepsilon>0$ be such that $B\big(x_0,\eps\big)\subset\Omega$. Let
$\varphi\in{\cal C}^\infty_0\big(B\big(x_0,\eps\big)\big)$, $\varphi=1$
on $B\left(x_0,\dfrac\eps2\right)$ and $\widetilde\varphi\in{\cal
C}^\infty_0(\Omega)$, $\widetilde\varphi=1$ on $\Supp\varphi$. Let
$a\in{\cal C}^\infty_0\left(\Br^d_x\times\BR^d_\xi\right)$ such that
$\pi_x\Supp a\subset B\left(x_0,\dfrac\eps2\right)$ and $\chi   \in{\cal
C}^\infty_0\left(\BR_t\times\BR_\tau\right)$. Recall that we have set
$W_k=\bun_{[0,T]}\bun_\Omega w_k$ with $w_k=h^{-\sfrac
12}_k\theta\big(h^2_kP_D\big)\widetilde u_k$ and that $\big(w_k\big)$ is
a bounded sequence in $L^2\left(\big[0,\tau\big],
L^2_{\loc}\big(\Br^d\big)\right)$ (see Proposition 4.1). Now we set
$$I_k=\left(a\big(x,h_kD_x\big)\chi   \big(t,h^2_kD_t\big)\varphi
h^2_k\big(D_t+P\big(x,D_x\big)\big)W_k, \widetilde \varphi
W_k\right)_{L^2(\BR^{d+1})}
\leqno{(9.4)}$$
We have
$h^2_k\left(D_t+P\big(x,D_x\big)\right)=h^2_kD_t+P_2
\big(x,h_kD_x\big)+h_kP_1\big(x,h_kD_x\big)+h^2_k 
P_0(X)$ where $P_j(x,\xi)$ are homogeneous in $\xi$ of order $j$.

Using the semi classical symbolic calculus and the fact that
$\big(\widetilde\varphi W_k\big)$ is bounded in $L^2\big(\Br^{d+1}\big)$
we see easily that the terms in $I_k$ corresponding to
$h_kP_1\big(x,h_kD_x\big)$ and $h^2_kP_0(x)$ tend to zero when $k\to
+\infty$. It remains to consider the term $P_2\big(x,h_kD_x\big)$. But
by the semi classical calculus
$$a\big(x,h_kD_x\big)\chi   \big(t,h^2_kD_t\big)\varphi
\big(h^2_kD_t+P_2\big(x,h_kD_x\big)\big)={\cal O}p\big(a\chi
(\tau+p)\big)+h_kR_k$$ 
where $R_k$ is a uniformly bounded semi classical pseudo-differential
operator in $L^2\big(\BR^{d+1}\big)$. Therefore the term in $I_k$
corresponding to $h_kR_k$ tends to zero.

It follows from Proposition 5.1 that
$$\Lim_{k\to +\infty}
I_{\sigma(k)}=\scal{\mu,(\tau+p)a\chi   }\leqno{(9.5)}$$
On the other hand we have since
$\big[D_t+P\big(x,D_x\big)\big]\widetilde
u_k=0$ in $\Omega\times\BR_t$ and $\varphi\in{\cal C}^\infty_0(\Omega)$,
$$\varphi\big(D_t+P\big(x,D_x\big)\big)W_k=\varphi\big(w_k(0)
\delta_{t=0}-w_k(T)\delta_{t=T}\big)\leqno{(9.6)}$$ 
\proclaim Lemma 9.5.
    Let $1\le p\le +\infty$, $\chi   \in{\cal C}^\infty_0(\BR\times\BR)$ and
$\ell\ge 1$. Then there exists $C>0$ such that
$$\Norme{\chi   \big(t,h^\ell D_t\big)\delta_{t=a}}_{L^p(\Br)}\le
C\,h^{\sfrac\ell p-\ell}$$
for every $0<h\le 1$.

\medskip
\noindent {\bf Proof}
\medskip
Let $\psi\in{\cal C}^\infty_0(\Br)$, $\psi(t)=1$ on a neighborhood of
$\pi_t\Supp\chi   $. Then $\psi\chi   =\chi   $. Now
$$\rondun=\chi   \big(t,h^\ell
D_t\big)\delta_{t=a}=\dfrac1{2\pi}\Int\erom^{i(t-a)\tau}\chi
\big(t,h^2\tau\big)\dtau$$ 
On the other hand
$\chi   (t,\tau)=\psi(t)\chi   (a,\tau)+\psi(t)(t-a)\widetilde\chi
(t,\tau,a)$ 
where $\widetilde\chi   \in{\cal C}^\infty$ has compact support in $\tau$.
It follows that
$$\left\{\eqalign{
\rondun&=\psi(t)h^{-\ell}\left(\conj{{\cal
F}}_x\chi
\right)\left(a,\dfrac{t-a}{h^\ell}\right)+\rondtwo\quad\hbox{with}\hc 
\rondtwo&=\dfrac 1{2\pi}\psi(t)\Int(t-a)\erom^{i(t-a)\tau}\widetilde
\chi   \big(t,h^\ell\tau,a\big)\dtau\hc}\right.
\leqno{(9.7)}$$
Noting that $(t-a)\erom^{i(t-a)\tau}=\dfrac 1i\dfrac{\pa}{\pa
\tau}\erom^{i(t-a)\tau}$ and making an integration by part, we see easily
that
$$\abs{\rondtwo}\le
C\abs{\psi(t)}h^\ell\Int\abs{\dfrac{\pa\widetilde\chi   }{\pa\tau}
\big(t,h^\ell\tau,a\big)}\dtau=C\abs{\psi(t)}\Int\abs{\dfrac{\pa\widetilde\chi
}{\pa\tau}(t,\tau,a)}\dtau 
\leqno{(9.8)}$$
Then the Lemma follows easily from (9.7) and (9.8).\hfill\cqfd

Now we see from (9.6) that $I_k$ is a sum of two terms of the form
$$J_k=\left(a\big(x,h_kD_x\big)\varphi w_k
(a)h^2_k\chi   \big(t,h^2_kD_t\big)\delta_{t=a},\widetilde\varphi
W_k\right),\quad a=0\hbox{ or } T$$
Since $\big(\widetilde\varphi W_k\big)$ is bounded in
$L^2\big(\BR^{d+1}\big)$ we see that
$$\abs{J_k}^2\le C\Norme{a\big(x,h_kD_x\big)\varphi
w_k(a)}^2_{L^2(\Br^d)}\Norme{h^2_k\chi   \big(t,h^2_kD_C\big)
\delta_{t=a}}^2_{L^2(\Br)}$$ 
so using Lemma 9.5 with $p=2$ and $\ell=2$ we deduce that
$$\abs{J_k}^2\le C\,h^{-1}_k \Norme{\widetilde
u_k(a)}^2_{L^2(\Omega)}h^2_k\le C\,h_k\Norme{\widetilde
u^0_k}^2_{L^2(\Omega)}$$
by the energy estimate. It follows from (3.7) that
$$\Lim_{k\to +\infty} I_k=0\leqno{(9.9)}$$
Using (9.5) and (9.9) we see that $\scal{\mu,(\tau+p)a\chi   }=0$. Since
$\tau_0+p\big(x_0,\xi_0\big)\ne 0$ and
$${\cal
C}^\infty_0\big(\Br^d_x\times\Br^d_\xi\big)\otimes{\cal
C}^\infty_0\big(\Br_t\times\Br_\tau\big)$$
is dense in ${\cal
C}^\infty_0\big(\Br^{d+1}\times\BR^{d+1}\big)$ we deduce that
$m_0=\big(x_0,t_0,\xi_0,\tau_0\big)\notin\Supp\mu$.

\vpt
{\bf Case 2}~: assume $x_0\in\pa\Omega$
\vpt
We would like to show that one can find a neighborhood $U_{x_0}$ of
$x_0$ in $\Br^d$ such that for any\break $a\in{\cal
C}^\infty_0\big(U_{x_0}\times\Br_t\times\Br^d_\xi\times\Br_\tau\big)$ we
have
$$\scal{\mu,(\tau+p)a}=0\leqno{(9.10)}$$
Indeed this will imply that the point
$m_0=\big(x_0,t_0,\xi_0,\tau_0\big)$ (with
$\tau_0+\big(x_0,\xi_0\big)\ne 0$) does not belong to the support of
$\mu$ as claimed.
%\vskip 5mm
  %\newpage

Now (9.10) will be implied, according to Proposition 5.1 and (4.1) by
$$\left\{\eqalign{
&\Lim_{k\to +\infty} I_{\sigma(k)}=0\quad\hbox{where}\hc
&I_k=\left(a\big(x,t,h_kD_x,h^2_kD_t\big)\varphi
h^2_k\big(D_t+P\big)W_k,W_k\right)_{L^2(\BR^{d+1})}\hc}\right.
\leqno{(9.11)}$$
where $\varphi\in{\cal C}^\infty_0\big(V_{x_0}\big)$, $\varphi=1$ on
$\pi_x\Supp a$.

Now we may choose $U_{x_0}$ so small that one can find a ${\cal
C}^\infty$ diffeomorphism $F$ from $U_{x_0}$ to a neighborhood $U_0$ of
the origin in $\Br^d$ such that
$$\left\{\eqalign{
F\big(U_{x_0}\cap\Omega\big)&=\big\{y\in U_0~:~y_1>0\big\}\hc
F\big(U_{x_0}\cap\pa\Omega\big)&=\big\{y\in U_0~:~y_1=0\big\}\hc
\big(P(x,D)W_k\big)\circ
F^{-1}&=\left(D^2_1+R\big(y,D'\big)\right)\left(W_k\circ
F^{-1}\right)\hc}\right.
\leqno{(9.12)}$$
where $R$ is a second order differential operator and
$D'=\big(D_2,\cdots,D_d\big)$.

Let us set
$$v_k=w_k\circ F^{-1},\quad
V_k=\bun_{[0,T]}\bun_{y_1>0}v_k\leqno{(9.13)}$$
then we will have
$$\left\{\eqalign{
&\left(D_t+D^2_1+R\big(y,D'\big)\right)v_k=0\hbox{ in }
U_0\times\BR_t,\quad y_1>0,\hc
&v_k|_{y_1=0}=0\hc}\right.
\leqno{(9.14)}$$
Making the change of variable $x=F^{-1}(y)$ in the right hand side of
the second line of (9.11) we see that
$$I_k=\left(b\big(y,t,h_tD_y, h^2_kD_t\big)\psi
h^2_k\big(D_t+D^2_1+R\big(y,D'\big)\big)V_k,V_k\right)_{L^2(\Br^{d+1})}$$
where $b\in {\cal
C}^\infty_0\big(U_0\times\BR_t\times\Br_\eta\times\Br_\tau\big)$ and
$\psi\in{\cal C}^\infty_0\big(U_0\big)$, $\psi=1$ near $\pi_y\Supp b$.

To prove (9.11) it is sufficient to prove that
$$
\Lim J_k=\Lim \left( T\psi_0\big(y_1\big)\psi_1\big(y'\big)h^2_k
\big(D_t+D^2_1+R\big(y,D'\big)\big)V_k,V_k\right)_{L^2(\Br^{d+1})}=0
\leqno{(9.15)}$$
where  $T=\theta\big(y_1,h_kD_1\big)\Phi\big(y',h_kD'\big)\chi   
\big(t,h^2_kD_t\big)$,
$\theta\Phi\chi   \in{\cal
C}^\infty_0\big(U_0\times\Br_t\times\Br^d_y\times\Br_\tau\big)$, 
$\psi_0\psi_1\in{\cal C}^\infty_0\big(U_0\big)$ and $\psi_0\psi_1=1$ on
$\pi_y\Supp\theta\Phi\xi$.

Now according to (9.14) we have
$$\eqalign{
\big(D_t+ & D^2_1+R\big(y,D'\big)\big)V_k\cr
&=-i\bun_{x_1>0}v_k(0,
\cdot)\delta_{t=0} 
+i\bun_{X_1>0}v_k(T,\cdot)\delta_{t=T}
-i\bun_{[0,T]}\big(D_1v_k|_{x_1=0}\big)\otimes\delta_{x_1=0}\cr}
\leqno{(9.16)}$$  
Therefore (9.11) will be proved in we can prove that
$$\left\{\eqalign{&
\Lim_{k\to +\infty}A^j_k=0,\quad j=1,2,\quad\hbox{where}\hc
A^1_k&=\left(\theta\big(y,h_kD_1\big)\Phi\big(y',h_kD'\big)\chi
\big(t,h^2_kD_t\big)\psi_0\psi_1h^2_k\bun_{y_1>0}v_k(a,\cdot)
\delta_{t=a},V_k\right),  
a=0,T\hc
A^2_k&=\left(\theta\big(y_1,h_kD_1\big)\Phi\big(y',h_kD'\big)\chi
\big(t,h^2_kD_t\big)\psi_0\psi_1h^2_k 
\bun_{[0,T]}\big(D_1v_k|_{y_1=0}\big)\otimes\delta_{y_1=0},V_k
\right)\hc}\right.  
\leqno{(9.17)}$$
Since  the operator $\theta\big(y_1,h_kD_1\big)\Phi\big(y',h_kD'\big)$
is uniformly bounded in $L^2\big(\Br^d\big)$ we can write with
$\psi_2\in{\cal C}^\infty_0\big(U_0\big)$, $\psi_2=1$ near
$\pi_y\Supp\theta\Phi$
$$\abs{A^1_k}^2\le
C\Norme{h^2_k\chi   \big(t,h^2_kD_t\big)\delta_{t=a}}^2_{L^2(\Br)}
\Norme{\psi_0\psi_1\bun_{y_1>0}v_k(a,\cdot)}^2_{L^2(\Br^d)}
\Norme{\psi_2V_k}^2_{L^2(\Br^{d+1})}$$
By (9.13), the energy estimate and Proposition 4.1 we have
$$\left\{\eqalign{
&\Norme{\psi_0\psi_1\bun_{y_1>0} v_k(a,\cdot)}^2_{L^2}\le
C\Norme{w_k(a)}^2_{L^2(\Omega)}\le C\,h^{-1}_k\Norme{\widetilde
u_k(a)}^2_{L^2}\le C\,h^{-1}_k\Norme{\widetilde u\,^0_k}^2_{L^2}\hc
&\Norme{\psi_2V_k}^2_{L^2}\le C\Int_0^T\Norme{\big(\psi_2\circ
F\big)w_k(t,\cdot)}^2_{L^2(\Omega)}\dt={\cal O}(1)\hc}\right.$$
Using Lemma 9.5 with $\ell=2$, $p=2$, we obtain
$$\abs{A^1_k}^2\le C\,h_k\Norme{\widetilde
u^0_k}^2_{L^2(\Omega)}\leqno{(9.18)}$$
To estimate the term $A^2_k$ we need a Lemma.

With $U_0$ introduced in (9.12) we set $U^+_0=\big\{y\in
U_0~:~y_1>0\big\}$. We shall consider smooth solution of the problem
$$\left\{\eqalign{
&\left(D_t+D^2_1+R(y,D'\big)\right)u=0\quad\hbox{in }
U^+_0\times\Br_t\hc
&u|_{y_1=0}=0\hc}\right.
\leqno{(9.19)}$$
\proclaim Lemma 9.6. 
    Let $\chi   \in{\cal C}^\infty_0\big(U_0\big)$ and $\chi   _1\in{\cal
C}^\infty_0\big(U_0\big)$ on $\Supp\chi   $. There exists $C>0$ such that
for any solution $u$ of (9.19) and all $h$ in $\big]0,1\big]$ we have
 $$\displaylines{
\Int_0^T\Norme{\big(\chi    h\pa_1u\big)_{|y_1=0}(t)}^2_{L^2}\dt\,\le\,
C\Biggl(\Int_0^T\Sum_{\abs\alpha\le 1}\Norme{\chi   _1(hD)^\alpha
u(t)}^2_{L^2(U^+_0)}
+\Norme{h^{\sfrac 12}\chi    u(0)}_{L^2(U^+_0)}\hfill\cr
\hfill\Norme{h^{\sfrac
12}\big(h\pa_1u\big)(0)}_{L^2(U^+_0)}+\Norme{h^{\sfrac
12}u(T)}_{L^2(U^+_0)}
\Norme{h^{\sfrac 12}\big(h\pa_1u\big)(T)}_{L^2(U^0_0)} \Biggl) \cr}$$

\proclaim Corollary 9.7.
    One can find a constant $C>0$ such that
$$\Int_0^T\Norme{(\chi    h_k\pa _1v_k\big)_{|y_1=0}(t)}^2_{L^2}\dt\,\le\,
C\left(\Int_0^T\Norme{\widetilde
\chi w_k(t)}^2_{L^2(\Omega)}\dt+\Norme{\widetilde
u^0_k}^2_{L^2(\Omega)}\right)={\cal O}(1)$$
 where $v_k$ has been defined in (9.13) and $\widetilde\chi   \in{\cal
C}^\infty_0\big(\Br^d\big)$.

\vpt
\noindent {\bf Proof of the Corollary}
\vpt
We use Lemma 9.6, (9.13), (9.14) the fact that $w_k=h^{-\sfrac
12}_k\theta\big(h^2_kP_D\big)\widetilde u_k$
Lemma 6.3$(ii)$ and the energy estimate for $\widetilde u_k$.\hfill\cqfd

\vpt
\noindent {\bf Proof of Lemma 9.6}
\vpt
Let us set with $L^2=L^2\left(\BR^+_{y_1}\times\BR^{d-1}_{y'}\right)$
$$\left\{\eqalign{
\hbox{I}&=\Sum_{\abs\alpha\le 1}\Norme{\chi   _1(hD)^\alpha
u(t)}^2_{L^2}\hc
\hbox{II}&=\Sum_{j=0}^1\Norme{h^{\sfrac 12}\chi   
u\big(a_j\big)}_{L^2}\cdot\Norme{h^{\sfrac
12}\chi   \big(h\pa_1u\big)\big(a_j\big)}_{L^2},\quad a_0=0,
a_1=T\hc}\right.
\leqno{(9.20)}$$
By (9.19) we have
$$\eqalign{2\Re\Int_0^T\left(\chi   
h\left(D^2_1u(t)+R\big(y,D'\big)u(t)\right),\chi   
h\pa_1u(t)\right)_{L^2}&\dt\hfill\cr
=-2\Im&\Int_0^T \big(\chi    h\pa_tu(t),\chi    h\pa _1
u(t)\big)_{L^2}\dt\cr}\leqno{(9.21)}$$
By integration by part we have
$$\displaylines{\Int_0^T\big(\chi    h\pa_tu(t),\chi    h\pa_1
u(t)\big)_{L^2}\dt-\Int_0^T\big(\chi    h\pa _1 u(t),\chi   h\pa
_tu(t)\big)_{L^2}\hfill\cr
\hfill=
\Int_0^T\left(\big(\pa_1\chi   ^2\big)h\,u(t),h\pa_tu(t)\right)\dt+{\cal
O}(\hbox{II})\cr}$$
Since $h\pa_tu(t)=-i hD^2_1u(t)-ihR\big(y,D'\big)u(t)$ integrating by
part and using the fact that $u_{|y_1=0}=0$ we find that
$$\Int_0^T\left(\big(\pa_1\chi
^2\big)h\,u(t),h\pa_tu(t)\right)_{L^2}\dt={\cal 
O}\hbox{(I)}$$
It follows that
$$\Im\Int_0^T\big(\chi    h\pa_t u(t),\chi    h\pa_1 u(t)\big)_{L^2}\dt={\cal
O}(\hbox{I}+\hbox{II})
\leqno{(9.22)}$$
Now
$$\displaylines{-\Int_0^T\left(\chi    h\pa^2_1u(t),\chi    h\pa _1
u(t)\right)\dt \hfill\cr
\hfill=\Int_0^T\Norme{\big(\chi    h\pa_1
u(t)\big)_{|y_1=0}}^2_{L^2(\Br^{d-1})}+{\cal O}(\hbox{I})
-\Int_0^T\left(\chi    h\pa_1 u(t),\chi    h\pa^2_1 u(t)\right)_{L^2}\dt\cr}$$
from which we deduce that
$$-2\Re\Int_0^T\left(\chi    h\pa^2_1 u(t),\chi    h\pa_1
u(t)\right)_{L^2}\dt=
\Int_0^T\Norme{\big(\chi    h\pa_1
u(t)\big)_{|y_1=0}}^2_{L^2(\BR^{d-1})}\dt+{\cal O}(\hbox{I})
\leqno{(9.23)}$$
Finally using again integration by parts, the fact that $R$ is
symmetric and $D'u_{|y_1=0}=0$ we find that
$$2\Re\Int_0^T\left(\chi    hR\big(y,D'\big)u(t),\chi    h\pa_1
u(t)\right)_{L^2}\dt={\cal O}(\hbox{I})
\leqno{(9.24)}$$
Then the Lemma follows from (9.20) to (9.24).\hfill\cqfd

Let us go back to the estimate of $A^2_k$ defined in (9.17). We have
$$\abs{A^2_k}^2\le C\,h^2_k
\Norme{\theta\big(y_1,h_xD_1\big)\delta_{y_1=0}}^2_{L^2(\Br)}
\Int_0^T\Norme{\psi_1\big(h D_1v_k(t)\big)_{|y_1=0}}^2_{L^2(\Br^{d-1})}
\Norme{\psi_2V_k}_{L^2(\Br^{d+1})}$$
Applying Lemma 9.5 with $p=2$, $\ell=1$, Corollary 9.7 and Proposition
4.1 we obtain
$$\abs{A^2_k}\le C\,h_k\leqno{(9.25)}$$
Using (9.18) and (9.25) we deduce (9.17) which implies (9.11) thus
(9.10). The proof of Theorem 5.2 is complete.\hfill\cqfd

\vpt
{\bf The measure on the boundary}
\vpt
Let us denote by $\dfrac{\pa}{\pa n}$ the normal derivative  at the
boundary $\pa\Omega$. By Corollary 9.6 we see that the sequence
$\displaystyle  \left(\bun_{[0,T]}h_k  \big({\pa w_k   \over\pa
n}\big) _{|\pa\Omega}  \right)$  
is
bounded in $L^2\big(\Br_t\times L^2(\pa\Omega)\big)$. Therefore with the
notations in (5.1) and Proposition~5.1 we have the following Lemma.

\proclaim Lemma 9.8.
 There exist a subsequence $\big(W_{\sigma_1(k)}\big)$ of
$\big(W_{\sigma(k)}\big)$ and a measure $\nu$ on
$T^*\big(\pa\Omega\times\BR_t\big)$ such that for every $a\in{\cal
C}^\infty_0\big(T^*\big(\pa\Omega\times\BR_t\big)\big)$ we have with
$$J_k=\left(a\left(x,t,h_kD_x,h^2_kD_t\right)h_{k}\dfrac 1i\dfrac{\pa
W_k}{\pa n},h_k\dfrac 1i\dfrac{\pa W_t}{\pa
n}\right)_{L^2(\pa\Omega\times\Br_t)}$$
$$\Lim_{k\to +\infty} J_{\sigma_1(k)}=\scal{\nu,a}\leqno{(9.26)}$$

\vpt
{\bf Proof of Theorem 5.3}
\vpt
We begin this proof by considering the case of points inside $T^*M$.

\proclaim Proposition 9.9.
Let  $m_0=\big(x_0,\xi_0,t_0,\tau_0\big)\in T^*M$ and $U_{m_0}$ a
neighborhood of this point in $T^*M$. Then for every $a\in{\cal
C}^\infty_0\big(U_{m_0}\big)$ we have
$$\scal{\mu,H_pa}=0\leqno{(9.27)}$$

\vpt
\noindent {\bf Proof}
\vpt
It is enough to prove (9.27) when
$a(x,t,\xi,\tau)=\Phi(x,\xi)\chi   (t,\tau)$ with $\pi_x\Supp\Phi\subset
V_{x_0}\subset\Omega$. Let $\varphi\in{\cal C}^\infty_0(\Omega)$ be such
that $\varphi=1 $ on $V_{x_0}$. We introduce
$$\eqalign{
A_k=\dfrac
i{h_k}&\Big[\left(\Phi\big(x,h_kD_x\big)\chi\big(t,h^2_kD_t\big)\varphi
h^2_k\big(D_t+P_D\big)\bun_{[0,T]}\sconj w_k,\bun_{[0,T]}\sconj
w_k\right)_{L^2(\Omega\times\Br)}\cr
&-\left(\Phi\big(x,h_kD_x\big)\chi   \big(t,h^2_kD_t\big)
\varphi\bun_{[0,T]}\sconj 
w_k,h^2_k\big(D_t+P_D\big)\bun_{[0,T]}\sconj
w_k\right)_{L^2(\Omega\times\Br)}\Big].}\leqno{(9.28)}$$

We claim that we have
$$\Lim_{k\to +\infty} A_k=0\leqno{(9.29)}.$$
The two terms in $A_k$ are of the same type and will tend both to zero.
Moreover since $\pi_x\Supp\Phi$ and $\varphi$ have compact supports
contained inside $\Omega$ we have in the scalar product
$\big(D_t+P_D\big)\sconj w_k=\big(D_t+P_D\big)w_k=0$. Since
$D_t\big(\bun_{[0,T]}w_k\big)=\dfrac 1iw_k(0)\otimes\delta_{t=0}-\dfrac
1i w_k(T)\otimes\delta_{t=T}$, (9.29) will be proved if we show that
$$\left\{\eqalign{
\Lim B_k&=0\quad\hbox{where}\hc
B_k&=\left(\chi   \big(t,h^2_kD_t\big)\delta_{t=a}h_k
\Phi\big(x,h_kD_x\big)\varphi(x) 
w_k(a),\bun_{[0,T]}w_k\right)_{L^2(\Omega\times\BR)}, a=0,T\hc}\right.
\leqno{(9.30)}$$
Now we have
$$\abs{B_k}\le\Int_0^T\abs{\chi   \big(t,h^2_kD_t\big)\delta_{t=a}}\left(\Int
h_k\abs{\Phi\big(x,h_kD_x\big)\varphi(x)
w_k(a)}\,\abs{w_k(t)}\dx\right)\dt$$
Using the Cauchy-Schwarz inequality in the second integral and the fact
that $w_k(t)=h^{-\sfrac 12}_k\theta\big(h^2_kP_D\big)\widetilde u_k(t)$
we obtain
$$\abs{B_k}\le
C\Int_0^T\abs{\chi   \big(t,h^2_kD_t\big)\delta_{t=a}}\Norme{\widetilde
u_k(a)}_{L^2(\Omega)}\Norme{\widetilde u_k(t)}_{L^2(\Omega)}\dt$$
It follows then from the energy estimate on $\big[0,T\big]$ and Lemma
9.4 with $\ell=2$, $p=1$ that
$$\abs{B_k}\le C\Norme{\widetilde
u^0_k}^2_{L^2(\Omega)}\le\dfrac{C'}k\quad\hbox{by } (3.7)$$
Thus (9.29) is proved.

Let us now compute $A_k$ in another manner. We write
$$\left\{\eqalign{
&A_k=A^1_k+A^2_k\hc
&A^j_k=\dfrac
i{h_k}\left[\left(\Phi\big(x,h_kD_x\big)\chi   \big(t,h^2_kD_t\big)\varphi
h^2_kQ_j\bun_{[0,T]}\sconj w_k,\bun_{[0,T]}\sconj
w_k\right)_{L^2(\Omega\times\BR)}\right.\hc
&\qquad\qquad
-\left.\left(\Phi\big(x,h_kD_x\big)\chi
\big(t,h^2_kD_t\big)\varphi\bun_{[0,T]}\sconj 
w_k, h^2_kQ_k\bun_{[0,T]}\sconj
w_k\right)_{L^2(\Omega\times\Br)}\right]\hc
&\hbox{with}\quad Q_1=D_t,\quad Q_2=P_D\big(x,D_x\big)\hc
}\right.
\leqno{(9.31)}$$
We claim that we have
$$\Lim_{k\to +\infty} A^1_k=0\leqno{(9.32)}$$
Indeed we have
$$A^1_k=-h_k\left(\Phi\big(x,h_kD_x\big)\dfrac{\pa\chi   }{\pa
t}\big(t,h^2_kD_t\big)\varphi\bun_{[0,T]}\sconj w_k,\bun_{[0,T]}\sconj
w_k\right)_{L^2(\Omega\times\Br)}$$
Therefore we have
$$\abs{A^1_k}\le C\,h_k\Int_0^T\Normes{\widetilde\varphi
w_k(t)}^2_{L^2(\Omega)}\dt$$
where $\widetilde\varphi\in{\cal C}^\infty_0(\Omega)$,
$\widetilde\varphi=1$ on $\Supp\varphi$ and (9.32) follows from (4.1).

Now since $P_D$ is self adjoint on $L^2(\Omega)$ we can write
$$A^2_k=\dfrac
i{h_k}\left(\big[\Phi\chi   \varphi,h^2_kP_D\big]\widetilde\varphi
\bun_{[0,T]}\sconj 
w_k,\widetilde\varphi\bun_{[0,T]}\sconj
w_k\right)_{L^2(\Omega\times\BR)}
\leqno{(9.33)}$$
It is easy to see that $h^2_kP_D=\Sum_{j=0}^2
h^{2-j}_kP_j\big(x,h_kD_x\big)$ where $P_j(x,\xi)$ is homogeneous in
$\xi$ of order $j$. %\break\vskip 10mm  %\newpage
 Moreover in the semi classical pseudo-differential
calculus we have $\big[P,Q\big]=\dfrac{h_k}i{\cal
O}p\big(\{p,q\}\big)+h^2_kR$ where~$R$ is $L^2$ bounded. Using the fact
that the sequence $\big(\widetilde\varphi\bun_{[0,T]}\sconj w_k\big)$ is
uniformly bounded in $L^2\big(\Br,L^2\big(\Br^d\big)\big)$ we see easily
that the terms in $A^2_k$ corresponding to $j=0,1$ tend to zero when
$k\to +\infty$. It follows that
$$A^2_k=\left({\cal
O}p\big(\{\Phi\chi   ,p\}\big)\widetilde\varphi\bun_{[0,T]}\sconj
w_k,\bun_{[0,T]}\sconj w_k\right)+o(1)$$
Using (5.1) and Proposition 5.1 we deduce that
$$\Lim_{k\to +\infty} A^2_{\sigma
(k)}=-\Scal{\mu}{H_p\big(\Phi\chi   \big)}
\leqno{(9.34)}$$
It follows from (9.29), (9.31), (9.32) and (9.34) that
$\Scal\mu{H_pa}=0$ if $a=\Phi\chi   $ which implies our
Proposition.\hfill\cqfd

We consider now the case of points $m_0=\big(x_0,t_0,\xi_0,\tau_0\big)$
with $x_0\in\pa\Omega$.

We take a neighborhood $U_{x_0}$ so small that we can perform the
diffeomorphism $F$ described in (9.12). Let $\mu$ and $\nu$ the measures
on $T^*\BR^{d+1}$ and $T^*\big(\pa\Omega\times\Br_t\big)$ defined in
Proposition 5.1 and Lemma 9.7. We shall denote by $\widetilde\mu$ and
$\widetilde\nu$ the measures on $T^*\big(U_0\times\Br_t\big)$ and
$T^*\big(U_0\cap\big\{y_1=0\big\}\times\Br_t\big)$ which are the pull back
of $\mu$ and $\nu$ by the diffeomorphism $\widetilde F$~: $(x,t)\mapsto
\big(F(x),t\big)$.

We first start a Lemma.

\proclaim Lemma 9.10.
Let $a\in{\cal C}^\infty_0\left(T^*\big(U_0\times\Br_t\big)\right)$. We
can find $a_j\in{\cal
C}^\infty_0\left(U_0\times\Br_t\times\Br^{d-1}_{\eta'}\times\Br_\tau\right)$,
$j=0,1$ and $a_2\in{\cal
C}^\infty\left(T^*\big(U_0\times\Br_t\big)\right)$ with compact support
in $\big(y,t,\eta',\tau\big)$ such that with the notations of (9.12)
$$a\big(y,t,\eta,\tau\big)=a_0\big(y,t,\eta',\tau\big)+a_1\big(y,t,\eta',
\tau\big)\eta_1+a_2(y,t,\eta,\tau)\big(\tau+\eta^2_1+r\big(y,\eta'\big)\big)$$
where $r$ is the principal symbol of $R\big(y,D'\big)$.

\vpt
\noindent {\bf Proof}
\vpt
We apply a version of the Malgrange preparation theorem given by Theorem
7.5.4 in   {H{\"o}rmander} [H{\"o}]. With the notations there, for fixed
$m'=\big(y,t,\eta',\tau\big)$ we shall take $t=\eta_1$,
$g_{m'}\big(\eta_1\big)=a(y,t,\eta,\tau)$, $k=2$, $b_1=0$,
$b_0=\tau+r\big(y,\eta'\big)$. According to this theorem we can write
$$a(y,t,\eta,\tau)=q\big(\eta_1,b_0,0,g_{m'}\big)\left(\eta^2_1+r
\big(y,\eta'\big)+\tau\right)+\widetilde 
a_0\big(b_0,0,g_{m'}\big)+\eta_1\widetilde a_1\big(b_0,0,g_{m'}\big)$$
If we multiply both sides by a function
$\varphi=\varphi\big(y,t,\eta',\tau\big)\in{\cal C}^\infty_0$ which is
equal to one on the support in $\big(y,t,\eta',\tau\big)$ of $a$ we
obtain the claim of the Lemma.\hfill\cqfd

In the following Remark we note that we can extend $\widetilde\mu$ to
symbols which are not with compact support in $\eta_1$.

\proclaim Remark 9.11.
Let $q(y,t,\eta,\tau)=\Sum_{j=0}^N q_j\big(y,t,\eta',\tau\big)\eta^j_1$
where $q_j\in{\cal C}^\infty_0\big(\Br^{2d+1}\big)$. Let  $\phi\in{\cal
C}^\infty_0(\Br)$, $\phi\big(\eta_1\big)=1$ if $\abs{\eta_1}\le 1$. Then
$\Scal{\widetilde\mu}{q\phi\left(\dfrac{\eta_1}R\right)}$ does not
depend on $R$ for large $R$. Indeed let $R_2>R_1\gg 1$. Then the symbol
$q\left(\phi\left(\dfrac{\eta_1}{R_1}\right)-\phi\left(\dfrac{\eta_1}
{R_2}\right)\right)$ 
has a support contained in the set $\big\{\abs\tau+\abs{\eta'}\le C,
\abs{\eta_1}\ge R_1\big\}$. Therefore
$$\displaylines{\Supp\widetilde\mu\cap\Supp\left(q\left(\phi\left(
\dfrac{\eta_1}{R_1}\right)-  
\phi\left(\dfrac{\eta_1}{R_2}\right)\right) 
\right)\hfill\cr
\hfill\subset
\left\{\tau+\eta^2_1+r\big(y,\eta'\big)=0,\abs\tau+\abs{\eta'}\le C,
\abs{\eta_1}\ge R\right\}\cr}$$
and the set in the right hand side is empty if $R_1$ is large enough.
\smallskip
We shall set $\Scal{\widetilde\mu}q=\Lim_{R\to+
\infty}\Scal{\widetilde\mu}{q\phi\left(\dfrac{\eta_1}R\right)}$.

  %\newpage

We can now state the analogue of Proposition 9.9 in the case of boundary
points.

\proclaim Proposition 9.12.
With the notations of Lemma 9.10 for any $a\in{\cal
C}^\infty_0\left(T^*\big(U_0\times\Br_t\big)\right)$ we have
$$\Scal{\widetilde\mu}{H_pa}=-\Scal{\widetilde\nu}{{a_1}_{|y_1=0}}$$
   
\vpt
\noindent {\bf Proof}
\vpt
Let us recall that $v_k$ has been defined in (9.13) which satisfies
$$\eqalign{\big(D_t+P_D\big)v_k=&0\hbox{ in } U^+_0=\big\{y\in
U_0~:~y_1>0\big\}\cr
{v_{k}}_{|y_1=0}=&0\cr
v_k=&h^{-\sfrac
12}_k\left(\theta\big(h^2_kP_D\big)w_k\right)\circ F^{-1}}$$
For sake of shortness whe shall set
$$\Lambda_k=\left(y,t,h_kD_y,h^2_kD_t\right),\quad
L^2_+=L^2\big(U_0^+\times\Br_t\big)
\leqno{(9.35)}$$
The Proposition will be a consequence of the following Lemmas.

\proclaim Lemma 9.13.
Let for $j=0,1$, $a_j=a_j\big(y,t,\eta',\tau\big)\in{\cal
C}^\infty_0\big(U_0\times\Br^{d+1}\big)$ and $\varphi\in{\cal
C}^\infty_0\big(U_0\big)$, $\varphi=1$ on $\pi_y\Supp a_j$. 
Then
$$\eqalign{
&\dfrac
i{h_k}\left[\left(\left(a_0\big(\Lambda_k\big)+a_1\big(\Lambda_k\big)
h_kD_1\right)\varphi
h^2_k\big(D_t+P_D\big)\bun_{[0,T]}v_k,\bun_{[0,T]}v_k\right)_{L^2_+}\right.\hc
&\hskip 19mm
\left.-\Int_{U_0^+}\Scal{\left(a_0\big(\Lambda_k\big)+a_1
\big(\Lambda_k\big)h_kD_1\right) 
\varphi\bun_{[0,T]}v_k}{h^2_k\big(D_t+P_D\big)\conj{\bun_{[0,T]}v_k}}
\dy\right]\hc
&=-\dfrac i{h_k}\left(\left[h^2_k\big(D_t+P_D\big),
\left(a_0\big(\Lambda_k\big)+a_1\big(\Lambda_k\big)h_kD_1\right)
\varphi\right]\bun_{[0,T]}v_k,\bun_{[0,T]}v_k\right)_{L^2_+}\hc 
&-\left(a_1\big(0,y',t,h_kD_{y'},h^2_kD_t\big)
\varphi_{|y_1=0}\bun_{[0,T]}\left(h_kD_1{v_k}_{|y_1=0}\right),
\bun_{[0,T]}\left(h_kD_1{v_k}_{|y_1=0}\right)\right)_{L^2(\BR^{d-1}\times
\BR)}\hc}\leqno{(9.36)}$$ 
Here $\scal{\pv,\pv}$ denotes the bracket in $\CD'\big(\BR_t\big)$.

\proclaim Lemma 9.14.
Let $b=b\big(y,t,\eta',\tau\big)\in{\cal
C}^\infty_0\big(U_0\times\Br^{d+1}\big)$ and $\varphi\in{\cal
C}^\infty_0\big(U_0\big)$, $\varphi=1$ on $\pi_y\Supp b$
For $j=0,1$ we set, with the same notations as in Lemma 9.13
\def\bunt{\bun_{[0,T]}}
$$\eqalign{
I^j_k&=\left(h^{-1}_k b\big(\Lambda_k\big)\varphi\big(h_kD_1\big)^j
h^2_k\big(D_t+P_D\big)\bun_{[0,T]}v_k,\bun_{[0,T]}v_k\right)_{L^2_+}\hc
J^j_k&=\Int_{U_0^+}\Scal{h^{-1}_k\,b\big(\Lambda_k\big)\varphi
\big(h_kD_1\big)^j\bunt 
v_k}{h^2_k\big(D_t+P_D\big)\bunt v_k}\dy\hc}$$
Then
$\Lim_{k\to +\infty} I^j_k=\Lim_{k\to +\infty} J^j_k=0$

\proclaim Lemma 9.15.
Let for $j=0,1,2$, $b_j=b_j\big(y,t,\eta',\tau\big)\in{\cal
C}^\infty_0\big(U_0\times\BR^d\big)$ and $\varphi\in {\cal
C}^\infty_0\big(U_0\big)$, $\varphi=1$ on $\pi_y\supp b_j$. Let us set
$$L^j_k=\left(b_j\big(\La_k\big)\varphi\big(h_kD_1\big)^j
\bunt v_k,\bunt  v_k\right)_{L^2_+}$$
Then with $\sigma(k)$ defined in Proposition 5.1 we have for $j=0,1,2$
$$\Lim_{k\to +\infty}
L^j_{\sigma(k)}=\Scal{\widetilde\mu}{b_j\eta^j_1}$$

\vpt
\noindent {\bf Proof of Proposition 9.12}
\vpt
Let $\phi$ be as in Remark 9.11. If $R$ is large enough we have
$H_pa=H_pa\phi\left(\dfrac{\eta_1}{R}\right)$ so by Lemma 9.10
$$\Scal{\widetilde\mu}{H_pa}=\Scal{\widetilde\mu}{H_p\big(a_0+a_1\eta_1\big)
\phi\left(\dfrac{\eta_1}R\right)}+\Scal{\widetilde\mu}{(\tau+p)H_pa_2
\phi\left(\dfrac{\eta_1}R\right)}$$ 
Since Theorem 5.2 implies that
$\Scal{\widetilde\mu}{(\tau+p)H_pa_2\phi\left(\dfrac{\eta_1}R\right)}=0$
we deduce from Remark 9.11 that
$$\Scal{\widetilde\mu}{H_pa}=\Scal{\widetilde\mu}{H_p\big(a_0+a_1\eta_1\big)}$$
Then Proposition 9.12 will be proved if we can show that
$$\Scal{\widetilde\mu}{H_p\big(a_0+a_1\eta_1\big)}=-\Scal{\widetilde\nu}
{{a_1}_{|y_1}=0}\leqno{(9.37)}$$
By Lemma 9.14 the left hand side of (9.36) tends to zero when $k\to
+\infty$.

Now by the semi classical symbolic calculus we can write
$$\dfrac
i{h_k}\left[h^2_k\big(D_t+P_D\big)\left(a_0\big(\Lambda_k\big)+a_1
\big(\Lambda_k\big)h_kD_1\right)\varphi\right]=\Sum_{j=0}^2 
b_j\big(\Lambda_k\big)\varphi_1\big(h_kD_1\big)^j$$
where $b_j\in{\cal C}^\infty_0\big(U_0\times\Br^{d+1}\big)$,
$\varphi_1=1$ on $\Supp\varphi$ and
$\big\{p,a_0+a_1\eta_1\big\}=\Sum_{j=0}^2 b_j \eta^j_1$. So using
(9.36), Lemma 9.15 and Lemma 9.8 we obtain
$0=-\Scal{\widetilde\mu}{H_p\big(a_0+a_1\eta_1\big)}-
\Scal{\widetilde\nu}{{a_1}_{|y_1=0}}$ 
which proves (9.37) and Proposition 9.12.\hfill\cqfd

\vpt
\noindent {\bf Proof of Lemma 9.13}
\vpt
To prove (9.36) we use integration by parts in $D_1$, $D_{y'}$ and
distribution derivative in $D_t$. Only the terms containing $D_1$ give a
boundary contribution. We treat them as follows for $j=0,1$.
$$\displaylines{\left(a_j\big(\Lambda_k\big)\varphi\big(h_kD_1\big)^j\bunt
v_k,h^2_kD^2_1\bunt v_k\right)_{L^2_+}\!\mkern -2mu=\!
\left(h^2_kD^2_1a_j\big(\Lambda_k\big)\varphi\big(h_kD_1\big)^j\bunt
v_k,\bunt v_k\right)_{L^2_+}\hfill\cr
+\dfrac{h_k}i\big(a_j\big(0,y',t,h_kD_{y'},h^2_kD_t\big)\varphi_{|y_1=0}\bunt
\big(h_kD_1\big)^j{v_k}_{|y_1=0},\mkern -2mu
\bunt\big(h_kD_1\big)^j{v_k}_{|y_1=0}  
\big)\mkern -2,5mu _{{L^2(\BR^{d-1}_{y'}\times\BR_t)}}\cr}$$
Here we have used the  boundary condition ${v_k}_{|y_1=0}$.\hfill\cqfd

\vpt
\noindent {\bf Proof of Lemma 9.14}
\vpt
It is enough to consider symbols $b$ of the form
$b=c\big(y,\eta'\big)\chi   (t,\tau)$.

We have $h^2_k\big(D_t+P_D\big)\bunt
v_k=\dfrac{h^2_k}i\big[v_k(0,\cdot)\delta_{t=0}-v_k(T,\cdot)\delta_{t=T}\big]$.
It follows that $I^j_k$ is a sum of two terms of the form
$$\widetilde I\,^j_k=\left(\chi   \big(t,h^2_kD_t\big)\delta_{t=a}h^{\sfrac
12}_kc\big(y,h_kD_{y'}\big)\varphi\big(h_kD_1\big)^j
v_k(a,\cdot),\bunt v_k\right)_{L^2_+}$$
We have with $a=0$ or $T$
$$\abs{\widetilde I\,^j_k}\le
C\Int_{\Br_t}\abs{\chi   \big(t,h^2_kD_t\big)\delta_{t=a}}\Norme{h^{\sfrac
12}_k\big(h_kD_1\big)^j v_k(a,\cdot)}_{L^2(U^+_0)}\Norme{h^{\sfrac
12}_k\bunt v_k(t,\cdot)}_{L^2(U^+_0)}\dt$$
Since for $t\in\big[0,T\big]$,
$\Norme{h^{\sfrac 12}_k\big(h_kD_1\big)^j v_k(t,\cdot)}_{L^2(U^+_0)}
\le C\Norme{\widetilde u\,^0_k}_{L^2(\Omega)}$, $j=0,1$, we deduce from
Lemma 9.5 with $\ell=2$, $p=1$ that $\Lim_{k\to +\infty} \widetilde
I\,^j_k=\Lim_{k\to +\infty} I^j_k=0$.

Let us consider the term $J^j_k$. It is sufficient to consider symbols
$b$ of the form
$$b=\psi(t)c\big(y,h_kD_{y'}\big)\chi   \big(h^2_kD_{\sigma}\big)$$
%
%  est-ce \sigma  ???????????????????????????
%
As before we have $\big(D_t+P_D\big)\bunt v_k=\dfrac
1i\left(v_k(0,\cdot)\delta_{t=0}-v_k(T,\cdot)\delta_{t=T}\right)$ so
$J^j_u$ will be a sum of two terms of the form
$$\widetilde
J^j_k=\Int_{U^+_0}\Scal{h^{-1}_k\psi
c(\cdots)\chi   (\cdots)\varphi\big(h_kD_1\big)^j\bunt
v_k}{\dfrac{h^2_k}i v_k\big(a,\cdot\big)\delta_{t=a}}\dy$$
Writing
$$\widetilde J\,^j_k=\Int_{U_0^+}\Scal
{c\big(y,h_kD_{y'}\big)\varphi\big(h_kD_1\big)^j\bunt v_k}
{\dfrac{h^2}{i}\psi(a)
v_k(a,\cdot)\chi   \big(h^2_kD_t\big)\delta_{t=a}}\dy$$
we see that $\widetilde J\,^j_k$ can be estimated exactly by the same
method as the term $\widetilde I\,^j_k$ above.\hfill\cqfd

\vpt
\noindent {\bf Proof of Lemma 9.15}
\vpt
The proof will be different for each $j$. We shall consider the case
$j=0$ then $j=2$ and finally $j=1$.

\vpt
{\bf Case} $j=0$~: we shall need the following Lemma.
\vpt
\proclaim Lemma 9.16.
Let $\Phi\in{\cal C}^\infty_0(\Br)$, $\phi\big(\eta_1\big)=1$ if
$\abs{\eta_1}\le 1$. 
\smallskip
Let $b=b\big(y,t,\eta',\tau\big)\in{\cal
C}^\infty_0\big(U_0\times\Br^{d+1}\big)$ and $\varphi\in{\cal
C}^\infty_0\big(U_0\big)$, $\varphi=1$ near $\pi_y\Supp b$. Then
$$\Lim_{R\to+\infty}\pv\dsp\mathop{\conj{\Lim}}_{k\to +\infty}
\Norme{\left(I-\Phi\left(\dfrac{h_kD_1}R\right)\right)b\big(\Lambda_k\big)
\varphi\bunt\bun_{y_1>0}v_k}_{L^2(\BR^{d+1})}=0$$

Let us assume this Lemma for a moment and show how it implies the case
$j=0$. We remark first that $b_0\big(\Lambda_k\big)\bunt\bun_{y_1>0}
v_k=\bun_{y_1>0}b_0\big(\La_k\big)\bunt v_k$ which allows us to write
$V_k=\bunt\bun_{y_1>0}v_k$
$$\eqalign{L^0_k
&=\left(b_0\big(\La_k\big)\varphi V_k,\varphi_1 V_k\right)_{L^2_+}\hc
&=\left(\left(I-\Phi\left(\dfrac{h_kD_1}R\right)\right)b_0\big(\La_k\big)
\varphi
V_k,\varphi_1 V_k\right)_{L^2_+}
+\left(\Phi\left(\dfrac{h_kD_1}R\right)b_0\big(\La_k\big)\varphi
V_k,\varphi_1V_k\right)_{L^2_+}\hc
&=A_k+B_k\hc}$$
where $\varphi_1\in{\cal C}^\infty_0\big(U_0\big)$, $\varphi_1=1$ on
$\Supp\varphi$. We have
$$\abs{A_k}\le\Norme{\left(I-\Phi\left(\dfrac{h_kD_1}R\right)\right)b_0
\big(\La_k\big)V_k}_{L^2(\Br^{2d+2})}\Norme{\varphi_1\bunt
\bun_{y_1>0}v_k}_{L^2(\Br^{2d+2})}$$ 
Since by (4.1) $\Norme{\varphi_1\bunt\bun_{y_1>0}v_k}_{L^2}\le C$
uniformly in $k$, Lemma 9.16 shows that \break
$\Lim_{R\to
+\infty}\pv\dsp\mathop{\conj\Lim}_{k\to +\infty} A_k=0$. Now by
Proposition 5.1 we have 
$\Lim_{k\to +\infty}
B_{\sigma(k)}=\Scal{\widetilde\mu}{\Phi\left(\dfrac{\eta_1}R\right)b_0}$
so $\Lim_{R\to +\infty}\Lim_{k\to +\infty}
B_{\sigma(k)}=\Scal{\widetilde\mu}{b_0}$.

\vpt
\noindent {\bf Proof of Lemma 9.16}
\vpt
If $v\in H^1\big(\Br^{d+1}\big)$ we can write
$$\eqalign{
\Int\Norme{\left(I-\Phi\left(\dfrac{h_kD_1}R\right)\right)v}^2_{L^2(\BR^d)}\dt
&=\Int\left(\Int\Bigg|{\dfrac{1-\Phi\left(\dfrac{h_k\eta_1}R\right)}
{h_k\eta_1}}\Bigg|^2\abs{h_k\widehat{D_1v}(\eta,t)}^2\drm\eta\right)\dt\hc 
&\le\dfrac c{r^2}\Norme{h_kD_1v}^2_{L^2(\Br^{d+1})}\hc}$$
Now

$h_kD_1b\big(\La_k\big)\varphi\bunt\bun_{y_1>0}v_k$

\hfill$\eqalign{&=\dfrac{h_k}i\left(\dfrac{\pa b}{\pa
y_1}\big(\La_k\big)\varphi+b\dfrac{\pa\varphi}{\pa
y_1}\right)\bunt\bun_{y_1>0}v_k+b\big(\La_k\big)\bunt\bun_{y_1>0}
\big(h_kD_1v_k\big)\hc 
&=A_k+B_k\hc}$

Here we have used $D_1\big(\bun_{y_1>0}v_k\big)=\bun_{y_1>0}D_1v_k$
since ${v_k}_{|y_1=0}=0$.

We have
$$\Norme{A_k}_{L^2(\Br^{d+1})}\le
C\,h_k\Int_0^T\Norme{\varphi_1w_k(t)}^2_{L^2(\Omega)}\dt$$
where $\varphi_1\in{\cal C}^\infty_0\big(\Br^d\big)$. It follows from
(4.1) that $\Lim_{k\to +\infty} A_k=0$.

Now since $v_k=w_k\circ F^{-1}$ and $w_k=h^{-\sfrac
12}_k\theta\big(h^2_kP_D\big)\widetilde u\,^0_k$ it follows from Lemma
6.3 that
$\Norme{B_k}_{L^2(\BR^{d+1})}\le
C\Int_0^T\Norme{\varphi_1w_k(t)}^2_{L^2(\Omega)}\dt\,\le\,C'$
since
$$\Norme{\left(I-\Phi\left(\dfrac{h_kD_1}R\right)b\big(\La_k\big)\varphi
\bunt\bun_{y_1>0}v_k}_{L^2(\Br^{d+1})}\right)\le\dfrac
C{R_2}\left(\Norme{A_k}_{L^2}+\Norme{B_k}_{L^2}\right)$$
the Lemma follows.\hfill\cqfd

\vpt
{\bf Case} $j=2$ (in Lemma 9.15)
\vpt
Since $\big(h_kD_1\big)^2=h^2_k\big(D_t+P_D\big)-h^2_k
D_t-R_2\big(y,h_kD'\big)-h_kR_1\big(y,h_kD'\big)-h^2_k R_0(y)$ we can
write
$$\left\{\eqalign{
L^2_k&=A_k+B_k\quad\hbox{where}\hc
A_k&=\left(b_2\big(\La_k\big)\varphi h^2_k\big(D_t+P_D\big)\bunt
v_k,\bunt v_k\right)_{L^2_+}\hc
B_k&=\left(c\big(\La_k\big)\varphi\bunt v_k,\bunt v_k\right)_{L^2_+}\hc
C_k&=h_k\left(d\left(h_k,y,t,h_kD_{y'},h^2_kD_t\right)\varphi\bunt
v_k,\bunt v_k\right)_{L^2_+}\hc
\hbox{and } c&=-b_2(\tau+r)\hc}\right.
\leqno{(9.38)}$$
By Lemma 9.14 we have $\Lim_{k\to +\infty} A_k=0$. By Lemma 9.15 for
$j=0$ (case proved above) we have $\Lim_{k\to +\infty}
B_{\sigma(k)}=\Scal{\widetilde\mu}{-b_2(\tau+r)}=
\Scal{\widetilde\mu}{b_2\eta^2_1}$ 
since
$\Scal{\widetilde\mu}{b_2\big(\tau+\eta^2_1+r\big)\phi
\left(\dfrac{\eta_1}R\right)}=0$ 
for all $R$ large enough.

Finally $\abs{C_k}\le M\,h_k\Int_0^T\Norme{w_k(t)}^2_{L^2(\Omega)}\dt\,\le\,
M'\,h_k$ so $\Lim_{k\to +\infty} C_k=0$.

\vpt
{\bf Case} $j=1$
\vpt
We have
$\bun_{y_1>0}\big(h_kD_1\big)v_k=\big(h_kD_1\big)\bun_{y_1>0}v_k$
because ${v_k}_{|y_1=0}=0$. It follows that with
$V_k=\bunt\bun_{y_1>0}v_k$
$$L^1_k=\left(b_1\big(\La_k\big)\varphi\big(h_kD_1\big)
V_k,V_k\right)_{L^2(\BR^{d+1})}$$
Let $\Phi\in{\cal C}^\infty_0(\Br)$, $\Phi\big(\eta_1\big)=1$ if
$\abs{\eta_1}\le 1$. Then we write with $\varphi_1=1$ on $\Supp\varphi$,
$$\eqalign{L^1_k
&=\left(\Phi\left(\dfrac{h_kD_1}R\right)b_1\big(\La_k\big)\varphi\big(h_kD_1
\big)V_k,\varphi_1V_k\right)_{L^2(\BR^{d+1})}\hc
&\hskip 25mm
+\left(\left(I-\Phi\left(\dfrac{h_kD_1}R\right)\right)b_1\big(\La_k
\big)\varphi\big(h_kD_1\big)V_k,\varphi_1 
V_k\right)_{L^2(\BR^{d+1})}\hc
&=A_k+B_k\hc}$$
It is easy to see from Proposition 5.1 that
$$\Lim_{k\to +\infty}
A_{\sigma(k)}=\Scal{\widetilde\mu}{\Phi\left(\dfrac{\eta_1}R\right)b_1\eta_1}$$
Now
$$\abs{B_k}\le
C\Norme{b_1\big(\La_k\big)\varphi\big(h_kD_1\big)V_k}_{L^2(\BR^{d+1})}
\Norme{\left(I-\Phi\left(\dfrac{h_kD_1}R\right)\right)
\varphi_1V_k}_{L^2(\Br^{d+1})}$$ 
As in the proof of Lemma 9.16 we have
$$\Norme{\left(I-\Phi\left(\dfrac{h_kD_1}R\right)\right)\varphi_1
v_k}_{L^2(\Br^{d+1})}\le\dfrac 1R\left(C+h^{\sfrac
12}_k\Norme{\widetilde u\,^0_k}_{L^2(\Omega)}\right)$$
and using the fact that $h_kD_1$ commutes with $\bun_{y_1>0}$ on $v_k$,
since ${v_k}_{y_1=0}=0$ and using Lemma 6.3$(i),(ii)$ we find easily
that
$\Norme{b_1\big(\La_k\big)\varphi\big(h_kD_1\big)V_k}_{L^2(\BR^{d+1})}\le
M$ uniformly in $k$.

\proclaim Lemma 9.17.
With ${\cal G}_d$ and ${\cal G}^k$ introduced in Definition 9.2 we have
$$\widetilde\nu\left({\cal G}_d\cup\left(\Cup_{k=3}^{+\infty} {\cal
G}^k\right)\right)=0$$
 
\vpt
\noindent {\bf Proof}
\vpt
Let us take in Lemma 9.10% Voir s'il ne faut pas mettre proposition 9.12
$a(y,t,\eta,\tau)=b\big(y,t,\eta',\tau\big)\eta_1$. Since
$p=\eta^2_1+r\big(y,\eta'\big)$ we will have
$$\Scal{\widetilde\mu}{\eta_1H_pb-b\dfrac{\pa r}{\pa
y_1}}=-\Scal{\widetilde\nu}{b_{|y_1=0}}
\leqno{(9.39)}$$
Let us make the change of variables
$$(y,t,\eta,\tau)\dsp\mathop{\mapsto}^{\Phi}\big(z=y,
s=t,\zeta=\eta,\sigma=\tau+r(y,\eta)\big)\leqno{(9.40)}$$
It is easy to see that
$$\left\{\eqalign{
H_pb&=X\big(b\circ\Phi^{-1}\big)\circ\Phi\hbox{\quad with}\hc
X&=2\zeta_1\dfrac\pa{\pa z_1}+2\zeta_1\dfrac{\pa r}{\pa
z_1}\dfrac\pa{\pa\sigma}-\dfrac{\pa r}{\pa
z_1}\dfrac\pa{\pa\zeta_1}+H'_r\hbox{\quad where}\hc
H'_r&=\dfrac{\pa r}{\pa\zeta'}\dfrac\pa{\pa z'}-\dfrac{\pa r}{\pa
z'}\dfrac\pa{\pa\zeta'}\hc}\right.
\leqno{(9.41)}$$
It we denote by $\widetilde\mu_1$, $\widetilde\nu_1$ the pull back of
$\widetilde\mu$ and $\widetilde\nu$ by $\Phi$ and $\widetilde
b=b\circ\Phi^{-1}$, the equality (9.39) becomes
$$\Scal{\widetilde\mu_1}{\zeta_1X\widetilde b-\widetilde b\dfrac{\pa
r}{\pa z_1}}=-\Scal{\widetilde\nu_1,}{\widetilde
b\big(0,z',s,\zeta',\sigma\big)}
\leqno{(9.42)}$$
Let us take $\widetilde b$ of the following form
$$\widetilde
b\big(z,s,\zeta',\sigma\big)=b_0\left(\dfrac{z_1}{\sqrt\eps},z',s,\zeta',
\dfrac\sigma\eps\right)\psi\left(\dfrac 
1{\sqrt\eps}\dfrac{\pa r}{\pa z_1}\big(z,\zeta'\big)\right)$$
where $b_0\in{\cal C}^\infty_0$, $b_0\ge 0$, and $\psi\in{\cal
C}^\infty(\Br)$ such that $\psi (t)=1$, $t\in\big(-\infty,0\big]$,
$\psi(t)=0$ for $ t\ge 1$ and $\eps>0$.
  %\newpage
According to (9.42) we can write
\let\wi=\widetilde
$$\left\{\eqalign{
&\Scal{\widetilde\mu_1}{\zeta_1X\widetilde b-\widetilde b\dfrac{\pa
r}{\pa z_1}}=\rondun+\rondtwo\hbox{\quad with}\hc
&\rondun=\Scal{\widetilde\mu_1}{-\dfrac{\pa r}{\pa z_1}b_0\psi}\hc
&\rondtwo=\Scal{\widetilde\mu_1}{f_\eps}\hbox{\quad with}\hc
&f_\eps=2\dfrac{\zeta^2_1}{\sqrt\eps}\dfrac{\pa b_0}{\pa
z_1}\psi+2\dfrac{\zeta^2_1}{\eps}\dfrac{\pa r}{\pa z_1}\dfrac{\pa
b_0}{\pa\sigma}\psi+\zeta_1H'_xb_0\psi+\zeta_1b_0\dfrac1{\sqrt\eps}
\left(X\dfrac{\pa 
r}{\pa z_1}\right)\psi'\hc}\right.
\leqno{(9.43)}$$
According to Theorem 5.2 and (9.40) we have
$\Supp\wi\mu_1\subset\big\{z_1\ge 0\hbox{ and }
\zeta^2_1+\sigma=0\big\}$. Therefore on $\Supp\wi\mu_1\cap\Supp b_0$ we
have $\abs{\zeta_1}^2\le\abs\sigma\le C\eps$.

This implies that $f_\eps\in{\cal C}_0^\infty$ is uniformly bounded in
$\eps\in\big]0,1\big]$.

Moreover the first and the third  term in $f_\eps$ tend to zero
uniformly with $\eps$. The second term can be written on $\Supp\wi\mu_1$
$$-2\dfrac\sigma\eps\dfrac{\pa r}{\pa z_1}\big(z,\zeta'\big)\dfrac{\pa
b_0}{\pa\sigma}\left(\dfrac{z_1}{\sqrt\eps},z',s,\zeta',
\dfrac\sigma\eps\right)\psi(\cdots)$$ 
Since $b_0$ has compact support in $\sigma$, for fixed $\sigma\ne 0$
this term is indentically zero for $\eps$ small enough and it also
vanish when $\sigma=0$.

Finally, since $\Supp\psi'=\big[0,1\big]$ and $\psi'(0)=0$,
$\psi'\left(\dfrac1{\sqrt\eps}\dfrac{\pa r}{\pa
z_1}\big(z,\zeta'\big)\right)$ vanishes if $\eps$ is small enough.
Therefore we can apply the Lebesgue dominated convergence theorem and
conclude that
$$\Lim_{\eps\to 0} \rondtwo=0\leqno{(9.44)}$$
Now let us set $A=\left\{(z,s,\zeta,\sigma)~: z_1=0, \sigma=0,
\dfrac{\pa r}{\pa z_1}\big(z,\zeta'\big)\le 0\right\}$ and write
$$\rondun=\Scal{\wi\mu_1}{-\dfrac{\pa r}{\pa
z_1}b_0\psi\bun_A}+\Scal{\wi\mu_1}{-\dfrac{\pa r}{\pa
z_1}b_0\psi\bun_{A^c}}
\leqno{(9.45)}$$
If we are in $A^c$ we have else $z_1\ne 0$ or $\sigma\ne 0$ or
$\dfrac{\pa r}{\pa z_1}\big(z,\zeta'\big)>0$. In all these case we have
$$\Lim_{\eps\to
0}b_0\left(\dfrac{z_1}{\sqrt\eps},z',s,\zeta',\dfrac\sigma\eps\right)
\psi\left(\dfrac 
1{\sqrt\eps}\dfrac{\pa r}{\pa z_1}\big(z,\zeta'\big)\right)=0$$
By the dominated convergence theorem the second term in the right hand
side of (9.45) tends to zero. It follows that for $\eps$ small enough we
have
$$\rondun=\Scal{\wi\mu_1}{-\dfrac{\pa r}{\pa
z_1}\big(0,z',\zeta'\big)b_0\big(0,z',s,\zeta',0\big)\bun_A}+o(1)$$
Using (9.43), (9.44) we conclude that
$$\Lim_{\eps\to 0}\Scal{\wi\mu_1}{\zeta_1X\wi b-\wi
b\dfrac{\pa r}{\pa z_1}}\ge 0\leqno{(9.46)}$$
On the other hand we have
$$\Scal{\wi\nu_1}{\wi
b\big(0,z',s,\zeta',\sigma\big)}=\Scal{\wi\nu_1}{b_0
\left(0,z',s,\zeta',\dfrac\sigma\eps\right)\psi\left(\dfrac 
1{\sqrt\eps}\dfrac{\pa r}{\pa z_1}\big(0,z',\zeta'\big)\right)}$$
We introduce $B=\big\{\big(z',s,\zeta',\sigma\big)\big\}$~:~ $
\sigma=0$, $\dfrac{\pa r}{\pa z_1}\big(0,z',\zeta'\big)\le0\big\}$
and write as before $1=\bun_B+\bun_{B^c}$. The term corresponding to
$\bun_{B^c}$ tends to zero. It follows that
$$\Scal{\wi\nu_1}{\wi
b_{|z_1=0}}=\Scal{\wi\nu_1}{b_0\big(0,z',s,\zeta',0\big)\bun_B}+o(1)$$
Therefore we have
$$\Lim_{\eps\to 0}\Scal{\wi\nu_1}{\wi
b_{|z_1=0}}=\Scal{\wi\nu_1}{b_0\big(0,z',s,\zeta',0\big)
\bun_{\left\{\sigma\,=\,0,\sfrac{\pa 
r}{\pa z_1}(0,z',\zeta')\le 0\right\}}}\ge 0
\leqno{(9.47)}$$
Using (9.42), (9.46) and (9.47) we conclude that both sides of (9.42)
vanish. Coming back to the coordinates $(y,t,\eta,\tau)$ by (9.40) we
conclude that
$$\Scal{\wi\nu}{b_0\big(0,y',t,\eta',0\big)\bun_{\left\{\tau+r(0,y',\eta')=0,
\sfrac{\pa
r}{\pa y_1}(0,y',\eta')\le0\right\}}}=0$$
for every $b_0\in{\cal C}^\infty_0$, $b_0\ge 0$.

Since ${\cal G}_d\cup\left(\Cup_{k=3}^{+\infty} {\cal
G}^k\right)=\left\{\big(y',t,\eta',\tau\big)~:~\tau+r\big(0,y',
\eta'\big)=0,\dfrac{\pa 
r}{\pa y_1}\big(0,y',\eta'\big)\le 0\right\}$, Lemma 9.17
follows.
\hfill\cqfd

\vpd
\noindent {\bf Proof of the propagation theorem 5.3 (continued)}
\vpd
From now on we follow closely [B], [B-G], [G-L] and we give the details
for the convenience of the reader.

Let us set, with the notations of section 9.1
$${\cal G}^0=T^*M,\quad {\cal G}^1={\cal H}\leqno{(9.48)}$$
We introduce for $k\in\BN$ the following proposition.
$$\left\{\eqalign{
&\hbox{Let $\zeta\in\Sigma_b$. If there exists $T>0$ such that for all
$s\in\big[0,T\big]$ we have}\hc
&\Gamma(s,\zeta)\in\Cup_{j=0}^k{\cal G}^j, \hbox{ then for all
$s_1,s_2$ in $\big[0,T\big]$ we have}\hc
&\pi^{-1}\big(\Gamma\big(s_1,\zeta\big)\big)\cap\Supp\mu\ne\vide\iff
\pi^{-1}\big(\Gamma\big(s_2,\zeta\big)\big)\cap\Supp\mu\ne\vide\hc}\right. 
\leqno{\big(P_k\big)}$$
If $\big(P_k\big)$ is true for all $k\in\BN$ then using Remark 9.4 and a
compacity argument we will obtain the conclusion of Theorem 5.3.

Now $\big(P_k\big)$ is of global nature but as usual using a connexity
argument we can reduce the proof by induction to the following result.
  
\proclaim Proposition 9.18.
   Let $k\ge 1$. Assume $\big(P_{k-1}\big)$ is true. Let $\zeta_0\in{\cal
G}^k$. If there exists $\eps>0$ such that
$\pi^{-1}\big(\Gamma\big(-s,\zeta_0\big)\big)\cap\Supp\mu=\vide$ for all
$s\in\big]0,\eps\big]$ then there exists $\delta>0$ such that
$\pi^{-1}\big(\Gamma\big(s,\zeta_0\big)\big)\cap\Supp\mu=\vide$ for all
$s\in\big[0,\delta\big]$.

Before giving the proof of Proposition 9.18 let us show that
$$ \big(P_0\big) \hbox{ is true}$$
Let $\zeta\in\Sigma_b$ and assume $\Gamma(s,\zeta)\subset{\cal
G}^0=T^*M$ for all $s\in\big[0,T\big[$. Then

\cc{$\zeta=(x,t,\xi,\tau)$ and
$\Gamma(s,\zeta)=\big(x(s),t,\xi(s),\tau\big)$}

 where
$\big(x(s),\xi(s)\big)=\ga\big(s,(x,\xi)\big)$ is the usual
bicharacteristic of $p$ in $T^*M$. Since by Proposition 9.9 we have
$\,^t H_p\mu=0$ in $\CD'\big(T^*M\big)$, the result follows.

\vpd
\noindent {\bf Proof of Proposition 9.18}
\vpd
{\bf Case 1}~: $k=1$
\vpd
Let $\zeta_0=\left(x'_0,t_0,\xi'_0,T_0\right)\in{\cal G}^1={\cal H}$.
Then
$\tau_0+r\left(0,x'_0,\xi'_0\right)=-A<0$.

Let us set
$\xi^0_1=\left(-\left(\tau+r\left(0,x_0',\xi'_0\right)\right)\right)^{\sfrac
12}$. For small $\delta>0$ we set
$$V^{\pm}=\left\{\left(x',t,\xi_1,\xi',\tau\right)~:~\abs{x'-x'_0}<\delta,
\abs{t-t_0}<\delta,\abs{\xi'-\xi'_0}<\delta,\abs{\tau-\tau_0}<\delta,
\abs{\xi_1\mp\xi_1^0}<\delta\right\}$$ 
If $\delta$ is small enough and
$\rho=\left(x_1,x',t,\xi_1,\xi',\tau\right)\in\big[0,\delta\big[\times
\big(V^+\cup V^-\big)$,
we have $\tau+r\left(x_1,x',\xi'\right)\le-\dfrac 12A$. If $p(\rho)=0$
and $\rho\in\big]0,\delta\big[\times V^-$ then $\rho\in T^*M={\cal G}^0$
and $x_1(s)=x_1+2\xi_1s+s^2g(s)$ where $\abs{g(s)}\le C$ and $C$ depends
only on $A$ and $p$. It follows that with $\eps=\dfrac 1{2C}\left(\dfrac
A2\right)^{\sfrac 12}$ we have $x_1(s)>0$ for $s\in\big]-\eps,0\big]$.
This shows that $\Gamma(-s,\rho)\in T^*M$ for $s\in\big]-\eps,0\big]$.
Now by the assumption in proposition 9.18 and continuity one can find
$\beta\in\big]0,\eps\big[$ and $\delta$ small such that
$\pi^{-1}\big(\Gamma(-\beta,\rho)\big)\cap\Supp\mu=\vide$ for all $\rho$
in $\big]0,\delta\big[\times V^-$. It follows then from
$\big(P_0\big)$ that
$$\rho\notin\Supp\mu\hbox{ for all $\rho$ in $\big]0,\delta\big[\times
V^-$}\leqno{(9.49)}$$
Since the hypersurface $\big\{x_1=0\big\}$ is non characteristic for the
vector field $\,^t H_p$ and $\,^t H_p\mu=0$ for $x_1>0$ (Proposition
9.9) the measure $\mu$ has a trace $\mu_{|x_1=0}$ which belongs to
$\CD'\big(V^+\cup V^-\big)$. It follows then that
$$\,^t H_p\mu=2\xi_1\mu_{|x_1=0}\otimes\delta_{x_1=0}\hbox{ in }
\CD'\left(\big]-\delta,\delta\big[\times\big(V^+UV^-\big)\right)
\leqno{(9.50)}$$
Moreover by (9.49) we have
$$\mu_{|x_1=0}=0\hbox{ in } V^-\leqno{(9.51)}$$
Our aim is to show that
$$\mu_{|x_1=0}=0\hbox{ in } V^+\leqno{(9.52)}$$
Indeed let us consider $j$~: $T^*\BR^{d+1}\to T^*\Br^{d+1}$,
$j\big(x,t,\big(\xi_1,\xi'\big),\tau\big)=\big(x,t,\big(-\xi_1,\xi'\big),
  \tau\big)$.
Its follows from Proposition 9.12 that $\Scal{\mu}{H_p(a\circ
j)}=-\Scal{\mu}{H_pa}$. Since $j^{-1}=j$ and $\abs{\det Dj}=1$ we have
$$\big(\,^t H_p\mu\big)\circ j=-\,^t H_p\mu$$
Using (9.50) we see that $-2\xi_1\big(\mu_{|x_1=0}\big)\circ
j\otimes\delta_{x_1=0}=-2\xi_1\mu_{|x_1=0}\otimes\delta_{x_1=0}$ on
$V^+\cup V^-$. Then (9.52) follows from (9.51). Using (9.50), (9.51) and
(9.52) we conclude that $\,^t H_p\mu=0$ on
$\big]-\delta,\delta\big[\times\big(V^+\cup V^-\big)$.

Let us set
$$V=\left\{\left(x',t,\xi_1,\xi',\tau\right)~:~\abs{x'-x'_0}<\delta_1,
\abs{t-t_0}<\delta_1,\abs{\xi'-\xi'_0}<\delta_1,\abs{\tau-\tau_0}<
\delta_1\right\}$$
where 
$0<\delta_1\ll\delta$. Since by theorem 5.2 we have 
$$\Supp\mu\subset\left\{(x,t,\xi,\tau)~:~x_1\ge 0\hbox{ and }
\tau+\xi^2_1+r\big(x,\xi'\big)=0\right\}\hbox{ and}$$
$$\big(\big]-\delta,\delta\big[\times V\big)\cap\left\{x_1\ge 0,
\tau+\xi^2_1+r\big(x,\xi'\big)=0\right\}\subset\big[0,\delta\big[\times
\big(V^+\cup V^-\big)$$
if $\delta_1$ is small enough we deduce that $\,^t H_p\mu=0$ on
$\big]-\delta,\delta\big[\times V$.

Let now $\ga^+(s)=\big(x^+(s),\xi^+(s)\big)$ be the bicharacteristic of
$p$ defined for $\abs s\le\eps$ and starting at the point
$\big(0,x'_0,\xi'_1,\xi'_0\big)$. Since $x^+_1(s)<0$ for
$s\in\big]-\eps,0\big[$ and $\Supp\mu\subset \big\{x_1\ge 0\big\}$ we
have
$\big(x^+(s),t_0,\xi^+(s),\tau_0\big)\notin\Supp\mu$
%\break
for $s\in\big]-\eps,0\big[$ since $\,^t H_p\mu=0$ on
$\big]-\delta,\delta\big[\times V$ one can find $\delta_0>0$ such that
for $s\in\big[0,\delta_0\big]$,
$\big(x^+(s),t_0,\xi^+(s),\tau_0\big)\notin\Supp\mu$ which implies that
$\pi^{-1}\big(\Gamma\big(s,\zeta_0\big)\big)\cap\Supp\mu=\vide$,
$s\in\big[0,\delta_0\big]$, and proves Proposition 9.1 in the case
$k=1$.

\vpt
{\bf Case 2}~:~ $k\ge 2$
\vpt
We shall need several preliminary results.

\proclaim Lemma 9.19.
    Let $\left(0,x'_0,t_0,0,\xi'_0,\tau_0\right)\in T^*\BR^{d+1}$ and
$\delta>0$. We set
$$V\!=\left\{(x,t,\xi,\tau)\in T^*\Br^{d+1},
 0\le
x_1<\delta,\abs{x'-x'_0}<\delta,\abs{t-t_0}<\delta,\abs{\xi'-\xi'_0}<\delta,
\abs{\tau-\tau_0}<\delta\right\}$$
We assume that
$$\eqalign{
(i)\quad& \Supp\mu\cap V\subset\left\{(x,t,\xi,\tau)\in
T^*\BR^{d+1}~:~x_1=\xi_1=0\right\}\hc
(ii)\quad
&\left(0,x'_g\big(s,x',\xi'\big),t,0,\xi'_g\big(s,x',\xi'\big),\tau\right)
\in V \hbox{ for all } s\in I=\big]-s_*,s^*\big[\hc}$$
Then for all $s_1,s_2\in I$
$$\displaylines{\left(0,x'_g\big(s_1,x',\xi'\big),t,0,\xi'_g
\big(s_1,x',\xi'\big),\tau\right)  
\in\Supp\mu\cr
\iff\cr
\left(0,x'_g\big(s_2,x',\xi'\big),t,0,\xi'_2\big(s_2,x',\xi'\big),
\tau\right)\in\Supp\mu\cr}$$

\vpt
\noindent {\bf Proof}
\vpt
By assumption $(i)$ there exists a measure
$\mu_1=\mu_1\big(x',t,\xi',\tau\big)$ on $T^*\BR^d$ such that
$\mu=\mu_1\otimes\delta_{x_1=0}\otimes\delta_{\xi_1=0}$. Moreover we can
extend that definition of $\scal{\mu,a}$ to smooth
$a=a\big(x_1,x',t,\xi_1,\xi',\tau\big)$  which have compact support in
$\big(x_1,x',t,\xi',\tau\big)$ contained in $V$. Indeed if $\chi   \in{\cal
C}^\infty_0(\Br)$, $\chi   \big(\xi_1\big)=1$ for $\abs{\xi_1}\le\dfrac
12$, $\chi   \big(\xi_1\big)=0$ for $\abs{\xi_1}\ge 1$ we can set
$\scal{\mu,a}=\Scal{\mu}{\chi   \big(\xi_1\big)a}$ and this definition does
not depend on $\chi   $. In particular we can take
$a=a\big(x,t,\xi',\tau\big)$. With the notation of Lemma 9.10, we have
$a_1=0$  so it follows from Proposition 9.11 that $\Scal{\mu}{H_pa}=0$.
By the above remark we have
$$0=\Scal{\mu}{H_pa}=\Scal{\mu_1}{H_pa_{|x_1=\xi_1=0}}=
\Scal{\mu_1}{H_{r_0}a_{|x_1=0}}=\Scal{\mu_1}{H_{r_0}\big(a_{|x_1=0}\big)}$$ 
Therefore the Lemma follows.\hfill\cqfd

\proclaim Remark 9.20.
$(i)$ Let $\rho_0=\left(0,x'_0,t_0,0,\xi'_0,\tau_0\right)\in
T^*\BR^{d+1}$. If $\rho_0\notin\Supp\mu$ and
$\tau_0+r\left(0,x'_0,\xi'_0\right)=0$ then
$$\left(x'_0,t_0,\xi'_0,\tau_0\right)\notin\Supp\nu$$
Indeed one can find $\delta>0$ such that
$B\big(\rho_0,\delta\big)\cap\Supp\mu=\vide$.
\smallskip
Let
$a(x,t,\xi,\tau)=b\left(x',t,\xi',\tau\right)\chi   \big(x_1\big)\xi_1$
with support in
\smallskip
$\left\{\abs{x'-x'_0}+\abs{t-t_0}+\abs{\xi'-\xi'_0}+\abs{\tau-\tau_0}
<\delta_1,\abs{x_1}<\delta_1\right\}$. 
Since
$$\Supp\mu\subset
\left\{\abs{\xi_1}^2=\abs{\tau+r\big(x,\xi'\big)}\right\}$$
we will have
our assumption
$$\Supp\mu\cap\Supp a\subset\Supp\mu\cap\supp
a\cap\left\{\abs{\xi_1}\le C\sqrt{\delta_1}\right\}$$
If $\delta_1$ is small enough we will have $\supp
a\cap\left\{\abs{\xi_1}\le C\sqrt{\delta_1}\right\}\subset
B\big(\rho_0,\delta\big)$ so $\supp\mu\cap\supp a=\vide$. With the
notation of Lemma 9.10 we have $a_1=\chi    b$ so, by Proposition 9.12,
$0=\scal{\mu,H_pa}=-\scal{\nu,b}$.
Since $b$ is arbitrary this shows that
$\left(x'_0,t_0,\xi'_0,\tau_0\right)\notin\supp\nu$.
\vpt
$(ii)$ Let $\left(x'_0,t_0,\xi'_0,\tau_0\right)\in{\cal H}$ and
$\rho^\pm_0=\left(0,x'_0,t_0,\pm\xi^0_1,\xi'_0,\tau_0\right)$ where
\smallskip
$\xi^0_1=\left(-\left(\tau_0+r\big(0,x'_0,\xi'_0\big)\right)\right)^{\sfrac
12}$. Then if $\left\{\rho^+_0,\rho^-_0\right\}\cap\supp\mu=\vide$ then
$\left(x'_0,t_0,\xi'_0,\tau_0\right)\notin\supp\nu$. Indeed one can find
$\delta>0$ such that
$$\left(B\big(\rho^+_0,\delta\big)\cup
B\big(\rho^-_0,\delta\big)\right)\cap
\supp\mu=\vide$$
If we take $a=b\chi   \xi_1$ as in $(i)$ we will have
$$\supp\mu\cap\supp a\subset\supp\mu\cap\supp
a\cap\left\{\abs{\xi^2_1-\big(\xi^0_1\big)^2}\le C\delta_1\right\}$$
Now, if $\delta_1\ll\delta^2$,
$\left\{\abs{\xi^2_1-\big(\xi^0_1\big)^2}\le
C\delta_1\right\}\subset\left\{\abs{\xi_1-\xi^0_1}<\delta\right\}\cup
\left\{\abs{\xi_1+\xi^0_1}<\delta\right\}$ 
so
$$\supp a\cap\supp\mu=\vide$$
As in $(i)$ we deduce that $0=\scal{\mu,H_pa}=-\scal{\nu,b}$.

\proclaim Lemma 9.21.
    Let $k\ge 2$ and $\zeta_0=\left(x'_0,t_0,\xi'_0,\tau_0\right)\in{\cal
G}^k$ be given as in Proposition 9.1. Set
$H^{k-2}_{r_0}\left(\dfrac{\pa r}{\pa
x_1}\right)\left(0,x'_0,\xi'_0\right)=A\ne 0$
 then\hfill\break
 $(i)$ one can find $\delta>0$ such that
$\abs{H^{k-2}_{r_0}\left(\dfrac{\pa r}{\pa
x_1}\right)\left(0,x',\xi'\right)}\ge\dfrac{\abs A}2$ in the set
$$V_1=\left\{\big(x',t,\xi',\tau\big)~:~\abs{x'-x_0'}<\delta,
\abs{\xi'-\xi'_0}<\delta,\abs{t-t_0}<\delta,\abs{\tau-\tau_0}<\delta\right\}$$ 
$(ii)$ one can find $\delta'>0$ and $\beta>0$ such that
$$\left\{\eqalign{
&\Gamma\big(s,\widetilde U\big)\cap T^*\pa M\subset V_1\hbox{ for all }
s\in\big[-\beta,0\big],\hc
&\pi^{-1}\big(\Gamma\big(-\beta,\wi U\big)\big)\cap\supp\mu=\vide,\hbox{
where}\hc}\right.$$
$$\displaylines{\wi U=\Big[ \left\{ \big(x',t,\xi',\tau\big)\in T^*\pa M~:
\abs{x'-x'_0}<\delta',
\abs{t-t_0}<\delta',
\abs{\xi'-\xi'_0}<\delta',
\abs{\tau-\tau_0}<\delta'\right\}\hfill\cr
\hfill\cap\left\{\big(x',\xi'\big)~:~
\tau_0+r\big(0,x',\xi'\big)\le 0\right\}\Big]\cr
\cup\Big[ \left\{
(x,t,\xi,\tau)\in T^*M : 0<x_1<\delta',
\abs{x'-x'_0}<\delta',
\abs{\xi'-\xi'_0}<\delta',
\abs{t-t_0}<\delta',
\abs{\tau-\tau_0}<\delta'\right\}\cr
\hfill\cap(\tau+p)^{-1}(0)\Big]\cr}$$
Moreover\hfill\break
{\bf case \rondun~:~ $k$ even, $A>0$}\hfill\break
$\forall\zeta\in{\cal G}^k\cap V_1,\quad \Gamma(s,\zeta)\in{\cal
G}g,\quad\forall s\in\big[-\beta,\beta\big]\moins\{0\},\quad (\forall
s\in\big[-\beta,\beta\big]\hbox{ if } k=2)$.\hfill\break
%\vpd
{\bf case \rondtwo~:~ $k$  even, $A<0$}\hfill\break
$\forall\zeta\in{\cal G}^k\cap V_1,\quad \Gamma(s,\zeta)\in
T^*M,\quad \forall s\in\big[-\beta,\beta\big]\moins\{0\},\quad
(\zeta\in{\cal G}_d\hbox{ if } k=2)$.\hfill\break
{\bf case \rondtrois~:~ $k$ odd, $k\ge 3$, $A<0$}\hfill\break
$\forall\zeta\in{\cal G}^k\cap V_1,\quad \Gamma(s,\zeta)\in
T^*M,\quad \forall s\in\big[-\beta,0\big[,\quad \Gamma(s,\zeta)\in{\cal
G}_g,\quad \forall s\in\big]0,\beta\big]$.\hfill\break
{\bf case \rondq~:~ $k$ odd, $k\ge 3$, $A>0$}\hfill\break
$\forall\zeta\in{\cal G}^k\cap V_1,\quad \Gamma(s,\zeta)\in{\cal
G}_g,\quad \forall s\in\big[-\beta,0\big[,\quad \Gamma(s,\zeta)\in
T^*M,\quad \forall s\in\big]0,\beta\big]$.

\vpt
\noindent {\bf Proof of Lemma 9.21}
\vpt
$(i)$ $H^{k-2}_{r_0}\left(\dfrac{\pa r}{\pa
x_1}\right)\big(0,x',\xi'\big)$ being continuous, the existence of
$\delta$ is clear. Moreover it has the same sign as $A$.
Let us set $e(s,\zeta)=\dfrac{\pa r}{\pa
x_1}\left(0,x'_g(s,\zeta),\xi'_g(s,\zeta)\right)$. Then since
$\zeta\in{\cal G}^k\cap V_1$ we have for small $\abs s$ (see [H], Chap.
24)
$$e(s,\zeta)=\dfrac 1{(k-2)!}H^{k-2}_{r_0}\left(\dfrac{\pa r}{\pa
x_1}\right)\big(0,x',\xi'\big)s^{k-2}+s^{k-1}g(s)$$
where $\abs{g(s)}\le C$, $C$ depending only on $p$. If
$\eps_1=\dfrac{\abs A}{2C}$, $e(s,\zeta)$ has a constant sign on each
interval $\big[-\eps_1,0\big[$ and $\big]0,\eps_1\big]$~; moreover
either $\Gamma(s,\zeta)\in{\cal G}_g$ or $\Gamma(s,\zeta)\in T^*M$ on
each interval which gives the four cases described in the statement of
the Lemma, (see [H], Chap. 24).
\vpt
$(ii)$ Let us prove that
$$\exists \beta_0, \exists\delta'~:~\forall\beta\in\big[0,\beta_0\big[,
\forall x\in\big[-\beta,0\big], \Gamma\big(s,\wi U\big)\cap T^*\pa
M\subset V_1$$
Otherwise one can find sequences $\beta_j\to0$, $\delta'_j\to 0$,
$s_j\in\big]-\beta_j,0\big]$ and $\xi_j\in\wi U$ such that
$\Gamma\big(s_j,\zeta_j\big)\in T^*\pa M$ and
$\Gamma\big(s_j,\zeta_j\big)\notin V_1$. If $\zeta_j\in T^*\pa M$,
$\zeta_j=\left(x'_j,t_j,\xi'_j,\tau_j\right)$ and if $\zeta_j\in T^*M$,
$\zeta_j=\left(x_1^j,x'_j,t_j,\xi^j_1,\xi'_j,\tau_j\right)$. In both
cases $\abs{x'_j-x'_0}<\delta'_j$, $\abs{\xi'_j-\xi'_0}<\delta'_j$,
$\abs{t_j-t_0}<\delta'_j$, $\abs{\tau_j-\tau_0}<\delta'_j$ and if
$\zeta_j\in T^*M$, $0<x^j_1<\delta'_j$. It follows that
$\tau_j+r\left(x^j_1,x'_j,\xi'_j\right)\to\tau_0+r\left(0,x'_0,\xi'_0\right)=0$
which implies that $\xi^j_1\to 0$ (since $\zeta_j\in p^{-1}(0)$ then).
Therefore $\zeta_j\to\zeta_0 $ in $T^*_bM$. Moreover since $s_j\to 0$ we
have $\Gamma\big(s_j,\zeta_j\big)\to\zeta_0$ in $T^*_bM$ so in $T^*\pa
M$ since $\Gamma\big(s_j,\zeta_j\big)\in T^*\pa M$. But $\zeta_0\in V_1$
so $\Gamma\big(s_j,\zeta_j\big)\in V_1$ for large $j$ and we obtain a
contradiction.

Let $\beta=\inf\left(\beta_0,\dfrac{\eps_1}2\right)$ and let us set that
if $\delta'$ is small enough then $\pi^{-1}\big(\Gamma\big(-\beta,\wi
U\big)\big)\cap\supp\mu=\vide$. We know that
$\pi^{-1}\big(\Gamma\big(-\beta,\zeta_0\big)\big)\cap\supp\mu=\vide$.
Let $V\subset T^*\BR^{d+1}$ be such that $V\cap\supp\mu=\vide$ and
$\pi^{-1}\big(\Gamma\big(-\beta,\zeta_0\big)\big)\subset V$. If
$$\forall\delta'>0,\quad \pi^{-1}\big(\Gamma\big(-\beta,\wi
U\big)\big)\cap\supp\mu\ne\vide$$
then there exists $\delta'_j\to 0$, $\zeta_j\in\wi U$ such that
$\rho_j\in\pi^{-1}\big(\Gamma\big(-\beta,\zeta_j\big)\big)$,
$\rho_j\in\supp\mu$.

We keep the notations in the beginning of $(ii)$. Then $x^j_1\to 0$,
$x'_j\to x'_0$, $t_j\to t_0$, $\xi'_j\to\xi'_0$, $\tau_j\to\tau_0$, so
$\zeta_j\to\zeta_0$ which implies that
$\Gamma\big(-\beta,\zeta_j\big)\to\Gamma\big(-\beta,\zeta_0\big)$. Let
us set
$$\Gamma\big(-\beta,\zeta_j\big)=\left\{\eqalign{
&\left(X'_j,T_j,\Xi'_j,\La_j\right)\hbox{ if
}\Gamma\big(-\beta,\zeta_j\big)\in T^*\pa M\hc
&\left(X^j_1,X'_j,T_j,\Xi^j_1,\Xi'_j,\La_j\right)\hbox{ if
}\Gamma(-\beta,\zeta_j\big)\in T^*M\hc}\right.$$
for $j=0$ and $j\ge 1$.

If $\Gamma(-\beta,\zeta_j\big)\in T^*\pa M$ one can find $\Xi^j_1$ such
that $\rho_j=\left(0,X'_j,T_j,\Xi'_1,\Xi'_j,\La_j\right)$ with
$$\big(\Xi^j_1\big)^2+r\left(0,X'_j,\chi   '_j,\Xi'_j\right)+\La_j=0$$
since
$\rho_j\in\supp\mu\cap(\tau+p)^{-1}(0)$.

If $\Gamma\big(-\beta,\zeta_j\big)\in T^*M$ we have the same thing with
$\rho_j=\left(X^j_1,X'_j,T_j,\Xi^j_1,\Xi'_j,\La_j\right)$.

Now in the case where $\Gamma\big(-\beta,\zeta_0\big)\in T^*M$ we have
$X^j_1\to X^0_1>0$, $X'_j\to X'_0$, $T_j\to T_0$, $\Xi^j_1\to\Xi_1^0$,
$\Xi'_j\to\Xi'_0$, $\La_j\to\La_0$. Then $\rho_j\in V$ if $j$ is large
enough which contradicts the fact that $\rho_j\in\supp\mu$ and
$V\cap\supp\mu=\vide$.

In the case where $\Gamma\big(-\beta,\zeta_0\big)\in T^*\pa M$ we have
$X'_j\to X'_0$, $\Xi'_j\to\Xi'_0$, $T_j\to T_0$, $\La_j\to\La_0$ and
$X^j_1\to 0$ (for the indices $j$ such that
$\Gamma\big(-\beta,\zeta_j\big)\in T^*M$).

In the both cases we have $\big(\Xi_1^j\big)^2\to
-r\left(0,X'_0,\Xi'_0\right)-\La_0=\xi^2_1$.

Now
$\left(0,X'_0,T_0,\pm\xi_1,\Xi'_0,\la_0\right)\in\pi^{-1}
\big(\Gamma\big(-\beta,\zeta_0\big)\big)\cap(\tau+p)^{-1}(0)$~; 
therefore \hfill\break $\left(0,X'_0,T_0,\pm\xi_1,\Xi'_0,\La_0\right)\in V$ so
$\rho_j\in V$ for $j$ large enough which is again a contradiction. The
proof of Lemma 9.21 is complete.\hfill\cqfd

\vpt
\noindent {\bf Proof of Proposition 9.18}
\vpt
We are going to consider separately the four cases introduced in Lemma
9.21.

\vpt
{\bf Case \rondun~:}  we have $\Gamma\big(s,\zeta_0\big)\in{\cal G}_g$
for $s\in\big[-\beta,\beta\big]\moins\{0\}$
Therefore
$$\Gamma\big(s,\zeta_0\big)=\left(x'_g\big(s,x'_0,\xi'_0\big),t_0,
\xi'_g\big(s,x'_0,\xi'_0\big),\tau_0\right), 
\quad s\in\big[-\beta,\beta\big]\leqno{(9.53)}$$
Let $U$ be the following set.
$$\displaylines{U=\left\{\big(x',t,\xi',\tau\big)\in T^*\pa
M,\ \abs{x'-x'_0}<\delta',\abs{t-t_0}<\delta,\abs{\xi'-\xi'_0}<\delta',
\abs{\tau-\tau_0}<\delta'\right\}\hfill\cr
\hfill\cup\left\{(x,t,\xi,\tau)\in
T^*M,\ 0<x_1<\delta',\abs{x'-x'_0}<\delta',\abs{t-t_0}<\delta',
  \abs{\xi'-\xi'_0}<\delta',\abs{\tau-\tau_0}<\delta'\right\}\cr}$$
Then $\wi U=U\cap\Sigma_b$ is the set introduced in Lemma 9.21. Moreover
$U$ being an open subset of $T^*_bM$, $\pi^{-1}(U)$ is open in
$T^*\Br^{d+1}$.

By continuity one can find $\eps_0>0$ such that
$\Gamma\big(s,\zeta_0\big)\in U$ for $s$ in$\big[-\eps_0,\eps_0\big]$.
Then one can find $\delta>0$ such that if we set
$$V\!=\left\{(x,t,\xi,\tau)\in T^*\Br^{d+1},0\le
x_1<\delta,\abs{x'-x'_0}<\delta,\abs{t-t_0}<\delta,\abs{\xi'-\xi'_0}
<\delta,\abs{\tau-\tau_0}<\delta\right\}$$ 
then for $s\in\big[-\eps_0,\eps_0\big]$ we have
$$\pi^{-1}\big(\Gamma\big(s,\zeta_0\big)\big)=\left(0,x'_g\big(s,x'_0,\xi'_0
\big),t_0,0,\xi'_g\big(s,x'_0,\xi'_0\big),\tau_0\right)\subset
V\subset\pi^{-1}(U)$$
Assume that we can prove
$$\supp\mu\cap V\subset\left\{(x,t,\xi,\tau)\in
T^*\Br^{d+1}~:~x_1=\xi_1=0\right\}\leqno{(9.54)}$$
then Proposition 9.18 follows immediately from Lemma 9.19.

Let $\rho=(x,t,\xi,\tau)\in\supp\mu\cap\pi^{-1}(U)$. By Theorem 5.2 we
have $\tau+p(x,\xi)=0$ i.e. $\rho\in\Sigma$.

Let
$\zeta=\pi(\rho)\in U\cap\Sigma_b=\wi U$.

If $\big\{\Gamma(-s,\zeta)~:~s\in\big[0,\beta\big]\big\}\cap T^*\pa
M\subset \Cup_{j=0}^{k-1} {\cal G}^j$ then since by Lemma 9.4$(ii)$ we
have
$$\pi^{-1}\big(\Gamma(-\beta,\zeta)\big)\cap\supp\mu=\vide$$
the hypothesis $\big(P_{k-1}\big)$ implies that
$\pi^{-1}(\zeta)\cap\supp\mu=\vide$ which contradicts our assumption
$\rho\in\supp\mu$.

Therefore one can find $s_1\in\big[0,\beta\big]$ such that
$\zeta_1=\Gamma\big(-s_1,\zeta\big)\in T^*\pa M$ but
$\zeta_1\notin\Cup_{j=0}^{k-1} {\cal G}^j$. Since
$$\zeta_1\in\Gamma\big(-s_1,\wi U\big)\cap T^*\pa M\subset V_1$$
 by
Lemma
9.21$(ii)$ we have $\zeta_1\in{\cal G}^k$. Then
$\Gamma\big(s,\zeta_1\big)\in{\cal G}_g$ if
$s\in\big[-\beta+s_1,0\big[\cup\big]0,s_1\big]$.

If $s_1\ne 0$ we have
$\Gamma\big(s_1,\zeta\big)=\Gamma\big(s_1,\Gamma\big(-s_1,\zeta\big)
\big)=\zeta\in{\cal 
G}_g$ and if $s_1=0$ we have $\zeta=\zeta_1\in{\cal G}^k$. In both
cases we have $\zeta=\big(x',t,\xi',\tau\big)$ and
$\rho=\big(0,x',t,0,\xi',\tau\big)$ because
$\tau+r\big(0,x',\xi'\big)=0$. It follows that
$$\supp\mu\cap\pi^{-1}(U)\subset\big\{(x,t,\xi,\tau)~:~x_1=\xi_1=0\big\}$$
as claimed in (9.54).

\vpt
{\bf Case \rondtwo~: } here for $\zeta\in{\cal G}^k\cap V_1$ we have
$\Gamma(s,\zeta)\subset T^*M$ when
$s\in\big[-\beta,\beta\big]\moins\{0\}$

We shall show that
$$\nu=0\hbox{ on } \wi U\cap T^*\pa M\leqno{(9.55)}$$
Since by Lemma 9.17 we have $\nu\left({\cal
G}_d\cup\left(\Cup_{k=3}^{+\infty} {\cal G}^k\right)\right)=0$
it is enough to prove that
$$\supp\nu\cap\wi U\cap\big({\cal H}\cup{\cal
G}_g\big)=\vide$$
Let $\zeta\in\wi U\cap\big({\cal H}\cup{\cal G}_g\big)\cap\supp\nu$.
$$\big\{\Gamma(-s,\zeta)~:~s\in\big[0,\beta\big]\big\}\cap T^*\pa
M\subset \Cup_{j=1}^{k-1}{\cal G}^j$$

since by Lemma 9.21$(ii)$ we have
$\pi^{-1}\big(\Gamma(-\beta,\zeta)\big)\cap\supp\mu=\vide$, by
$\big(P_{k-1}\big)$ we have $\pi^{-1}(\zeta)\cap\supp\mu=\vide$ so, by
Remark 9.20, we have $\zeta\notin\supp\nu$ which is a contradiction. It
follows that we can find $s_1\in\big[0,\beta\big]$ such that
$\zeta_1=\Gamma\big(-s_1,\zeta\big)\in{\cal G}^k$ (${\cal G}_d$ if
$k=2$) (since $\Gamma(-s,\zeta)\in V_1$ by Lemma 9.21). Moreover by
Lemma 9.21, case \rondtwo, we have $\Gamma\big(s,\zeta_1\big)\in T^*M$
if $s\in\big[-\beta+s_1,0\big[\cup\big]0,s_1\big]$.
If $s_1\ne 0$ then  $\zeta=\Gamma\big(s_1,\xi_1\big)\in T^*M$ which
contradicts our assumption. If $s_1=0$ we have $\zeta_1=\zeta\in{\cal
G}^k$ (${\cal G}_d$ if $k=2$) which again is impossible. It follows that
$\wi U\cap\supp\nu\cap\big({\cal H}\cup{\cal G}_g\big)=\vide$ which
proves (9.55).

It follows from Proposition 9.12 that $\!\,^t H_p\mu=0$ on
$\pi^{-1}\big(\wi U\big)$. This implies that the support of $\mu$
propagates along the bicharacteristics of $p$ (with $(t,\tau)=$constant).
Now by assumption
$\pi^{-1}\big(\Gamma\big(-s,\zeta_0\big)\big)\cap\supp\mu=\vide$ and
$\pi^{-1}\big(\Gamma\big(s,\zeta_0\big)\big)\cap(\tau+p)^{-1}(0)=\big(x(s),t_0,
\xi(s),\tau_0\big)$. It follows that for small $s>0$ we have
$\pi^{-1}\big(\Gamma\big(s,\zeta_0\big)\big)\cap\supp\mu=\vide$
which proves Proposition~9.18.

\vpt
{\bf Case \rondtrois~:~ } here $k$ is odd, $k\ge 3$
\vpt
We claim that
$$\supp\mu\cap\pi^{-1}(U)\subset\left\{(x,t,\xi,\tau)\in
T^*\BR^{d+1}~:~x_1=\xi_1=0\right\}\leqno{(9.56)}$$
where $U$ is been defined in case \rondun.

Let $\rho=(x,t,\xi,\tau)\in\pi^{-1}(U)\cap \supp\mu$. Then
$\tau+p(x,\xi)=0$. Let $\zeta=\pi(\rho)\in U\cap\Sigma_b=\wi U$. If
$$\big\{\Gamma(s,\zeta)~:~s\in\big[0,\beta\big]\big\}\subset
\Cup_{j=1}^{k-1}{\cal 
G}^j$$
then $\big(P_{k-1}\big)$ and the fact that
$\pi^{-1}\big(\Gamma(-\beta,\zeta)\big)\cap\supp\mu=\vide $ imply that
$\pi^{-1}(\zeta)\cap\supp\mu=\vide$ which is in contradiction with
$\rho\in\supp\mu$.

Therefore one can find $s_1\in\big[0,\beta\big]$ such that
$\zeta_1=\Gamma\big(-s_1,\zeta\big)\in{\cal G}^k$ (since $\zeta_1\in
V_1$). Since we are in case~\rondtrois\ we have
$$\Gamma\big(s,\zeta_1\big)=\left\{\eqalign{
\left(x'_g(s),t,\xi'_g(s),\tau\right)\quad &s\in\big]0,s_1\big]\hc
\big(x(s),t,\xi(s),\tau\big)\quad
&s\in\big[-\beta+s_1,0\big[\hc}\right.$$
 If $s_1\ne 0$,
$\zeta=\Gamma\Big(s_1,\Gamma\big(-s_1,\zeta\big)\Big)=\Gamma
\big(s_1,\zeta_1\big)\in{\cal 
G}_g$ so $\rho=\big(0,x',t,0,\xi',\tau\big)$.

\noindent If $s_1=0$ then $\zeta=\zeta_1\in{\cal G}^k$ and
$\rho=\big(0,x',0,\xi',\tau\big)$.

This proves (9.56). Therefore we can use Lemma 9.19 and its conclusion
with $V$ such that
$\pi^{-1}\big(\Gamma\big(s,\zeta_0\big)\big)\subset
V\subset\pi^{-1}(U)$ for $s\in\big[-\beta,\beta\big]$.

  %\newpage
Now
$\wi\zeta_0=\left(x'_g\left(-s,x'_0,\xi'_0\right),t_0,\xi'_g
\left(-s,x'_0,\xi'_0\right)\right)\in{\cal 
G}_d$, when $s\in\big]0,\beta\big]$ and $\wi\zeta_0\to\zeta_0$ if $s\to
0$ so $\wi \zeta_0\in\wi U$ if $s$ is small enough, it follows from
Lemma 9.21$(ii)$ that
$\pi^{-1}\big(\Gamma\big(-\beta,\wi\zeta_0\big)\big)\cap\supp\mu=\vide$
since $\big\{\Gamma\big(s,\wi\zeta_0\big)~:
x\in\big[ -\beta,0\big]\big\}\subset{\cal G}_d\cup T^*M$, it follows from
$\big(P_2\big)$, $\pi^{-1}\big(\wi\zeta_0\big)\cap\supp\mu=\vide$. By
Lemma 9.21 we deduce that
$\left(0,x'_g\left(s,x'_0,\xi'_0\right),t_0,0,\xi'_g
\left(s,x'_0,\xi'_0\right),\tau_0\right)\notin\supp\mu$ 
for small $s$, which proves Proposition 9.18 in this case.

\vpt
{\bf Case \rondq}
\vpt
We claim in this case that
$$\supp\nu\cap\wi U=\vide\leqno{(9.57)}$$
As in case \rondtwo\ it is enough to prove that $\supp\nu\cap\wi
U\cap\big({\cal H}\cup{\cal G}_g\big)=\vide$. Let
$\zeta\in\supp\nu\cap\wi U\cap\big({\cal H}\cup{\cal G}_g\big)$. If
$$\big\{\Gamma(-s,\zeta)~:~s\in\big[0,\beta\big]\big\}\cap T^*\pa
M\subset\Cup_{j=1}^{k-1}{\cal G}^j$$
then $\big(P_{k-1}\big)$ and the fact that
$\pi^{-1}\big(\Gamma(-\beta,\zeta)\big)\cap\supp\mu=\vide $ (Lemma
9.21$(ii)$) imply that\break $\pi^{-1}(\zeta)\cap\supp\mu=\vide$. Then
by Remark 9.20 we have $\zeta\notin\supp\nu$ which is a contradiction.
Therefore one can find $s_1\in\big[0,\beta\big]$ such that $\zeta_1=
\Gamma\big(-s_1,\zeta\big)\in{\cal G}^k$ (since $\zeta_1\in V_1)$ and
$$\Gamma\big(s,\zeta_1\big)=\left\{\eqalign{
\big(x(s),t,\xi(s),\tau\big)&\subset T^*M,\quad s\in\big]0,s_1\big],\hc
\left(x'_g(s),t,\xi'_g(s),\tau\right)&\subset T^*M,\quad
s\in\big[-\beta+s_1,0\big[.\hc}\right.$$
If $s_1\ne 0$ then $\zeta=\Gamma\big(s_1,\zeta_1\big)\in T^*M$ which
contradicts our assumption. If $s_1=0$ then $\zeta=\zeta_1\in{\cal G}^k$
which again contradicts the fact that $\zeta\in{\cal H}\cup{\cal G}_g$
since $k$ is odd, $k\ge 3$, in this case. Thus (9.57) is proved. It
follows then, from Proposition 9.12, that $\,^t H_p\mu=0$ on
$\pi^{-1}\big(\wi U\big)$ therefore on a complete neighborhood of
$\pi^{-1}\big(\zeta_0\big)$ in $T^*\Br^{d+1}$ since $\mu\equiv 0$ in
$x_1<0$. Now
$\big(x\big(s,x_0,\xi_0\big),t_0,\xi\big(s,x_0,\xi_0\big),\tau_0\big)$
is contained in $\big\{(x,t,\xi,\tau)~: x_1<0\big\}$ when $\beta+s<0$.
By propagation along the bicharacteristics of $p$ (since
$\,^t H_p\mu=0$) we deduce that
$\Gamma\big(s,\zeta_0\big)=\big(x\big(s,x_0,\xi_0\big),t_0,\xi
\big(s,x_0,\xi_0\big),\tau_0\big)$
does not intersect $\supp\mu$ when $s>0$ is small enough. The proof of
Theorem 5.2 is complete.\hfill\cqfd

\vpt
{\bf 9.3 Proofs of the technical Lemmas}
\vpt
We shall need the following elementary result.

\proclaim Lemma 9.22. 
    Let $P=\Sum_{j,k=1}^d D_ja^{jk}(x)D_k+V$ where $P$ satisfies conditions
(2.3), (2.5) and $V\ge 1$. Then there exists $C\ge 0$ such that for any
$z\in\BC$ such that $\Im z\ne 0$, any $h$ in $\big]0,1\big]$ and any
solution $u\in H^1_0(\Omega)$ of the problem $h^2Pu-zu=f$ with $f\in
L^2(\Omega)$ we have
$$\Norme{h^2Pu}^2+\Sum_{j=1}^n\norme{hD_ju}^2+\Normes{hV^{\sfrac
12}u}^2+\norme u^2\le C\dfrac{\scal{\abs z}^2}{\abs{\Im z}^2}\norme
f^2$$
 where $\norme{\cdot}$ is the norm in $L^2(\Omega)$.

\vpt
\noindent {\bf Proof}
\vpt
Taking the scalar product in $L^2(\Omega)$ of the equation with $u$ we
obtain
$$\big(h^2Pu,u\big)-\Re z\norme u^2-i\Im z\norme
u^2=(f,u)\leqno{(9.58)}$$
Since $\big(h^2Pu,u\big)$ is real, taking the imaginary part of (9.58),
we obtain
$$\norme u\le\dfrac{\norme f}{\abs{\Im z}}\leqno{(9.59)}$$
Now we have $\big(h^2Pu,u\big)\ge
C\left(\Sum_{j=1}^d\Norme{hD_ju}^2+\Normes{hV^{\sfrac 12} u}^2\right)$
so taking the real part in (9.58) and using (9.59) we obtain
$$\Sum_{j=1}^d\Norme{hD_ju}^2+\Normes{hV^{\sfrac 12}
u}^2\le\dfrac{C\abs z}{\abs{\Im z}^2}\norme f^2\leqno{(9.60)}$$
Finally we have $\Norme{h^2Pu}^2\le 2\left(\norme z^2\norme u^2+\norme
f^2\right)$ so using (9.59) and (9.60) we obtain the claim in the
Lemma.\hfill\cqfd

In that follows, we shall make a great use of the so called
{Helffer-Sj{\"o}strand} formula (see [Da]) which will recall now.

Let $\theta\in{\cal C}^\infty_0(\Br)$. We defined an \og almost analytic
extension\fg\ of $\theta$ as follows. Let $\varphi\in{\cal
C}^\infty_0(\Br)$ be such that $\varphi(t)=1$ if $\abs t\le 1$,
$\varphi(t)=0$ if $\abs t\ge 2$.

We set
$$\wi\theta(x,y)=\Sum_{\ell=1}^2\dfrac{\theta^{(\ell)}(x)}{\ell!}
(iy)^\ell\varphi\left(\dfrac 
y{\scal x}\right)\leqno{(9.61)}$$
Then $\wi\theta$ is a ${\cal C}^\infty$ function on $\Br\times\Br$ and
satisfies
$$\left\{\eqalign{
&\abs{\conj\pa\,\wi\theta(x,y)}\le C_N\abs y^2\hbox{ as } \abs y\to
0,\hbox{ where }\hc
&\conj\pa\,\wi\theta(x,y)=\dfrac 12\left(\dfrac{\pa\wi\theta}{\pa x}+i\
\dfrac{\pa\wi\theta}{\pa y}\right)(x,y)\hc}\right.
\leqno{(9.62)}$$
Let $P_D$ be our self adjoint operator defined in (2.1). Then the
Helffer-Sj{\"o}strand formula asserts that
$$\theta\big(h^2P_D\big)=-\dfrac
1\pi\Int_{\Br^2}\conj\pa\,\wi\theta(x,y)\left(z-h^2P_D\right)^{-1}
\dx\dy\leqno{(9.63)}$$ 
where $z =x+iy$.

\vpt
\noindent {\bf Proof of Lemma 6.3}
\vpt
$(i)$ According to (9.63) we have (writing $P$ instead $P_D$)
$$\big[\theta\big(h^2P\big),\chi   \big]=-\dfrac
1\pi\Int_{\Br^2}\conj\pa\,\wi\theta(x,y)\left[\big(z-h^2P\big)^{-1},\chi
\right]\dx\dy\leqno{(9.64)}$$ 
Now $\big(z-h^2P\big)\left[\big(z-h^2P\big)^{-1},\chi   \right]f=\chi   
f+h^2\big[P,\chi   \big]\big(z-h^2P\big)^{-1}f-\chi    f$. Thus
$$\left[\big(z-h^2P\big)^{-1},\chi
\right]=\big(z-h^2P\big)^{-1}h^2\big[P, \chi   \big]\big(z-h^2P\big)^{-1}f$$
Let us set
$\rondun=\Norme{\left[\big(z-h^2P\big)^{-1},\chi   \right]f}_{L^2}$. By
(9.59) we have
$$\rondun\le\dfrac1{\abs{\Im
z}}\Norme{h^2\big[P,\chi   \big]\big(z-h^2P\big)^{-1}f}_{L^2}$$
Now $\big[P,\chi   \big]=\Sum_{j=1}^d b_jD_j+b_0$ where $b_j\in{\cal
C}^\infty_0\big(\conj\Omega\big)$, $j=0,\cdots,d$. It follows that
$$\rondun\le\dfrac{Ch}{\abs{\Im
z}}\left(\Sum_{j=1}^d\Norme{hD_j\big(z-h^2P\big)^{-1}f}_{L^2}+h
\Norme{\big(z-h^2P\big)^{-1}f}_{L^2}\right)$$ 
Using Lemma 9.22 we deduce that
$$\rondun\le\dfrac{C'h\scal{\abs z}}{\abs{\Im z}^2}\norme f_{L^2}$$
It follows from (9.54) that, with $z=x+iy$, we have
$$\Norme{\left[\theta\big(h^2P\big),\chi   \right]f}_{L^2}\le
Ch\Int_{\BR^2}\dfrac{\scal{\abs z}}{\abs{\Im
z}^2}\abs{\conj\pa\,\wi\theta(x,y)}\dx\dy$$
Using formula (9.62) and the fact that $\wi\theta$ has compact support
in $x$ and $y$ we obtain $(i)$.

\vpt
$(ii)$ Again the formula (9.63) we have
$$\Norme{hD_j\theta\big(h^2P\big)f}_{L^2}\le\dfrac
1\pi\Int\abs{\conj\pa\,\wi\theta(x,y)}\,\Norme{h\pa
_j\big(z-h^2P\big)^{-1}f}_{L^2}\dx\dy$$
so using Lemma 9.22 we obtain
$$\Norme{hD_j\theta\big(h^2P\big)f}_{L^2}\le
C\Int\abs{\conj\pa\,\wi\theta(x,y)}\dfrac{\scal{\abs z}}{\abs{\Im
z}}\dx\dy\norme f_{L^2}\le C'\norme f_{L^2}$$
%
%
%\vpt
$(iii)$ We have as above
$$hD_j\left[\theta\big(h^2P\big),\chi   \right]=-\dfrac
1\pi\Int_{\Br^2}\conj\pa\,\wi\theta(x,y)hD_j\left[\big(z-h^2P\big)^{-1},\chi
\right]\dx\dy$$ 
and
$$hD_j\left[\big(z-h^2P\big)^{-1},\chi   \right]
=hD_j\big(z-h^2P\big)^{-1}h^2\big[P,\chi   \big]\big(z-h^2P\big)^{-1}$$
Using again Lemma 9.22 we obtain
$$\eqalign{\Bigg\|&  hD_j\Big[\big(z-h^2P\big)^{-1},\chi  \Big]
f \Bigg\|_{L^2}
\le C\dfrac{\scal{\abs z}}{\abs{\Im
z}}\Norme{h^2\big[P,\chi   \big]\big(z-h^2P\big)^{-1}f}_{L^2}\cr
&\le \dfrac{Ch\scal{\abs z}}{\abs{\Im z}}\left(\Sum_{j=1}^d \Norme{h
D_j\big(z-h^2P\big)^{-1}f}_{L^2}+h\Norme{\big(z-h^2P\big)^{-1}f}_{L^2}
\right)
\le C'h\dfrac{\scal{\abs z}^2}{\abs{\Im z}}\norme f_{L^2}\cr}$$
and we conclude as before.

\vpt
\noindent {\bf Proof of Lemma 8.2}
\vpt

We proceed as above. We have using (9.63) and Lemma 9.22
$$\displaylines{
\Norme{\left[\theta\big(h^2P\big),\chi_0P^{\sfrac
14}\right]v}_{L^2}\hfill\cr
\le \dfrac
1\pi\Int_{\Br^2}\abs{\conj\pa\,\wi\theta(x,y)}\Norme{\big(z-h^2P
\big)^{-1}h^2\big[P,\chi   _0\big]P^{\sfrac 
14}\big(z-h^2P\big)^{-1}v}_{L^2}\dx\dy\hfill\cr
\le\! C\!\!\Int_{\Br^2}\abs{\conj\pa\,\wi\theta(x,y)}\dfrac h{\abs{\Im
z}}\!\left(\Sum_{j=1}^n \Norme{hD_jP^{\sfrac
14}\big(z-h^2P\big)^{-1}v}_{L^2}+\Norme{hP^{\sfrac
14}\big(z-h^2P\big)^{-1}v}_{L^2}\right)\dx\dy\cr
}$$

Now we have with $u=\big(z-h^2P\big)^{-1}v\in D(P)$
$$\eqalign{\Norme{hD_jP^{\sfrac 14}u}_{L^2}
&=h^{-\sfrac 12}\Norme{hD_j\big(h^2P\big)^{\sfrac 14} u}_{L^2}
 \le Ch^{-\sfrac 12}\Norme{\big(h^2P\big)^{\sfrac
12}\big(h^2P\big)^{\sfrac 14} u}_{L^2}\hc
&\le Ch^{-\sfrac 12}\dfrac{\scal{\abs z}}{\abs{\Im z}}\norme
v_{L^2}\hc}$$
by interpolation using Lemma 9.22. It follows that
$$\Norme{\left[\theta\big(h^2P\big),\chi   _0P^{\sfrac
14}\right]v}_{L^2}\le
C\Int_{\Br^2}\abs{\conj\pa\,\wi\theta(x,y)}h^{\sfrac
12}\dfrac{\scal{\abs z}}{\abs{\Im z}^2}\dx\dy\ \norme v_{L^2}\le
Ch^{\sfrac 12}\norme v_{L^2}$$
\hfill\cqfd

  %\newpage
%\partienu{Bibliographie}
%\vskip 20mm
%\leftskip 40bp
\vskip 15pt

\centerline{\soustitre Bibliography}

\vskip 15pt
\def\bi#1{\item {[#1]}}

\bi{BA-D} {\gcap M. Ben Artzi, A. Devinatz},
{\it  Local smoothing and 
convergence properties of Schr{\"o}\-din\-ger type equations},
J. Funct. Anal. 101 (1991), 231-254.

\vpt
\bi{B1} {\gcap N. Burq},
{\it  Mesures semi classiques et mesures de d{\'e}faut},
S{\'e}minaire Bourbaki, Ast{\'e}ris\-tique n$^\circ$245 (1997), p. 167-195.

\vpt
\bi{B2} {\gcap N. Burq},
{\it  Smoothing Effect for Schr{\"o}dinger Boundary Value Problems},
Duke Math. Journal 123 (2004), 403-427.

\vpt
\bi{B3} {\gcap N. Burq},
{\it  Semi-classical estimates for the resolvant in non trapping
geometries},
IMRN n$^\circ$5 (2002), p. 221-241.

\vpt
\bi{B-G} {\gcap N. Burq, P. G{\'e}rard},
{\it  Condition n{\'e}cessaire et suffisante pour la controlabilit{\'e}
exacte des ondes},
CRAS 325 (1997),  749-752.

\vpt
\bi{C-S} {\gcap P. Constantin, J-C. Saut},
{\it  Local smoothing properties of dispersive equations},\hfill\break
J. Amer. Math. Soc. (1988), 413-439.

\vpt
\bi{Da} {\gcap E.B. Davies},
{\it  Spectral theory and differential operators},
Cambridge studies in advanced mathematics, 42, Cambridge Univers. press

\vpt
\bi{D1} {\gcap S. Doi},
{\it  Smoothing effects for Schr{\"o}dinger evolution group on Riemannian
manifolds},
Duke Math. J. 82 (1996), 679-706.

\vpt
\bi{D2} {\gcap S. Doi},
{\it  Smoothing effects for Schr{\"o}dinger evolution equation and global
behavior of geodesic flow},
Math. Ann. 318 (2000), 355-389.

\vpt
\bi{D3} {\gcap S. Doi},
{\it  Smoothing of solutions for Schr{\"o}dinger equations with unbounded
potentiel}.

\vpt
\bi{D4} {\gcap S. Doi},
{\it  Remarks on the Cauchy problem for Schr{\"o}dinger type equations},
Comm. in PDE 21 (1\& 2), (1996), p. 163-178.

\vpt
\bi{G-L} {\gcap P. G{\'e}rard, E. Leichtnam},
{\it  Ergodic properties of eigenfunctions for the Dirichlet problem},
Duke Math. J. 71 n$^\circ$2 (1993), p. 559-607.

\vpt
\bi{H} {\gcap T. Hoschiro},
{\it  Mourre's method and smoothing properties of dispersive
equations},\hfill\break
Comm. Math. Phys. 202 (1999), 255-265.

\vpt
\bi{H{\"o}} {\gcap L. H{\"o}rmander},
{\it  The analysis of Linear Partial Differential Operators I, III},
Springer Verlag, Berlin, Heidelberg, New-York (1985).

\vpt
\bi{K} {\gcap T. Kato},
{\it  On the Cauchy problem for the (generalized) KdV equation},
Stud. Appl. Math. Adv. Math. Suppl. Stud. 18 (1983), 93-128.

  %\newpage
\vpt
\bi{L} {\gcap G. Lebeau},
{\it  \'Equation des ondes amorties. Algebraic and Geometric methods in
math. physics},
Math. Phys. Math. Studies, vol. 19, Kluwer Acad. Publ. Dovdrecht (1996), p.
73-109.

\vpt
\bi{M-S} {\gcap R.B. Melrose, J. Sj{\"o}strand},
{\it  Singularities of boundary value problems I},
Comm. Pure Appl. Math 31 n$^\circ$ 5 (1978), p. 593-617.

\vpt
\bi{Mi} {\gcap L. Miller},
{\it  Refraction of high-frequency waves density by sharp interfaces
and semi classical measures at the boundary},
J. Math. Pures Appl. (9) 79 n$^\circ$ 3 (2000), p. 227-269.

\vpt
\bi{R-Z} {\gcap L. Robbiano, C. Zuily},
{\it  Strichartz estimates for Schr{\"o}dinger equations with variables
coefficients},
Memoires SMF 101-102, 1-206.

\vpt
\bi{Sj} {\gcap P. Sj{\"o}lin},
{\it  Regularity of solutions to the Schr{\"o}dinger equation},
Duke Mat. J. 55 (1987), 699-715.

\vpt
\bi{V} {\gcap L. Vega},
{\it  Schr{\"o}dinger equations~: pointwise convergence to the initial
data},
Proc. Amer. Math. Soc 102 (1988), 874-878.

\vpt
\bi{Y} {\gcap K. Yajima},
{\it  On smoothing property of Schr{\"o}dinger propagator},
Lecture Notes in Math. 1450, Springer Verlag (1990), 20-35.

%\vfill

\bye